\UseRawInputEncoding

\documentclass[]{article}

\usepackage{amssymb,amsmath,amsthm,mathtools}
\usepackage{graphicx}
\usepackage{subcaption}
\usepackage[toc,page]{appendix}
\usepackage{xcolor}
\usepackage{tikz}
\usepackage{dsfont}
\usetikzlibrary{positioning}
\usepackage{tkz-euclide}
\usetikzlibrary{arrows,shapes,positioning,plotmarks}
\usepackage{pgfplots}
\usepackage{grffile}
\pgfplotsset{compat=newest}
\usetikzlibrary{plotmarks}
\usetikzlibrary{arrows.meta,decorations.text}
\usepgfplotslibrary{patchplots}
\usepackage{adjustbox}

\usepackage{subcaption}

\usepackage{authblk}
\usepackage{hyperref}
\usepackage[numbers,sort&compress]{natbib}
\usepackage{setspace}

\usepackage{mathtools}

\usepackage{multirow}

\mathtoolsset{showonlyrefs=true}
\newtheorem{corollary}{Corollary}
\newtheorem{proposition}{Proposition}
\newtheorem{lemma}{Lemma}
\newtheorem{theorem}{Theorem}
\newtheorem{remark}{Remark}
\newtheorem{assumption}{Assumption}
\newtheorem{definition}{Definition}

\title{Heavy-Traffic Universality of Redundancy Systems \\ with Assignment Constraints}

\author[1]{Ellen~Cardinaels\thanks{\href{mailto:e.cardinaels@tue.nl}{e.cardinaels@tue.nl}}}
\author[1]{Sem~Borst}
\author[2]{Johan~S.H.~van~Leeuwaarden}
\affil[1]{Eindhoven University of Technology, The Netherlands}
\affil[2]{Tilburg University, The Netherlands}

\newcommand{\keywords}[1]{\textbf{Keywords:} #1}

\begin{document}
\maketitle

\begin{abstract}
Service systems often face task-server assignment-constraints due to skill-based routing or geographical conditions. Redundancy scheduling responds to this limited flexibility by replicating tasks to specific servers in agreement with these assignment constraints. We gain insight from product-form stationary distributions and weak local stability conditions to establish a state space collapse in heavy traffic. In this limiting regime, the parallel-server system with redundancy scheduling operates as a multi-class single-server system, achieving full resource pooling and exhibiting strong insensitivity to the underlying assignment constraints. In particular, the performance of a fully flexible (unconstrained) system can be matched even with rather strict assignment constraints.
\end{abstract}
\keywords assignment constraints, heterogeneity, parallel-server systems, load balancing, redundancy scheduling, heavy-traffic limit, state space collapse, resource pooling.

\section{Introduction}
Task-server assignment constraints are ubiquitous in a broad range of
everyday service systems. 
In contrast to fully flexible scenarios, where any task can be carried out by any server, assignment of tasks in these systems is restricted to a subset of the servers depending on the underlying features of the tasks and servers. This may potentially degrade the system performance and raises the question how much flexibility in the task-server assignment constraints is needed to achieve performance levels comparable to those in a fully flexible system.

Assignment constraints play a particularly critical role in dynamic matching scenarios.
For example, customers in a ride-sharing network will mainly be matched with drivers in their vicinity, and blood transfusions or organ transplants can only take place whenever there are both available and compatible donors and patients.
Assignment constraints are also prevalent in
customer contact centers. Skills of the various agents, like the spoken language or particular problem-solving skills, narrow the options for forwarding an incoming call.
Skill-based resource management also plays a crucial role in healthcare operations, where there are patients with specific conditions and specialized medical staff members.
Conceptually similar to these skill-based service models are computer systems, such as data center environments and cloud computing platforms.
These systems support a continually evolving variety of applications which involve not just increasing amounts of traffic but also a growing diversity in task types.
Even when tasks or servers are not intrinsically different, some servers tend to be better equipped to perform particular tasks because of data locality and network topology constraints.
This notion is even more pronounced in manufacturing systems where available supply, expected demand and involved costs may be dominant factors to decide which production plants should be able to produce which products.




Motivated by the above observations, we set out to explore how assignment constraints impact the system performance in terms of queue lengths and delays.
We assume that the assignment constraints between task types and servers are described in terms of a general bipartite graph, a so-called compatibility graph, and additionally allow for heterogeneous server speeds.
In particular, we focus on parallel-server systems in a redundancy scheduling setting where replicas are created for each arriving task and then assigned to different servers.
These could either be a subset of $d$~servers selected uniformly at random in power-of-$d$ policies or an arbitrary subset of eligible servers in the context of the assignment constraints as described above.
As soon as the first of these replicas either starts service or completes service, the remaining ones are abandoned (referred to as the `cancel-on-start' (c.o.s.) and `cancel-on-completion' (c.o.c.) versions, respectively).

Dispatching replicas of the same task to several servers increases the chance for one of the replicas to find a short queue and thus start service fast. A well-known application of this paradigm is \textit{multiple listing} where patients in need for an organ transplant register at multiple waiting lists.
It is for instance shown by Zheng et al.~\cite{Zheng2021} that multiple listing for lung transplants reduces the waiting times and improves the probability of finding a suitable donor without affecting the waiting list mortality.
The c.o.s.\ version of redundancy scheduling indeed resembles a Join-the-Smallest-Workload (JSW) policy with partial selection of servers~\cite{Adan2018,ayesta2018unifying}.
The c.o.c.\ version additionally increases the chance for one of the replicas to have a short run time (assuming independent run times on different servers). These relatively small jobs that experience a disproportionally longer run time in computer systems are referred to as \textit{stragglers}, see for instance \cite{Ananthanarayanan2013,Joshi2018}.
On the flip side, the possibly concurrent execution of replicas creates a risk of potential wastage of capacity, depending on run time distributions. Moreover, the cost to cancel jobs in progress can be nonnegligible, urging caution in the implementation of the c.o.c.\ version, see for instance \cite{Shah2016,Joshi2016,Lee2017} where performance trade-offs are investigated.

The launchpad for our analysis is provided by product-form distributions for redundancy systems with arbitrary compatibility graphs as obtained in the seminal papers by Gardner et al.~\cite{Gardner2016queueing} and Ayesta et al.~\cite{ayesta2018unifying}.
While closed-form results for such complex systems are a rare luxury, the expressions in the literature depend on the compatibility graph in a highly intricate fashion, and are unfortunately not particularly transparent.
We therefore adopt a heavy-traffic perspective in order to extract the essential elements and obtain explicit insight into the impact of the assignment constraints on the performance.
The heavy-traffic results reveal a remarkable universality property, and in particular indicate that the performance of a fully flexible system can asymptotically be matched even with rather strict assignment constraints as further discussed below. 

\paragraph{Related literature.}
As mentioned above, product-form distributions for the c.o.c.\ and c.o.s.\ versions of redundancy systems with general compatibility graphs were established in~\cite{Gardner2016queueing} and~\cite{ayesta2018unifying}, respectively.
The latter results extended product-form distributions derived earlier under more restrictive conditions by Visschers et al.~\cite{visschers2012product}.
Related product-form distributions for similar systems with assignment constraints and Assign-to-Longest-Idle-Server-First (ALIS) policies were obtained by Adan et al.~\cite{adan2014skill}.
In fact, all these results turn out to be connected to product-form
distributions for Order-Independent queues~\cite{Krzesinski2011}, the concept of Balanced Fairness~\cite{Bonald2017b,Bonald2017} and token-based central queues~\cite{ayesta2019token}.
A recent overview of queueing models with task-server assignment constraints that yield product-form distributions is provided in~\cite{Gardner2020}.

While such product-form distributions are specified in closed form for general compatibility graphs, the expressions are unwieldy and yield no illuminating formulas for performance metrics like mean queue lengths or delays.
The expressions do not even readily lend themselves for computational purposes, except in specific scenarios where the compatibility graphs satisfy particular structural properties.
For instance, for `nested structures' the mean delays can be computed in an inductive fashion as shown in~\cite{Gardner2017scheduling} and~\cite{Gardner2020}.

Fairly tractable expressions for response time distributions and mean response times have also been derived in~\cite{gardner2017redundancy} and~\cite{Hellemans2018} for the supermarket model with full flexibility, identical servers and power-of-$d$ policies which replicate jobs to each of the subsets of $d$~servers with equal probability. 
Even though this notion of selective randomized replication is conceptually different from replication under assignment constraints, it can mathematically be described as a special instance. 
Since these scenarios are inherently symmetric and constrained to a specific structural family, however, they do not provide generic insight in the performance impact of assignment constraints.

In a different and broader strand of work, stochastic systems have extensively been considered in heavy-traffic regimes to provide greater tractability.  
Indeed, heavy-traffic limits have been established for a wide variety of multi-class parallel-server systems, with a broad range of both task assignment (routing, load balancing) policies and server allocation (scheduling, sequencing) strategies~\cite{Bramson1998,Harrison1998,Harrison1999,Bell2001}.
This body of literature is vast, and an exhaustive review is beyond the scope of this paper.
In particular, the notion of state space collapse has been observed as a common phenomenon in a heavy-traffic regime in, for instance, stochastic processing networks and switched networks \citep{Maguluri2015,Shah2012,Sharifnassab2020}. A stochastic system is said to exhibit state space collapse if its multi-class limiting process has a lower-, often one-dimensional description in contrast to the original pre-limit process.
The notion of state space collapse has often been observed under natural Complete Resource Pooling (commonly abbreviated as CRP) conditions that guarantee the system to behave asymptotically as a single-server queue with pooled resources. These combined concepts are studied, for instance, in \cite{Dai2008,Harrison1999,Stolyar2004,Stolyar2005}.

By comparison, heavy-traffic results for redundancy systems have remained scarce, and the c.o.c.\ version with concurrent execution of tasks in fact seems to move beyond the conventional dynamics considered in the heavy-traffic literature.
The few heavy-traffic results that do exist pertain to the c.o.s.\ version in power-of-$d$ settings, which can be interpreted as special instances with highly symmetric compatibility graphs as noted above.
Atar et al.~\cite{Atar2019replicate,Atar2019} obtain process-level Brownian limits for these models and demonstrate a state space collapse.
Ayesta et al.~\cite{ayesta2018unifying} use the product-form distributions to establish that the total queue length when properly scaled tends to a unit-exponential random variable, but do not consider joint queue length dynamics.
Af\`eche et al.~\cite{Afeche2021} investigate the general setting with assignment constraints, but focus on the expected delay and do not consider multi-dimensional queue-length processes and the associated state space collapse.

As mentioned earlier, the c.o.s.\ version of redundancy scheduling essentially mimics a JSW policy, which in turn closely resembles a Join-the-Shortest-Queue (JSQ) strategy, especially in heavy-traffic conditions.
We will discuss further related literature threads in the realm of JSQ strategies after presenting a detailed model description in Section~\ref{sec:modeldescr}. 

\paragraph{Main contributions.} 
We examine how assignment constraints impact the performance of redundancy systems in terms of queue lengths and delays. We demonstrate that the performance impact tends to be limited, provided the assignment constraints and traffic composition satisfy a mild and natural assumption comparable to the above-mentioned CRP conditions. In particular, if the assignment constraints leave sufficient flexibility for the full service capacity to be used given the load proportions of the various task types, then it is ensured that the assignment constraints create no local capacity bottlenecks and that no subset of the servers can get overloaded as long as the total load is less than the total service capacity.

We establish that when the latter condition holds and traffic is Markovian, the system occupancy exhibits state space collapse in heavy traffic, and asymptotically behaves as in a multi-class single-server FCFS queue.
Informally speaking, the number of tasks of each type remains in strict proportion to the arrival rates in the limit while the total number of tasks, properly scaled, weakly converges to an exponential random variable. Thus the number of tasks of each type has an exponential limiting distribution as well.
By virtue of the distributional form of Little's law, this means that job delay, after scaling, also has an exponential limiting distribution.
Moreover, we extend the above results to scenarios where local capacity bottlenecks could occur, and the queue lengths in a critically loaded subsystem exhibit state space collapse.

In order to prove the above results for the c.o.c.\ mechanism, we start from the product-form distributions for arbitrary assignment constraints as studied in~\cite{Gardner2016queueing}.
We construct a specific enumeration of all the possible task configurations to write the joint~PGF in a convenient form that facilitates a heavy-traffic analysis.
The results for the c.o.s.\ mechanism are established by exploiting the relation between the above product-form expressions and those studied in~\cite{ayesta2018unifying}.

The above results reveal a remarkable universality property, in the sense that the system achieves complete resource pooling and exhibits the same behavior across a broad range of scenarios, as long as no local capacity bottlenecks occur.
In particular, the performance of a fully flexible system can be asymptotically matched, even under quite stringent assignment constraints.
These results translate into several practical implications and guidelines.
First, they indicate that a limited degree of flexibility, when properly designed, is sufficient to achieve full resource pooling. 
Adding greater flexibility provides only limited performance gains, though it may improve robustness if there is uncertainty in the server speeds or load proportions of the various job types.
Second, under a fairly mild condition, the system obeys qualitatively similar scaling laws as if it were fully flexible, so that dimensioning approaches for such scenarios can be adopted without accounting for the assignment relations in great detail.


\paragraph{Organization of the paper.}
In Section~\ref{sec:modeldescr} we present a detailed model description and discuss some broader context and preliminaries, such as the product-form distributions which provide the starting point for our analysis.
The main results are stated in Subsections~\ref{subsec:max_stable} and~\ref{subsec:general}. Numerical results and observations are provided in Subsection~\ref{subsec:discussion}.
Section~\ref{sec:HT} contains the proofs for the c.o.c.\ mechanism. Some details and the proofs for the c.o.s.\ mechanism are deferred to the appendix. We conclude with an outlook for further research in Section~\ref{sec:conclusion}.

\section{Model description and preliminaries}
\label{sec:modeldescr}
\subsection{Model description}
We consider a system with $N$~parallel servers with speeds $\mu_1, \dots, \mu_N$
and several job types that correspond to (non-empty) subsets $S \subseteq \{1, \dots, N\}$ of the servers.
Type-$S$ jobs arrive as a Poisson process of rate~$\lambda_S$, and are replicated to the servers in the subset~$S$.
This set-up fits the premise that job assignment is subject to some constraints, like compatibility relations or data locality issues, as discussed in the introduction.

If we denote by ${\mathcal S} = \{S \in 2^{\{1, \dots, N\}}: \lambda_S > 0\}$ the collection of job types, then the assignment constraints can be represented in terms of a bipartite graph, with nodes for each of the $K = |\mathcal{S}|$ different job types on the one hand and nodes for each of the $N$~servers on the other hand.
A job-type node and a server node are connected by an edge in this compatibility graph whenever a job of this type is replicated to this server~\cite{visschers2012product,adan2014skill,Gardner2016queueing}.

We distinguish between two different versions of redundancy scheduling, referred to as cancel-on-completion (c.o.c.) and cancel-on-start (c.o.s.), respectively.
In the c.o.c.\ version, as soon as the first replica of a particular job finishes service, the remaining replicas are discarded. The sizes of the replicas are independent and exponentially distributed with unit mean.
In the c.o.s.\ version, the redundant replicas are already abandoned as soon as the first replica starts its service. The sizes of the replicas are also exponentially distributed with unit mean, but do not need to be independent.
In either case, each of the servers follows a First-Come First-Served (FCFS) discipline.

For compactness, denote by $\lambda_{\text{tot}} \coloneqq \sum_{S \subseteq \{1, \dots, N\}} \lambda_S$ the total arrival rate and define $\lambda \coloneqq \lambda_{\text{tot}} / N$.
Because of the Poisson splitting property, we can equivalently think of the system as receiving jobs that arrive at rate $\lambda_{\text{tot}}$ and are replicated to the servers in $S \subseteq \{1, \dots, N\}$ with probability $p_S \coloneqq \lambda_S / \lambda_{\text{tot}}$. 

\begin{remark}[Power-of-$d$ policies]\label{remark:pod} \normalfont
We observe that so-called power-of-$d$ policies are subsumed as an important special case of our set-up.
These policies replicate jobs to a randomly selected subset of $d$~servers and have been widely considered in the context of the supermarket model with full flexibility and without any assignment constraints~\cite{gardner2017redundancy}.
From a modeling perspective however, they can be recovered as the special case where the job types correspond to all $K = {N \choose d}$ different subsets of $\{1, \dots, N\}$ of size~$d \leq N$ and $p_S = 1/ K$ for all such~$S$.
Our analysis is entirely cast in terms of the probabilities~$p_S$, and it is immaterial whether these arise from selective replication of jobs (possibly non-uniform and/or randomized), underlying assignment constraints, or a combination of these two factors.
\end{remark}

\begin{remark}[Fully flexible system] \normalfont
We will explore how the system performance is impacted by the heterogeneity of the server speeds $\mu_1, \dots, \mu_N$ and the assignment constraints in terms of the probabilities $(p_S)_{S\in\mathcal{S}}$.
 We examine under what conditions on $\mu_1, \dots, \mu_N$ and $p_S$ performance of a fully flexible system can be approached.
The fully flexible system corresponds to a scenario with homogeneous jobs that can be replicated to all servers, i.e., $\lambda_{\text{tot}} = \lambda_{\{1, \dots, N\}}$ and $p_{\{1, \dots, N\}} = 1$.
The system then behaves as either an M/M/1 FCFS queue with service rate $\mu_{\text{tot}}\coloneqq \sum_{n=1}^{N} \mu_n$ under the c.o.c.\ mechanism or as an M/M/$N$ FCFS queue with heterogeneous server speeds under the c.o.s.\ mechanism. 
\end{remark}

The first-order performance criterion is the stability condition~\cite{adan2014skill,Gardner2016queueing}.

\begin{assumption}[Stability conditions]
Throughout we assume that
\begin{equation}
N \lambda \sum_{S \in {\mathcal T}} p_S <
\sum\limits_{n\in \bigcup\limits_{S\in\mathcal{T}} S}\mu_n
\label{stabcond1}
\end{equation}
for all \textup{(}non-empty\textup{)} ${\mathcal T} \subseteq {\mathcal S}$, or equivalently,
\begin{equation}
N \lambda \sum_{S: S \subseteq U} p_S < \sum\limits_{n\in U} \mu_n
\label{stabcond2}
\end{equation}
for all \textup{(}non-empty\textup{)} $U \subseteq \{1, \dots, N\}$, which have been shown to be necessary and sufficient
conditions for the system to be stable under the c.o.c.\ or c.o.s.\ mechanisms.
\end{assumption}

In particular, taking ${\mathcal T} = {\mathcal S}$
or $U = \{1, \dots, N\}$, we see that $\lambda < \mu$ is a necessary
condition with $\mu \coloneqq \mu_{\text{tot}}/N$ denoting the average service rate across all servers.
Note that the left-hand side of~\eqref{stabcond1} represents the total arrival rate of job types $S \in {\mathcal T}$, while the right-hand side measures the aggregate service rate of the servers that can help in handling jobs of these types.
Likewise, the right-hand side of~\eqref{stabcond2} represents the aggregate service rate of the servers in the set~$U$, while the left-hand side captures the total arrival rate of job types that can be handled by these servers only.

As noted earlier, the c.o.s.\ version amounts to a JSW policy with partial selection of servers, which is covered by the framework of state-dependent assignment strategies by Foss and Chernova~\cite{Foss1998}.
Their stability criteria for systems with `partial accessibility' indicate that the above conditions are in fact necessary and sufficient for arbitrary job size distributions.
Because of the possibly concurrent execution of replicas under the c.o.c.\ policy, the corresponding stability conditions for non-exponential job size distributions are challenging and have remained elusive so far.

In the power-of-$d$ scenario discussed above the stability conditions in~\eqref{stabcond1} and~\eqref{stabcond2} reduce to a set of just $N - d + 1$ inequalities
\begin{equation}
N \lambda \frac{\binom{j}{d}}{\binom{N}{d}} < \sum\limits_{n=1}^{j} \mu_{(n)}
\end{equation}
for $j = d, \dots, N$, with $\mu_{(n)}$ denoting the $n$th smallest service rate among the $N$~servers.
These inequalities reflect that the aggregate service rate of the $j$~slowest servers should exceed the total arrival rate of jobs that are replicated to these servers only.
In the case of identical service rates $\mu_n \equiv \mu$ for all $n = 1, \dots, N$, the inequality for $j = N$ is the most stringent one, yielding the stability condition $\lambda < \mu$, independent of the value of~$d$, which is consistent with the result obtained by Gardner et al.~\cite{gardner2017redundancy} and Anton et al.~\cite{Anton2021}.

\subsection{Further related literature}
The c.o.s.\ version of redundancy scheduling basically emulates a JSW policy, which in turn roughly behaves as a (weighted) JSQ policy, especially in a heavy-traffic regime.
This is reflected by the fact that the above-mentioned framework of state-dependent assignment strategies in \cite{Foss1998} not only covers the JSW policy but also the JSQ policy.
Indeed, the stability conditions for the JSQ policy coincide with~\eqref{stabcond1} and~\eqref{stabcond2} for the JSW policy, see also \cite{Bramson2011} and \cite{Cruise2020}.
The resemblance further manifests itself in the similarity between the process-level limits and state space collapse results for the c.o.s.\ mechanism in power-of-$d$ settings and those for JSQ($d$) policies, i.e., power-of-$d$ versions of JSQ strategies~\cite{Atar2019replicate,Atar2019}. Further heavy-traffic results for JSQ($d$) policies are obtained by Hurtado-Lange and Maguluri~\cite{HurtadoLange2020} in a discrete-time set-up, as well as by Chen and Ye~\cite{Chen2012} and Sloothaak et al.~\cite{Sloothaak2019} for non-uniform sampling variants.

Unlike the situation for redundancy scheduling, the existing literature does contain some results on the performance impact of assignment constraints on JSQ strategies that go beyond power-of-$d$ settings, albeit in a many-server scenario rather than a heavy-traffic regime.
Specifically, \cite{Turner1998,gast2015power,Mukherjee2018,Budhiraja2019} consider JSQ policies in network scenarios where the servers are arranged in a graph structure and each receive arriving jobs which can be forwarded to their neighbors.
In other words, the selection of subsets of servers is not done uniformly at random as in power-of-$d$ policies, but governed by the neighborhood sets in a network graph with the servers as nodes.
These network models can be viewed as a further class of special instances within
the framework that we consider, with identical server speeds, uniform loads across the various job types, and a one-to-one correspondence between job types and servers.

The results in \cite{Turner1998,He2008,gast2015power} pertain to certain fixed-degree graphs, in particular line graphs~\citep{He2008} and ring topologies \citep{Turner1998,gast2015power}.
Their results demonstrate that the performance sensitively depends on the underlying graph topology, and that sampling from fixed neighborhood sets is typically outperformed by re-sampling the same number of alternate servers across the entire system.

In contrast, the results by Mukherjee et al.~\cite{Mukherjee2018} and Budhiraja et al.~\cite{Budhiraja2019} focus on cases where the degrees may grow with the total number of servers.
Their results establish for a many-server regime conditions in terms of the density and topology of the network graph in order for JSQ and JSQ($d$) policies to achieve asymptotically similar performance as in a fully connected graph. From a high level, conceptually related graph conditions for asymptotic optimality were examined using quite different techniques in~\cite{TX13b,TX17} in a dynamic scheduling framework (as opposed to a load balancing context).

The recent papers by Rutten and Mukherjee~\cite{Rutten2020} and Weng et al.~\cite{Weng2020} consider the same general set-up with a bipartite compatibility graph as in the present paper, but pursue a many-server regime and obtain results extending those by \cite{Mukherjee2018} and \cite{Budhiraja2019} to this more general set-up.
Informally speaking, both studies identify conditions in terms of the connectivity properties of the bipartite compatibility graph for similar performance to be achievable as in a fully flexible system.
More specifically, the results in~\cite{Rutten2020} focus on scenarios with identical server speeds and uniform loads across the various job types, and establish process-level limits indicating convergence of the system occupancy under JSQ($d$) policies to that in the supermarket model with full flexibility.
While the results in~\cite{Weng2020} allow for heterogeneous server speeds and arbitrary load distributions, and demonstrate that speed-aware extensions of the JSQ and JIQ strategies achieve vanishing waiting times and minimum expected sojourn times.
Interestingly, the results in \cite{Rutten2020} and \cite{Weng2020} also entail a certain notion of universality, with similar achievable performance as in a fully flexible system under relatively sparse assignment constraints.
However, in the many-server regime this universality property does not manifest itself in terms of a state space collapse of the queue lengths, but rather the fluid-scaled system occupancy showing no queue build-up.

In summary, we map out the existing work in Figure~\ref{fig:diagram_literature} along two dimensions, with JSQ vs redundancy policies as the vertical axis and the many-server vs heavy-traffic regime as the horizontal axis. Literature focusing on the power-of-$d$ setting for the JSQ and redundancy policies, demarcated by a dashed line, is present in all four quadrants and reveals interesting similarities and contrasting features. The results by Eschenfeldt and Gamarnik~\cite{Eschenfeldt2016} and Hellemans and Van Houdt~\cite{Hellemans2021} are positioned at the border of both limiting regimes as they cover systems operating under the heavy-traffic many-server regime for the JSQ($d$) policy and c.o.s. policy, respectively. However, the performance impact of general assignment constraints and heterogeneous server speeds, which falls outside the dashed lines in Figure~\ref{fig:diagram_literature}, has only been pursued for JSQ like strategies in a many-server regime, and has not received any attention so far for redundancy strategies or in a heavy-traffic regime. The only exception is the work by Af\`eche et al.~\cite{Afeche2021}, 
who consider the design of reward-optimal bipartite compatibility graphs. However, their study is focused on expected delay as optimization criterion, and they do not consider convergence of multi-dimensional queue-length processes for arbitrary assignment constraints and associated state space collapse properties.


\def\l{5.75}
\def\h{2.45}
\def\k{0.3}
\def\a{0.45*\l}
\def\aa{2.6}
\def\pod{0.75}
\definecolor{blue}{rgb}{0.00000,0.44700,0.74100}
\definecolor{paars}{rgb}{0.49400,0.18400,0.55600}%
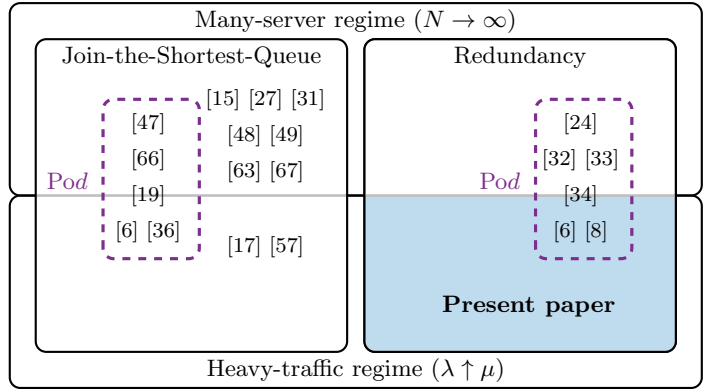
\begin{figure}
\centering
\begin{tikzpicture}[scale=0.8]
\begin{scope}
\clip[rounded corners] (-\a+0.475*\l,-\aa) rectangle (\a+0.475*\l,\aa);
\clip (-\l,0) rectangle (\l,-\k-\h);
\fill[color = blue!100] (-\l,0) rectangle (\l,-2.5*\k-\h);
\end{scope}
\begin{scope}[shift={(0.475*\l,0)}]
\draw[thick,rounded corners] (-\a,-\aa) rectangle (\a,\aa);
\end{scope}
\draw[thick,rounded corners] (-\l,0) rectangle (\l,2.5*\k+\h);
\node  at (0,\h+1.5*\k) {\small Many-server regime ($N\rightarrow\infty$)};
\draw[thick,rounded corners] (-\l,-0) rectangle (\l,-2.5*\k-\h);
\node[]  at (0,-\h-1.5*\k) {\small Heavy-traffic regime ($\lambda\uparrow\mu$)};
\begin{scope}[shift={(-0.475*\l,0)}]
\draw[thick,rounded corners,fill = white,fill opacity=0.75] (-\a,-\aa) rectangle (\a,\aa); 
\node  at (0,\aa-\k) {\small Join-the-Shortest-Queue};
\end{scope}
\begin{scope}[shift={(0.475*\l,0)}]
\draw[thick,rounded corners,fill = white,fill opacity=0.75] (-\a,-\aa) rectangle (\a,\aa);
\node  at (0,\aa-\k) {\small Redundancy};
\end{scope}
\begin{scope}[shift={(-0.6*\l,0)}]
\draw[very thick,dashed, paars,rounded corners] (-0.75,-1.05) rectangle (0.85,1.6);
\node[paars, left] at (-0.85,0.25) {\small Po$d$};
\node at (0,0) {\footnotesize \cite{Eschenfeldt2016}};
\node[anchor=center] at (0,1.2) {\footnotesize \cite{mitzenmacher2001power}};
\node[anchor=center] at (0,0.6) {\footnotesize \cite{vvedenskaya1996queueing}};
\node[anchor=center] at (0,-0.6) {\footnotesize \cite{Atar2019}~\cite{HurtadoLange2020}};
\end{scope}
\begin{scope}[shift={(-1.5,0)}]
\node[draw = none] at (0,-0.85) {\footnotesize \cite{Chen2012}~\cite{Sloothaak2019}};
\node[draw = none] at (0,1.6) {\footnotesize \cite{Budhiraja2019}~\cite{gast2015power}~\cite{He2008}};
\node[draw = none] at (0,1) {\footnotesize \cite{Mukherjee2018}~\cite{Rutten2020}};
\node at(0,0.4) {\footnotesize \cite{Turner1998}~\cite{Weng2020}};
\end{scope}
\begin{scope}[shift={(0.65*\l,0)}]
\draw[very thick,dashed, paars,rounded corners] (-0.75,-1.05) rectangle (0.85,1.6);
\node[paars, left] at (-0.85,0.25) {\small Po$d$};
\node at (0,1.2) {\footnotesize \cite{gardner2017redundancy}};
\node at (0,0.6) {\footnotesize \cite{Hellemans2019}~\cite{Hellemans2018}};
\node at (0,0) {\footnotesize \cite{Hellemans2021}};
\node at (0,-0.6) {\footnotesize \cite{Atar2019}~\cite{ayesta2018unifying}};
\end{scope}
\node at (0.45*\l,-0.25*\h) {\footnotesize \cite{Afeche2021}};
\node[] at (0.5*\l,-0.75*\h) {\small{\textbf{Present paper}}};
\end{tikzpicture}
\caption{A taxonomy of the literature on scheduling policies with assignment constraints in asymptotic regimes. Literature on power-of-$d$ (Po$d$) settings is demarcated by the dashed lines.} \label{fig:diagram_literature}
\end{figure}

\subsection{Product-form distributions}
\label{sec:prodform}
In preparation for our analysis, we review in this subsection the existing product-form distributions for redundancy systems with assignment constraints as described above.

\subsubsection{Redundancy cancel-on-completion}

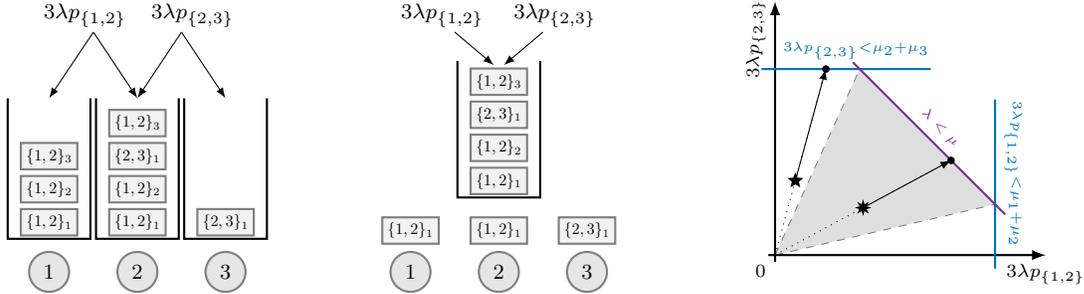
\begin{figure}
\centering
\begin{subfigure}[t]{0.27\textwidth}
  \centering
  \def\h{1.1}
  \def\g{0.8}
\begin{tikzpicture}[>=latex,scale =1]
\tikzstyle{serv} = [circle,draw = gray, thick,fill = gray!20,scale = 0.9]
\tikzstyle{job} = [rectangle,draw = gray, thick,fill = gray!10,scale = 0.6]
\tikzstyle{arr} = [circle,draw = none, text centered]
\begin{scope}[shift={(0,0.1*\h)}]
\node[serv] (s1) at (-2*\h*\g,0) {\small 1};
\node[serv] (s2) at (-0.666*\h*\g,0) {\small 2};
\node[serv] (s3) at (0.666*\h*\g,0) {\small 3};
\end{scope}
\begin{scope}[shift={(-2*\h*\g,0.5*\h)}]
\draw[thick] (-0.5*\h,1.7*\h) -- (-0.5*\h,0) -- (0.5*\h,0) -- (0.5*\h,1.7*\h);
\node[job] at (0,0.2*\h) {$\{1,2\}_1$};
\node[job] at (0,0.6*\h) {$\{1,2\}_2$};
\node[job] at (0,1*\h) {$\{1,2\}_3$};
\end{scope}
\begin{scope}[shift={(-0.666*\h*\g,0.5*\h)}]
\draw[thick] (-0.5*\h,1.7*\h) -- (-0.5*\h,0) -- (0.5*\h,0) -- (0.5*\h,1.7*\h);
\node[job] at (0,0.2*\h) {$\{1,2\}_1$};
\node[job] at (0,0.6*\h) {$\{1,2\}_2$};
\node[job] at (0,1*\h) {$\{2,3\}_1$};
\node[job] at (0,1.4*\h) {$\{1,2\}_3$};
\end{scope}
\begin{scope}[shift={(0.666*\h*\g,0.5*\h)}]
\draw[thick] (-0.5*\h,1.7*\h) -- (-0.5*\h,0) -- (0.5*\h,0) -- (0.5*\h,1.7*\h);
\node[job] at (0,0.2*\h) {$\{2,3\}_1$};
\end{scope}
\draw[->] (-1.33*\h*\g,3*\h) -- (-2*\h*\g,2.2*\h);
\draw[->] (-1.33*\h*\g,3*\h) -- (-0.666*\h*\g,2.2*\h);
\draw[->] (0,3*\h) -- (-0.666*\h*\g,2.2*\h);
\draw[->] (0,3*\h) -- (0.666*\h*\g,2.2*\h);
\node[arr] at (-1.34*\h*\g-0.1,3.2*\h) {\footnotesize $3\lambda p_{\{1,2\}}$};
\node[arr] at (0.15,3.2*\h) {\footnotesize $3\lambda p_{\{2,3\}}$};
\end{tikzpicture}
  \caption{Representation with one queue per server and replicas.}
  \label{fig:example2}
\end{subfigure}~
\begin{subfigure}[t]{0.32\textwidth}
  \centering
\def\h{1.1}
\def\g{0.8}
\begin{tikzpicture}[>=latex,scale =1]
\tikzstyle{serv} = [circle,draw = gray, thick,fill = gray!20,scale = 0.9]
\tikzstyle{job} = [rectangle,draw = gray, thick,fill = gray!10,scale = 0.6]
\tikzstyle{arr} = [circle,draw = none, text centered]
\begin{scope}[shift={(0,0.1*\h)}]
\node[serv] (s1) at (-2*\h*\g,0) {\small 1};
\node[serv] (s2) at (-0.666*\h*\g,0) {\small{2}};
\node[serv] (s3) at (0.666*\h*\g,0) {\small{3}};
\end{scope}
\begin{scope}[shift={(-0.666*\h*\g,1*\h)}]
\draw[thick] (-0.5*\h,1.65*\h) -- (-0.5*\h,0) -- (0.5*\h,0) -- (0.5*\h,1.65*\h);
\node[job] at (0,0.2*\h) {$\{1,2\}_1$};
\node[job] at (0,0.6*\h) {$\{1,2\}_2$};
\node[job] at (0,1*\h) { $\{2,3\}_1$};
\node[job] at (0,1.4*\h) { $\{1,2\}_3$};
\end{scope}
\begin{scope}[shift={(0,0.6*\h)}]
\node[job] at (-2*\h*\g,0) {$\{1,2\}_1$};
\node[job] at (-0.666*\h*\g,0) {$\{1,2\}_1$};
\node[job] at (0.666*\h*\g,0) { $\{2,3\}_1$};
\end{scope}
\begin{scope}[shift={(-0.666*\h*\g,0)}]
\draw[->] (-0.666*\h*\g,3*\h) -- (-0.075*\h*\g,2.6*\h);
\draw[->] (0.666*\h*\g,3*\h) -- (0.075*\h*\g,2.6*\h);
\node[arr] at (-0.85*\h*\g,3.2*\h) {\footnotesize $3\lambda p_{\{1,2\}}$};
\node[arr] at (0.85*\h*\g,3.2*\h) {\footnotesize $3\lambda p_{\{2,3\}}$};
\end{scope}
\end{tikzpicture}
  \caption{Central queue representation. ~~~~ ~~~~~}
  \label{fig:example3}
\end{subfigure}
~
\begin{subfigure}[t]{0.35\textwidth}
  \centering
\def\as{4.5}
\def\e{0.5}
\definecolor{blue}{rgb}{0.00000,0.44700,0.74100}
\definecolor{mycolor1}{rgb}{0.49400,0.18400,0.55600}
\definecolor{orange}{rgb}{0.85000,0.32500,0.09800}
\tikzstyle{starsymbol}=[star, star points=6, star point ratio=2.25,green]
\begin{tikzpicture}[>=latex,scale =0.9]
\fill[color = gray!25] (0,0)--(1.25,2.75)--(3.25,0.75) -- cycle;
\draw[dashed,gray] (0,0) -- (3.25,0.75);
\draw[dashed,gray] (0,0) -- (1.25,2.75);
\draw[->,thick] (-0.1,0) -- (4,0);
\node[below]  at (\as-\e,0) {\scriptsize $3\lambda p_{\{1,2\}}$};
\draw[->,thick] (0,-0.1) -- (0,3.75);
\node[above,rotate =90]  at (0,4-0.8) {\scriptsize $3\lambda p_{\{2,3\}}$};
\node[below left]  at (0,0) {\scriptsize 0};
\draw[thick,blue] (-0.2,2.75) --   (2.3,2.75);
\node[above,blue] at (1.2,2.75) {\tiny $3\lambda p_{\{2,3\}}{<}\mu_2{+}\mu_3$};
\draw[thick,blue] (3.25,2.3) -- (3.25,-0.2);
\node[above,rotate = -90,blue] at (3.25,1.2) {\tiny $3\lambda p_{\{1,2\}}{<}\mu_1{+}\mu_2$};
\draw[thick,mycolor1] (1.15,2.85) -- node[above,sloped] {\tiny  $\lambda<\mu$}(3.4,0.6);
\draw[thin,dotted] (0,0) -- (0.75,2.75);
\draw[thin,->] (0.3,1.1) -- (0.75,2.75);
\node[circle,fill=black,scale=0.3] at (0.75,2.75) {};
\node[star,star points=5,fill=black,star point ratio=2.25,scale=0.3] at (0.3,1.1) {};
\draw[thin,dotted] (0,0) -- (2.6,1.4);
\draw[thin,->] (1.3,0.7) -- (2.6,1.4);
\node[circle,fill=black,scale=0.3] at (2.6,1.4) {};
\node[star,star points=8,fill=black,star point ratio=2.25,scale=0.3] at (1.3,0.7) {};
\end{tikzpicture}
  \caption{The capacity region enclosed by the stability conditions in~\eqref{stabcond1} and~\eqref{stabcond2}. }
  \label{fig:example4}
\end{subfigure}
\caption[]{Representations of the redundancy c.o.c.\ policy. With the assignment constraints as indicated in Figure~(a) with $N=3$ servers, Figures~(a) and (b) give two different representations of the state $\boldsymbol{c} = (\{1,2\},\{1,2\},\{2,3\},\{1,2\})$ of the system operating according to the redundancy c.o.c.\ policy. The notation $\{i,j\}_k$ stands for the $k$th arrival of a type-$\{i,j\}$ job. For the initial arrival rate vectors $3\lambda(p_{\{1,2\}}, p_{\{2,3\}})$ given by \tikz{\node[star,star points=8,fill=black,star point ratio=2.25,scale=0.3] {};} and \tikz{\node[star,star points=5,fill=black,star point ratio=2.25,scale=0.3] {};} in Figure~(c), the critical normalized arrival rate $\lambda^*$ is equal to $\mu$ and $(\mu_2+\mu_3)/(3 p_{\{2,3\}})$, respectively.} \label{fig:example}
\end{figure}

The occupancy of the system at time~$t$ under the redundancy c.o.c.\ policy may be represented in terms of a vector $(c_1, \dots, c_{M(t)})$, with $M(t)$ denoting the total number of jobs, including the ones in service, in the system at time~$t$ and $c_m \in {\mathcal S}$ indicating the type of the $m$th oldest job at that time.
It was shown by Gardner et al.~\cite{Gardner2016queueing} that, if the stability conditions~\eqref{stabcond1} and~\eqref{stabcond2} are satisfied, the stationary distribution of the system occupancy is
\begin{equation}
\label{eq:statdistr}
\pi_{\text{c.o.c.}}(\boldsymbol{c}) = \pi_{\text{c.o.c.}}(c_1, \dots, c_M) = C \prod_{i = 1}^{M} \frac{N\lambda p_{c_i}}{\mu(c_1, \dots, c_i)},
\end{equation}
with $C$ a normalization constant corresponding to the state without any job present and
\begin{equation}
\label{eq:statdistr_mu}
\mu(c_1, \dots, c_i) = \sum\limits_{n \in \bigcup_{m=1}^{i} \{c_m\}} \mu_n.
\end{equation}
Figure~\ref{fig:example} visualizes two different representations of a particular state for a system with $N=3$ servers and assignment constraints depicted at the top of Figure~\ref{fig:example2}. Figure~\ref{fig:example2} shows how the replicas, belonging to the jobs in state~$\boldsymbol{c}$, are stored in separate queues in front of each of the compatible servers. Figure~\ref{fig:example3} visualizes the same state~$\boldsymbol{c}$ from a modeling perspective where all the jobs are stored in a virtual central queue in order of arrival.

\subsubsection{Redundancy cancel-on-start}
\label{subsec:cos}

As soon as a copy starts service at one of its compatible servers, the remaining replicas are instantaneously discarded from the system under the redundancy c.o.s.\ policy.
In a situation where all servers are busy, a server that becomes available selects the longest waiting compatible job.
Whenever a new job arrives and several compatible servers are idle, an assignment rule must be specified.
Within the literature three prominent rules can be distinguished: (i) assign uniformly at random; (ii) assign according to the so-called \textit{assignment condition}~\cite{visschers2012product}; (iii) Assign to the Longest Idle Server (ALIS). In the present paper we will focus on assignment rule~(iii).
The product form of a FCFS-ALIS policy for general assignment constraints was first derived by Adan and Weiss~\cite{adan2014skill}, later it was argued by Adan et al.~\cite{Adan2018} and Ayesta et al.~\cite{ayesta2018unifying} that it is equivalent to the product form under the redundancy c.o.s.\ policy with assignment rule~(iii).\\

The occupancy of the system at time~$t$ under assignment rule~(iii) can be represented as 
$$(c_1,\dots,c_{\tilde{M}(t)};u_1,\dots,u_{L(t)}),$$
 with $\tilde{M}(t)$ and $L(t)$ denoting the total number of {\em waiting} jobs and the number of idle servers at time $t$, respectively.
The type of the $m$th oldest {\em waiting} job is given by $c_m\in\mathcal{S}$ and $u_l \in \{1, \dots, N\}$ represents the $l$th longest waiting idle server.
Note that none of the waiting jobs can be compatible with one of the idle servers and that the state descriptor omits the types of jobs that are in service.
It was shown in~\cite{Gardner2020} that, if the stability conditions~\eqref{stabcond1} and~\eqref{stabcond2} are satisfied, the stationary distribution of the system occupancy is
\begin{equation}
\label{eq:product_form_cos}
\pi_{\text{c.o.s.}}(\boldsymbol{c},\boldsymbol{u}) = \pi_{\text{c.o.s.}}(c_1, \dots, c_M; u_1, \dots, u_L) = C' \prod\limits_{i=1}^{\tilde{M}} \frac{N \lambda p_{c_i}}{\mu(c_1, \dots, c_i)} \prod_{l=1}^{L} \frac{\mu_{u_l}}{\lambda_{\mathcal{C}(u_1, \dots, u_l)}},
\end{equation}
with $C'$ a normalization constant corresponding to the state where all servers are busy and no jobs are waiting, $\mu(c_1, \dots, c_i)$ as in~\eqref{eq:statdistr_mu} and $\lambda_{\mathcal{C}(u_1, \dots, u_l)}$ the arrival rate of jobs that are compatible with (at least one of) the idle servers $(u_1, \dots, u_l)$, i.e.,
\[
\lambda_{\mathcal{C}(u_1,\dots,u_l)} = N \lambda \sum\limits_{S: S \cap \{u_1, \dots, u_l\} \neq \emptyset} p_S.
\]

\section{Main results}
\label{sec:main_results}

We now provide an overview of the main results, relegating the proofs to later sections.
We first introduce some useful notation.
Let $(Q_S)_{S \in {\mathcal S}}$ and $(\tilde{Q}_S)_{S \in {\mathcal S}}$ be random vectors with the stationary distribution of the total number of jobs of each type and the total number of waiting jobs of each type, respectively. Let $(R_n)_{n = 1, \dots, N}$ be a random vector with the joint stationary distribution of the number of replicas assigned to each server. 
Note that not all $R_n$ replicas assigned to server $n$ will necessarily receive (let alone complete) service at that server before being discarded from the system due to the redundancy policy. The random variables with the stationary distribution of the sojourn time and waiting time of an arbitrary type-$S$ job are denoted by $V_S$ and $W_S$, respectively.

For compactness, with minor abuse of notation, define $p_{{\mathcal T}} \coloneqq \sum_{S \in {\mathcal T}} p_S$ as the fraction of jobs that belong to the subset of job types ${\mathcal T} \subseteq {\mathcal S}$.
Also, define 
\[
\mu_{{\mathcal T}}  \coloneqq \sum_{n \in \bigcup_{S \in \mathcal{T}} S} \mu_n
\] 
as the aggregate service rate of the servers that are compatible with the subset of job types ${\mathcal T} \subseteq {\mathcal S}$.

\begin{definition}[Critical subset and critical arrival rate]\label{def:critical_subset_arrrate}
The ratio $\rho_{{\mathcal T}} \coloneqq \frac{p_{{\mathcal T}}}{\mu_{{\mathcal T}}}$ is called the capacity ratio of the subset of job types ${\mathcal T} \subseteq {\mathcal S}$.
The subset ${\mathcal T}^* \coloneqq \arg\max_{{\mathcal T} \subseteq {\mathcal S}}
\rho_{{\mathcal T}}$ with the maximum capacity ratio is called the critical subset, and $\lambda^* = \frac{1}{N \rho_{{\mathcal T}^*}} = \frac{\mu_{{\mathcal T}^*}}{N p_{{\mathcal T}^*}}$ is called the critical (normalized) arrival rate.
\end{definition} 
Note that the stability conditions~\eqref{stabcond1} and~~\eqref{stabcond2} may be equivalently written as $\lambda < \lambda^*$.

\subsection{Heavy-traffic regime under the local stability conditions}\label{subsec:max_stable}
In this subsection we assume that the probabilities $(p_S)_{S\in\mathcal{S}}$ and the server speeds $\mu_1, \dots, \mu_N$ satisfy the so-called local stability conditions.

\begin{assumption}[Local stability conditions]
For all \textup{(}non-empty\textup{)} ${\mathcal T} \subsetneq {\mathcal S}$:
\begin{equation}
p_{{\mathcal T}} < \frac{1}{N \mu} \mu_{{\mathcal T}}
\label{stabcond3}
\end{equation}
 with $\mu = \frac{1}{N}\mu_{\text{tot}}$ the average service rate as defined earlier, or equivalently, for all \textup{(}non-empty\textup{)} $U \subsetneq \{1, \dots, N\}$:
\begin{equation}
\sum_{S: S \subseteq U} p_S < \frac{1}{N \mu} \sum\limits_{n \in U} \mu_n.
\label{stabcond4}
\end{equation}
\end{assumption}

Note that we restrict to strict subsets ${\mathcal T} \subsetneq {\mathcal S}$ and $U \subsetneq \{1, \dots, N\}$, as \eqref{stabcond3} and~\eqref{stabcond4} hold with equality in case ${\mathcal T} = {\mathcal S}$ and $U = \{1, \dots, N\}$, respectively, by definition of the average service rate~$\mu$.
Then $\rho_{{\mathcal T}} < \frac{1}{N \mu} = \rho_{{\mathcal S}}$ for all
${\mathcal T} \subsetneq {\mathcal S}$, so that the critical `subset' is ${\mathcal T}^* = {\mathcal S}$ and the critical normalized arrival rate is $\lambda^* = \mu$.

Also, the inequalities in~\eqref{stabcond3} and~\eqref{stabcond4} imply that the stability conditions in~\eqref{stabcond1} and~\eqref{stabcond2} are satisfied for any $\lambda < \mu$, i.e., as long as the total arrival rate is strictly less than the aggregate service rate.
More specifically, the strict inequalities in~\eqref{stabcond3} ensure that as $\lambda\uparrow\mu$, for all (non-empty) $\mathcal{T} \subsetneq \mathcal{S}$, the total arrival rate of the job types in~$\mathcal{T}$ remains bounded away from the aggregate service rate of the servers that can help in handling jobs of these types. Or more formally, it can be deduced that there exists an $\epsilon>0$ such that for all $\mathcal{T}\subsetneq\mathcal{S}$ we have that $N\mu p_{\mathcal{T}} +\epsilon < \mu_{\mathcal{T}}$, and hence also $N\lambda p_{\mathcal{T}} +\epsilon < \mu_{\mathcal{T}}$ for all $\lambda\le\mu$. 
While the above expression is closely related to the form of the inequalities in~\eqref{stabcond3} and~\eqref{stabcond4}, it can be rewritten in a form that is better suited for application in the proofs later on: There exists an $\epsilon>0$ such that for all $\mathcal{T}\subsetneq\mathcal{S}$ we have that $\frac{N\mu p_{\mathcal{T}}}{\mu_{\mathcal{T}}}<1-\epsilon$. Similarly, the strict inequalities in~\eqref{stabcond4} guarantee that as $\lambda\uparrow\mu$, for all (non-empty) $U \subset\{1, \dots, N\}$, the aggregate service rate of the servers in the set~$U$ remains bounded away from the total arrival rate of job types that can be handled by these servers only.
Thus, the inequalities in~\eqref{stabcond3} and~\eqref{stabcond4} ensure that there are no local capacity bottlenecks and that as $\lambda\uparrow\mu$, only the inequalities in~\eqref{stabcond1} and~\eqref{stabcond2} for ${\mathcal T} = {\mathcal S}$ and $U = \{1, \dots, N\}$, respectively, are tight in the limit.
We will henceforth refer to the inequalities in~\eqref{stabcond3} and~\eqref{stabcond4} as the \textit{local stability conditions}.\\

For later use we observe that the inequalities in~\eqref{stabcond3} may also be written in the less intuitive but equivalent form
\begin{equation}
\frac{1}{N \mu} \sum_{n \in \big(\bigcup\limits_{S: S \not\in {\mathcal T}'} S\big)^c} \mu_n < \sum_{S \in {\mathcal T}'} p_S
\end{equation}
for all (non-empty) ${\mathcal T}' = {\mathcal T}^c \subsetneq {\mathcal S}$.
These inequalities reflect that for all (non-empty) ${\mathcal T}'$ the total arrival rate of job types $S \in {\mathcal T}'$ as $\lambda \uparrow \mu$ should become strictly higher than the aggregate service rate of the servers that are able to handle jobs of these types only.
In particular, taking ${\mathcal T}' = \{S_0\}$, we obtain
\begin{equation}
\frac{1}{N \mu} \sum_{n \in {\mathcal N}_{S_0}^c} \mu_n <  p_{S_0},
\label{eq:local_stab_interpret}
\end{equation}
with ${\mathcal N}_{S_0} = \bigcup_{S \in {\mathcal S} \setminus {S_0}} S$ and $S_0$ an arbitary job type belonging to $\mathcal{S}$.
Likewise, the conditions in~\eqref{stabcond4} may be rewritten as
\begin{equation}
\frac{1}{N \mu} \sum_{n \in U'} \mu_n < \sum_{S: S \cap U' \neq \emptyset} p_S
\end{equation}
for all (non-empty) $U' = U^c \subsetneq \{1, \dots, N\}$.
These inequalities indicate that for all (non-empty) $U' \subsetneq \{1, \dots, N\}$, the aggregate service rate of the servers in the set~$U$ should become strictly lower than the total arrival rate of the job types that can be handled by at least one of these servers as $\lambda \uparrow \mu$.
In particular, taking $U' = \{n_0\}$, we see that any server~$n_0$ (with a strictly positive speed $\mu_{n_0} > 0$) must be able to handle at least one job type~$S$ (with a strictly positive arrival probability $p_S > 0$).

The next theorem establishes that detailed performance metrics, such as the queue lengths and delays, are not strongly affected by the exact values of the probabilities $(p_S)_{S\in\mathcal{S}}$ as long as the local stability conditions~\eqref{stabcond3} and~\eqref{stabcond4} hold. We will write $\xrightarrow{d}$ to denote convergence in distribution of random variables, and $\mathrm{Exp}(1)$ will represent a unit-mean exponentially distributed random variable.

\begin{theorem}
\label{statespacecollapse}
If the local stability conditions~\eqref{stabcond3} and~\eqref{stabcond4} hold, then for both the c.o.c.\ and c.o.s.\ mechanisms
\[
\left(1 - \frac{\lambda}{\mu}\right) (Q_S)_{S \in {\mathcal S}} \xrightarrow{d}
\mathrm{Exp}(1) (p_S)_{S \in {\mathcal S}},
\]
as $\lambda \uparrow \mu$. 
\end{theorem}

The above theorem shows that the joint distribution of the number of jobs of each type, under both the c.o.c.\ and c.o.s.\ mechanisms, exhibits state space collapse in a heavy-traffic regime.
Moreover, the joint stationary distribution coincides in the limit with the joint distribution of a multi-class M/M/1 queue with arrival rate $N \lambda$, service rate $N \mu$, and class probabilities~$(p_S)_{S\in\mathcal{S}}$, provided the local stability conditions~\eqref{stabcond3} and~\eqref{stabcond4} are satisfied.
In particular, the total number of jobs, properly scaled, tends to an exponential random variable in the limit, and  thus the number of jobs of each type has an exponential limiting distribution as well.\\ 

Theorem~\ref{statespacecollapse} may be interpreted as follows.
It is highly unlikely that any server is idling when the total number of jobs is large because of the natural diversity in job types and the fact that any server is able to handle at least one job type as noted above.
In fact, although the parallel-server system with a c.o.c.\ mechanism is not work-conserving, the probability of any server being idle will vanish in the heavy-traffic limit.
This can be formalized based on the arguments used to prove Theorem~\ref{statespacecollapse}.
For the c.o.s.\ mechanism this is already shown in~\cite[Theorem~3.1]{adan2014skill}.
Hence, the system will operate at the full aggregate service rate $N \mu$ with high probability when the total number of jobs is sufficiently large.
This explains why the total number of jobs behaves asymptotically the same as in an M/M/1 queue with arrival rate $N \lambda$ and service rate $N \mu$, and in particular follows the well-known scaled exponential distribution in the limit. In the special case of the power-of-$d$ policy, this result was already shown in~\cite{ayesta2018unifying} for the c.o.s.\ mechanism.

What is far less evident though, is the state space collapse, i.e., the proportion of type-$S$ jobs present in the system coincides with the corresponding arrival probability $p_{S}$ for each job type $S\in\mathcal{S}$ like in a multi-class M/M/1 queue with a FCFS discipline.
In order to provide an informal explanation, suppose that for a particular type~$S_0$ the proportion of jobs is significantly lower than the corresponding arrival probability~$p_{S_0}$ while the total number of jobs is large.
Hence, in order to observe this unexpectedly low fraction of type-$S_0$ jobs, a large number of these type-$S_0$ jobs have been completed that arrived after jobs of other types that are still in the system.
This in turn implies that type-$S_0$ jobs will only receive a non-vanishing fraction of the capacity of the servers in $\mathcal{N}^c_{S_0}$ that cannot handle jobs of any other types, as all other compatible servers are processing earlier arrived jobs of different types. Here $\mathcal{N}_{S_0}$ is defined as $\cup_{S \in \mathcal{S} \setminus \{S_0\}} S$, the set of all servers compatible with the job types $\mathcal{S} \setminus \{S_0\}$. 
Since the type-$S_0$ jobs are not accumulating in the system, the aggregate service rate of the servers in $\mathcal{N}_{S_0}^c$ is no less than the arrival rate of the type-$S_0$ jobs, i.e., the inequality in~\eqref{eq:local_stab_interpret} is violated, which yields a contradiction.
Hence, we conclude that for none of the job types the proportion of jobs can be significantly lower than the corresponding arrival probability, meaning that the state space collapse occurs.\\
 
The local stability conditions~\eqref{stabcond3} and~\eqref{stabcond4} may be interpreted as so-called \textit{Complete Resource Pooling} (CRP) conditions.
Such conditions have emerged as an ubiquitous concept in the heavy-traffic behavior of stochastic networks, and in particular play a paramount role in the occurrence of a (one-dimensional) state space collapse.
CRP conditions have been formulated in roughly three different, yet related, ways in the literature:
(i) linear programming characterizations; (ii) geometric interpretations; and (iii) systems of linear inequalities.
All three notions pertain to the relative position of the critical load vector on the boundary of the capacity region of the system, and obviously involve the relevant system parameters (e.g. service rates, compatibility constraints between job types and servers, proportions of job types).
However, the various representations differ in the degree to which either the parameter values or the intrinsic properties of interest are explicitly covered.

The first type of characterization as introduced by Harrison~\cite{Harrison1998} and Harrison and L\'opez~\cite{Harrison1999} requires that the solution to a certain linear program (the dual version of a linear program describing feasible resource allocation options) is unique.
The second representation stipulates that the normal vector to the boundary of the capacity region at the critical load vector is unique with strictly positive components, see for instance~\cite{Mandelbaum2004,Stolyar2004,hurtado2020transform}.
The third kind of characterization involves a system of linear equalities in terms of the system parameters, see for instance the seminal papers by Laws~\cite{Laws1992} and Kelly and Laws~\cite{Kelly1993} which refer to these inequalities as `generalized cut constraints', and the more recent studies in~\cite{Shi2019,Banerjee2020,Varma2021} which develop related notions.

The local stability conditions~\eqref{stabcond3} and~\eqref{stabcond4} are similar in spirit as the ones in the latter category, which take a particularly convenient form in the context of a parallel-server system.
We have opted to state the CRP condition in this form, since the linear equalities provide an explicit characterization in terms of the system parameters which is straightforward to verify and additionally captures how the CRP condition is used in the proof arguments.
However, the local stability conditions can be shown to be equivalent to a conceptually similar linear programming construction as in~(i) or a geometric representation in a similar vein as in~(ii).
The latter two characterizations, both specialized to the setting of a parallel-server system with compatibility constraints, are provided and further discussed in Appendix~\ref{app:CRP}, where the equivalence with the inequalities in~\eqref{stabcond3} and~\eqref{stabcond4} is also established.

 \subsection{Heavy-traffic regime in a general scenario}\label{subsec:general}
While the local stability conditions~\eqref{stabcond3} and~\eqref{stabcond4} are natural and desirable design objectives, one can definitely conceive scenarios where these inequalities may not be satisfied. In order to illustrate this,
 Figure~\ref{fig:example4} visualizes the capacity region for the system depicted in Figure~\ref{fig:example2}. As can be seen here, for probabilities $(p_S)_{S\in\mathcal{S}}$ for which the vector $\lambda(p_S)_{S\in\mathcal{S}}$ lies within the cone region (indicated in gray), this vector $\lambda(p_S)_{S\in\mathcal{S}}$ will indeed reach the boundary of the capacity region corresponding to the stability condition $\lambda<\mu$ when $\lambda$ increases. In other words, the local stability conditions are met. This implies that the critical subset is $\mathcal{T}^*=\mathcal{S}$ and the normalized critical relative arrival rate is $\lambda^*=\mu$. 
For probabilities $(p_S)_{S\in\mathcal{S}}$ with corresponding arrival rate vectors outside this cone region, however, the vector $\lambda(p_S)_{S\in\mathcal{S}}$ could already reach the boundary of the capacity region when $\lambda$ approaches some $\lambda^*<\mu$, and hence, $\mathcal{T}^*$ is a strict subset of $\mathcal{S}$.
In general, starting from the stability conditions in~\eqref{stabcond1} and increasing~$\lambda$ will result in a subset $\mathcal{T}^*\subseteq\mathcal{S}$ of job types such that the inequality $N\lambda p_{\mathcal{T}^*} < \mu_{\mathcal{T}^*}$ becomes tight when $\lambda \uparrow \lambda^*$. Note that $\mathcal{T}^*$ and $\lambda^*$ are precisely the critical subset and critical normalized arrival rate as defined in Definition~\ref{def:critical_subset_arrrate}. 

We now extend Theorem~\ref{statespacecollapse} to such a scenario with $\lambda^*<\mu$ and $\mathcal{T}^*\neq\mathcal{S}$ under the following mild assumption.


\begin{assumption}\label{assumption:non_local_stab}
Let the capacity ratio~$\rho_{\mathcal{T}}$ for all $\mathcal{T}\subseteq\mathcal{S}$ and the critical subset~$\mathcal{T}^*$ be as specified in Definition~\ref{def:critical_subset_arrrate}. Then the subset~$\mathcal{T}^*$ that maximizes the capacity ratio $\rho_{\mathcal{T}}$ is unique,
i.e., the vector of  critical arrival rates can reach any point on the boundary of the capacity region except for its extreme vertices.
\end{assumption}

\begin{theorem}
\label{statespacecollapse_general}
Let $\mathcal{T}^*$, $p_{\mathcal{T}^*}$ and $\lambda^*$ be as in Definition~\ref{def:critical_subset_arrrate}, satisfying Assumption~\ref{assumption:non_local_stab}. Then for both the c.o.c.\ and c.o.s.\ mechanisms
\[
\left(1 - \frac{\lambda}{\lambda^*}\right) (Q_S)_{S \in {\mathcal S}} = \left(1 - \frac{\lambda}{\lambda^*}\right) \left((Q_S)_{S \in {\mathcal T}^*},(Q_S)_{S \notin {\mathcal T}^*}\right) \xrightarrow{d}
\left(\mathrm{Exp}(1) \left(\frac{p_S}{p_{\mathcal{T}^*}}\right)_{S \in {\mathcal T}^*},(0)_{S \notin {\mathcal T}^*}\right),
\]
as $\lambda \uparrow \lambda^*$. 
\end{theorem}

Hence, a similar state space collapse occurs as in Theorem~\ref{statespacecollapse} for the number of jobs of each type in~$\mathcal{T}^*$. On the other hand, the number of jobs of types not in $\mathcal{T}^*$, i.e., job types for which the aggregate arrival rate to the compatible servers is sub-critical, becomes negligible after scaling.


From Theorem~\ref{statespacecollapse_general} we can derive the following two corollaries. 

\begin{corollary}\label{cor:R_i}
Let $\mathcal{T}^*$, $p_{\mathcal{T}^*}$ and $\lambda^*$ be as in Definition~\ref{def:critical_subset_arrrate}, satisfying Assumption~\ref{assumption:non_local_stab}. Then for both the c.o.c.\ and c.o.s.\ mechanisms
\[
\left(1 - \frac{\lambda}{\lambda^*}\right) (R_n)_{n = 1, \dots, N} \xrightarrow{d} \mathrm{Exp}(1) (q_n)_{n = 1, \dots, N},
\]
as $\lambda \uparrow \lambda^*$, with 
$q_n = \sum_{S\in \mathcal{T}^*: n \in S} \frac{p_S}{p_{\mathcal{T}^*}}$, $n = 1, \dots, N$. 
\end{corollary}
Note that the summation in the definition of $q_n$ is equal to zero when server $n$ is not compatible with any of the job types in $\mathcal{T}^*$. Hence its queue length becomes negligible, after scaling, in the heavy-traffic regime.

By virtue of the distributional form of Little's law~\cite{Keilson1988}, the sojourn time and waiting time of a type-$S$ job, properly scaled, also have an exponential distribution in the limit when $S\in\mathcal{T}^*$. When  $S\notin\mathcal{T}^*$, we know from the previous result that the (scaled) queue lengths at its compatible servers are 0, hence its (scaled) sojourn time and waiting time are 0 as well.

\begin{corollary}\label{cor:SojournWaiting}
Let $\mathcal{T}^*$, $p_{\mathcal{T}^*}$ and $\lambda^*$ be as in Definition~\ref{def:critical_subset_arrrate}, satisfying Assumption~\ref{assumption:non_local_stab}. Then for both the c.o.c.\ and c.o.s.\ mechanisms, for $S\in \mathcal{T}^*$
\[
\left(1 - \frac{\lambda}{\lambda^*}\right) V_S \xrightarrow{d} \rho_{\mathcal{T}^*}\mathrm{Exp}(1),
\hspace{.4in}
\left(1 - \frac{\lambda}{\lambda^*}\right) W_S \xrightarrow{d} \rho_{\mathcal{T}^*}\mathrm{Exp}(1),
\]
and for $S\notin\mathcal{T}^*$
\[
\left(1 - \frac{\lambda}{\lambda^*}\right) V_S \xrightarrow{d} 0,
\hspace{.4in}
\left(1 - \frac{\lambda}{\lambda^*}\right) W_S \xrightarrow{d}        0,
\]
as $\lambda \uparrow \lambda^*$. 
\end{corollary}

This result shows that, for instance when the local stability conditions in~\eqref{stabcond3} and~\eqref{stabcond4} hold, i.e., $\mathcal{T}^*\equiv\mathcal{S}$, the sojourn time and waiting time of a particular type of job asymptotically do not depend on the probabilities $(p_S)_{S\in\mathcal{S}}$ in any way as $\rho_{\mathcal{T}^*} = \frac{p_{\mathcal{T}^*}}{\mu_{\mathcal{T}^*}} = \frac{1}{N\mu}$. Thus, surprisingly, jobs of types with more compatible servers do not enjoy significantly shorter waiting times or sojourn times in the limit. \\



The proof of Theorem~\ref{statespacecollapse_general} for the c.o.c.\ mechanism, presented in Section~\ref{sec:HT}, relies on a specific
enumeration of all possible job configurations, which yields a particularly
convenient form of the joint Probability Generating Function~(PGF) of the
number of jobs of each type as provided in the next proposition.

\begin{proposition}\label{prop:pgf}
Assuming that the stability conditions~\eqref{stabcond1} and~\eqref{stabcond2} are satisfied, the joint~PGF of the number of jobs of each type for the redundancy c.o.c.\ policy is given by 
\begin{equation}\label{eq:pgf}
\mathbb{E}\left[\prod\limits_{S \in {\mathcal S}} z_S^{Q_S}\right] = \frac{f(\boldsymbol{z})}{f(\boldsymbol{1})},
\end{equation}
where $\boldsymbol{z}$ and $\boldsymbol{1}$  are $|\mathcal{S}|$-dimensional vectors with entries $|z_S| \le 1$ and 
\begin{equation}\label{eq:pgf_f}
 f(\boldsymbol{z}) = 1{+} \sum\limits_{m=1}^{|\mathcal{S}|}\sum\limits_{\boldsymbol{S}\in \mathcal{S}_m} \prod\limits_{j=1}^m \frac{N\lambda p_{S_j}z_{S_j}}{\mu(S_1,\dots,S_j)}\prod\limits_{j=1}^m\left(1{-}\frac{N\lambda}{\mu(S_1,\dots,S_j)}\sum\limits_{i=1}^j p_{S_i}z_{S_i}\right)^{-1}.
\end{equation}
The $m$-dimensional vector $\boldsymbol{S}$ consists of $m$ different job types, and the set consisting of all these vectors is denoted by $\mathcal{S}_m$.
\end{proposition}

An interpretation of this PGF in terms of the ordered vectors $\boldsymbol{S}$ of job types and geometrically distributed random variables can be found in Appendix~\ref{app:interpretation_pgf}.
The proof of Proposition~\ref{prop:pgf} for the c.o.c.\ mechanism and the close relationship between the product-form expressions~\eqref{eq:statdistr} and~\eqref{eq:product_form_cos} can be exploited to prove Theorem~\ref{statespacecollapse_general} for the c.o.s.\ mechanism. This derivation is given in Appendix~\ref{app:cos_indirect}.
Moreover, a proof technique similar to the one outlined for the c.o.c.\ mechanism can also be applied to prove the results for the c.o.s.\ mechanism, where first the joint~PGF of the number of waiting jobs is derived before taking the appropriate limit. This proof is deferred to Appendix~\ref{app:cos_direct}.\\

Having exact expressions for the PGFs of the number of jobs of each type as in Propositions~\ref{prop:pgf} and~\ref{prop:pgfcos} allows us to establish convergence of the (scaled) $n$th moments of both $Q$ and $Q_S$, besides the convergence in distribution in Theorem~\ref{statespacecollapse}.
\begin{theorem}\label{th:conv_in_mean}
If the local stability conditions~\eqref{stabcond3} and~\eqref{stabcond4} hold, then for both the c.o.c.\ and c.o.s.\ mechanisms
\begin{equation}
\lim_{\lambda\uparrow\mu} \mathbb{E}\left[\left(\left(1-\frac{\lambda}{\mu}\right)Q\right)^n \right] = n! = \mathbb{E}\left[\left(\mathrm{Exp}(1)\right)^n \right] 
\end{equation}
and
\begin{equation}
\lim_{\lambda\uparrow\mu} \mathbb{E}\left[\left(\left(1-\frac{\lambda}{\mu}\right)Q_S\right)^n \right] = n!(p_S)^n = \mathbb{E}\left[\left(\mathrm{Exp}(1)p_S\right)^n \right] 
\end{equation}
for any $n\ge 1$ and $S\in\mathcal{S}$. 
\end{theorem}
The proofs of Theorem~\ref{th:conv_in_mean}, including explicit expressions for the $n$th moments of $Q$ and $Q_S$, are deferred to Appendix~\ref{app:conv_moments}.\\

The above results focus on the job types that experience critical load, and show that the contribution of the non-critical job types to the total number of jobs becomes negligible after scaling. Hence, the full system comprises two subsystems. The critical subsystem consists of all job types in $\mathcal{T}^*$ and their compatible servers and the non-critical subsystem consists of all job types $\mathcal{S}\setminus\mathcal{T}^*$ and their compatible servers. Note that the server sets of both subsystems do not need to be disjoint.
However, along the same lines as the proof of Theorem~\ref{statespacecollapse_general}, it can be shown that the non-critical subsystem in the heavy-traffic regime operates as an isolated system with servers that are only compatible with types outside of $\mathcal{T}^*$, referred to as the \textit{truncated} system, as formalized in the following theorem.
\begin{theorem}\label{th:ssc_nonscaled}
Let $\mathcal{T}^*$, $p_{\mathcal{T}^*}$ and $\lambda^*$ be as in Definition~\ref{def:critical_subset_arrrate}, satisfying Assumption~\ref{assumption:non_local_stab}.
For the c.o.c. policy 
\begin{equation}
\left(\left(1-\frac{\lambda}{\lambda^*}\right)(Q_S)_{S\in \mathcal{T}^*},(Q_S)_{S\notin \mathcal{T}^*}\right) \xrightarrow{d} \left(\mathrm{Exp}(1)\left(\frac{p_S}{p_{\mathcal{T}^*}}\right)_{S\in \mathcal{T}^*}, (Q_S^*)_{S\notin \mathcal{T}^*}\right)
\end{equation}
as $\lambda\uparrow\lambda^*$, and for the c.o.s. policy 
\begin{equation}
\left(\left(1-\frac{\lambda}{\lambda^*}\right)(\tilde{Q}_S)_{S\in \mathcal{T}^*},(\tilde{Q}_S)_{S\notin \mathcal{T}^*}\right) \xrightarrow{d} \left(\mathrm{Exp}(1)\left(\frac{p_S}{p_{\mathcal{T}^*}}\right)_{S\in \mathcal{T}^*}, (\tilde{Q}_S^*)_{S\notin \mathcal{T}^*}\right)
\end{equation}
as $\lambda\uparrow\lambda^*$. Here $Q_S^*$ and $\tilde{Q}_S^*$ denote the number of type-$S$ jobs and the number of waiting type-$S$ jobs, respectively, in a truncated system as described above. 
\end{theorem}
The above theorem implies that the full system decomposes into two disjoint and independent subsystems in the heavy-traffic regime. 
A more formal definition of the truncated system, additional information and a sketch of the proof are deferred to Appendix~\ref{sec:non_critial_jobtypes}. 

\subsection{Numerical results and discussion}\label{subsec:discussion}
In this subsection we present numerical results to illustrate the heavy-traffic limits and discuss some design implications.


\paragraph{Numerical illustration.}
We focus on a system with three
identical servers and two job types labeled $\{1, 2\}$ and $\{2, 3\}$ as schematically represented in Figure~\ref{fig:example2}.
For each arriving job, replicas are assigned to either servers~$1$ and~$2$ or servers~$2$ and~$3$, respectively, depending on its type.
The joint stationary distribution of the number of jobs under the c.o.c.\ mechanism can be derived through some straightforward but tedious calculations from the more detailed stationary distribution~\eqref{eq:statdistr}.
Specifically, $\mathbb{P}\{(Q_{\{1,2\}},Q_{\{2,3\}}) = (q,q')\}$ equals  \begin{equation}
\left\{
\begin{array}{lcl}
C & & \text{if } q=q'=0 \\
C \left(\frac{N\lambda p_{\{1,2\}}}{2\mu}\right)^{q} & &  \text{if } q\ge 1, q'=0 \\
C \left(\frac{N\lambda p_{\{2,3\}}}{2\mu}\right)^{q'} & &  \text{if } q=0, q'\ge 1 \\
C \left(\frac{N\lambda p_{\{1,2\}}}{3\mu}\right)^{q}\left(\frac{N\lambda p_{\{2,3\}}}{3\mu}\right)^{q'}\left[
\sum\limits_{k=1}^{q}\binom{q+q'-k-1}{q'-1}\left(\frac{3}{2}\right)^{k} +
\sum\limits_{k=1}^{q'}\binom{q+q'-k-1}{q-1}\left(\frac{3}{2}\right)^{k}
\right] & & \text{if } q,q'\ge 1.
\end{array}
\right.
\end{equation}
The normalization constant~$C$ is given by
\[
\frac{3 (2\mu-N\lambda p_{\{1,2\}})(2\mu-N\lambda p_{\{2,3\}})(\mu-\lambda)}{(N \lambda)^2 p_{\{1,2\}} p_{\{2,3\}}\mu + 12\mu^2(\mu-\lambda)},
\]
as was also established in~\cite{Gardner2016queueing}.
In order to produce numerical results, we set the values of the job type fractions to $p_{\{1,2\}}=7/12$ and $p_{\{2,3\}}=5/12$ satisfying the local stability conditions in~\eqref{stabcond3} and~\eqref{stabcond4}, which in this particular case read $p_{\{1,2\}}, p_{\{2,3\}} < 2/3$.

\begin{figure}
\centering
\begin{subfigure}[t]{0.23\textwidth}
  \centering 
  \includegraphics[scale=0.5]{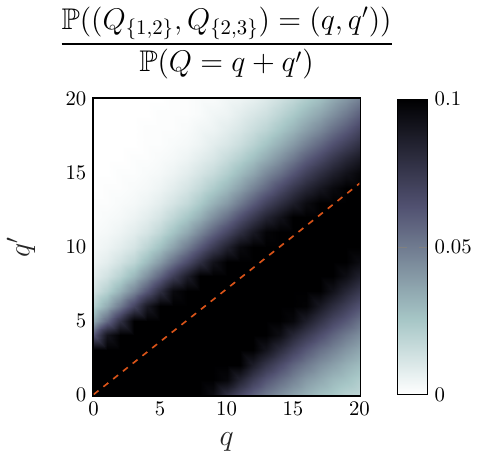}
  \caption{Joint probability mass function with $\lambda = 0.6\mu$.}
  \label{fig:HT_ssc06}
\end{subfigure}%
\hspace{0.2cm}
\begin{subfigure}[t]{0.23\textwidth}
\centering 
\includegraphics[scale=0.5]{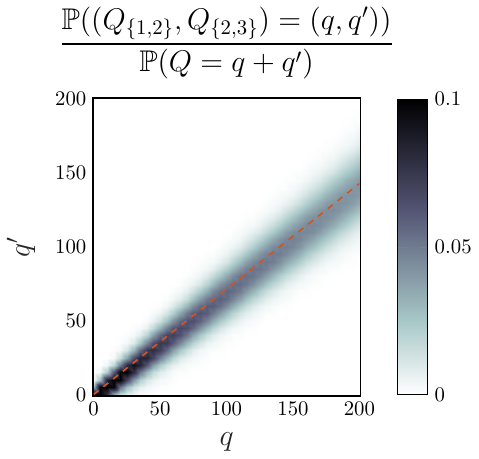}
  \caption{Joint probability mass function with $\lambda = 0.99\mu$.}
  \label{fig:HT_ssc}
\end{subfigure}
\hspace{0.2cm}
\begin{subfigure}[t]{0.23\textwidth}
  \centering 
  \includegraphics[scale=0.5]{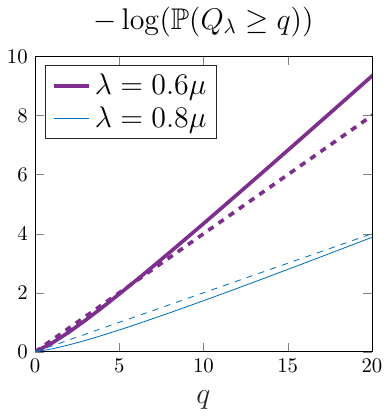}
  \caption{Exponential approximation.}
  \label{fig:HT_ExpAppr_small}
\end{subfigure}%
\hspace{0.2cm}
\begin{subfigure}[t]{0.23\textwidth}
\centering 
\includegraphics[scale=0.5]{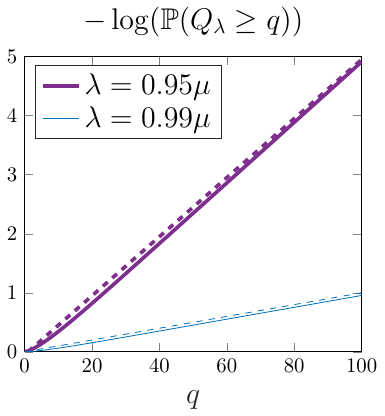}
  \caption{Exponential approximation.}
  \label{fig:HT_ExpAppr}
\end{subfigure}
\caption{Numerical illustration of the result in Theorem~\ref{statespacecollapse}. The system consists of three identical servers, has type fractions $(p_{\{1,2\}},p_{\{2,3\}}) = (\frac{7}{12},\frac{5}{12})$ and a compatibility graph as indicated in Figure~\ref{fig:example}. Figures~(a) and~(b) visualize the (scaled) probability mass function $\mathbb{P}\{(Q_{\{1,2\}},Q_{\{2,3\}}) = (q,q')\}$ when $\lambda = 0.6\mu$ and $\lambda = 0.99\mu$, respectively. The (dashed) line 
with slope $5/7$ depicts the limiting regime. Figures~(c) and~(d) give a comparison between the cumulative distributions of the total number of jobs (solid line) and the random variable $\mathrm{Exp}(1)(1-\frac{\lambda}{\mu})^{-1}$ (dashed line) for various values of $\lambda$.} \label{fig:HT}
\end{figure}

Figures~\ref{fig:HT_ssc06} and~\ref{fig:HT_ssc} corroborate the state space collapse result: The (scaled) probability mass function strongly concentrates around the dashed line $\{(p_{\{1,2\}}, p_{\{2,3\}}) x \colon x \ge 0\}$ with slope $5/7$, which represents the heavy-traffic limit.
Remarkably, this holds not only for a high load of~$0.99$, but already manifests itself for a moderate load value of~$0.6$ (note the ten-fold difference in scale). 

Figures~\ref{fig:HT_ExpAppr_small} and~\ref{fig:HT_ExpAppr} confirm that the total number of jobs, properly scaled, converges to a unit exponential random variable, i.e., $(1-\frac{\lambda}{\mu}) Q_{\lambda} \rightarrow \mathrm{Exp}(1)$ as $\lambda\uparrow\mu$, with $Q_{\lambda} = Q_{\{1,2\}}+Q_{\{2,3\}}$.
We observe that the heavy-traffic approximation
\begin{equation}\label{eq:ht_approximation_numerical}
\mathbb{P}\{Q_{\lambda}\ge q\} = \mathbb{P}\left\{\left(1-\frac{\lambda}{\mu}\right)Q_{\lambda}\ge \left(1-\frac{\lambda}{\mu}\right)q\right\} \approx \mathbb{P}\left\{\mathrm{Exp}(1)\ge\left(1-\frac{\lambda}{\mu}\right)q\right\} = \mathrm{e}^{- \left(1-\frac{\lambda}{\mu}\right)q}
\end{equation}
is not only nearly exact for high load values like~$0.99$ and~$0.95$, but also fairly close for moderate load values like~$0.8$ and~$0.6$ outside the asymptotic regime.
Only for the tail probabilities the accuracy starts to diminish.
Note that for a load value of~$0.6$, the probability that for instance $q = 15$ is exceeded is around $4.7 \times 10^{- 4}$.

The above system with three servers and two job types is admittedly a toy model.
The product-form distributions in~\eqref{eq:statdistr} in fact hold for systems of any size and arbitrary compatibility constraints.
As noted earlier, however, these expressions do not lend themselves well for computational purposes in more complex situations.

\paragraph{Simulation results.} 
For the c.o.s.\ mechanism we conducted simulations for systems with $N=4, 10$ or 20 homogeneous servers and a variety of compatibility graphs. These compatibility graphs range from fairly sparse (e.g. $N$ job types with type-$n$ jobs only compatible with servers $n$ and $n+1$), to fairly dense (e.g. power-of-$d$ setting with $d$ close to $N$). It is observed that for a small system, the heavy-traffic limit in~\eqref{eq:ht_approximation_numerical} provides a reasonable approximation for the total number of jobs in the system, i.e., $Q_{\lambda} \approx \mathrm{Exp}(1)(1-\frac{\lambda}{\mu})^{-1}$. 
For larger systems, this still seems to hold true when the compatibility graph is rather dense, meaning that each job type is compatible with a large fraction of servers.
However, when a type-$S$ job arrives with probability $p_{S}$, we observe that the number of type-$S$ jobs is approximately equal to $p_{S}Q_{\lambda}$.
More detailed simulation results are provided in Appendix~\ref{app:simulations}.

\paragraph{Design implications.} 
The universality property embodied in Theorem~\ref{statespacecollapse} implies that the performance of a fully flexible system can asymptotically be matched as long as the local stability conditions are satisfied.
In practical terms, this means that creating a limited amount of flexibility in the assignment constraints, if done judiciously, is sufficient, and that an excess of flexibility yields only minor performance gains. Similar observations have been made in the context of process flexibility in manufacturing systems, see for instance \cite{Jordan1995,Graves2003,Simchi-Levi2012, Shi2019}.

To illustrate the above observation, suppose that we set the probabilities as
\[
p_n =
\frac{N - 1 - \epsilon}{N (N - 1) \mu} \mu_n \mbox{ for all } n = 1, \dots, N
\quad\text{and}\quad
p_{\{1, \dots, N\}} = \frac{\epsilon}{N - 1}.
\]
Hence, most jobs are assigned to a single server and a few jobs are replicated to all servers.
Then the average replication degree is $1 + \epsilon$ while for all (non-empty) subsets $U \subsetneq \{1, \dots, N\}$ we have for any $\epsilon > 0$,
\[
\sum_{S: S \subseteq U} p_S =
\sum_{n \in U} p_n = \frac{N - 1 - \epsilon}{N (N - 1) \mu} \sum_{n \in U} \mu_n <
\frac{1}{N \mu} \sum_{n \in U} \mu_n,
\]
implying that the local stability conditions in~\eqref{stabcond4} are satisfied.

Note that for $\epsilon = 0$ and a setting with identical servers, the above-described scenario reduces to a system of $N$~independent M/M/1 queues with arrival rate~$\lambda$ and service rate~$\mu$.
In that case the number of jobs at each server has a geometric stationary distribution with parameter $\lambda / \mu$, and when scaled with $1 - \lambda / \mu$, tends to a unit-exponential distribution as $\lambda\uparrow\mu$.  Thus the scaled total number of jobs converges to the sum of $N$~independent unit-exponential random variables.
In contrast, for any $\epsilon > 0$, Theorem~\ref{statespacecollapse} implies that the total number of jobs tends to a single unit-exponential random variable, and is hence smaller by a factor~$N$.

In the above example we observed that the performance at high load is roughly similar for all values $\epsilon > 0$, and significantly better than for $\epsilon = 0$.
This observation is reminiscent of the finding that a power-of-two policy provides a significant improvement over purely random assignment to a single server in a regime where the number of servers~$N$ grows large while $\lambda$ remains fixed~\cite{gardner2017redundancy}.
(Interestingly, the significant benefit of `just a little' flexibility has also been encountered in different resource sharing contexts, see for instance \cite{FS99,TX13a}.)
The fact that the heavy-traffic regime and a many-server scenario point to similar behavior also suggests that it would be interesting to explore joint scalings.

\section{Proofs}\label{sec:HT}


The proof of the heavy-traffic result stated in Theorem~\ref{statespacecollapse_general}, which implies the result in Theorem~\ref{statespacecollapse}, relies on a well-suited expression for the joint~PGF of the number of jobs of each type. This expression, which may be of independent interest, is provided in Proposition~\ref{prop:pgf}, and involves a specific enumeration of all possible job configurations as explained in the proof below.

\begin{proof}[Proof of Proposition~\ref{prop:pgf}.]
Using the stationary distribution given in \eqref{eq:statdistr}, the joint~PGF of the number of jobs of each type may be written as
\begin{equation}\label{eq:general_gen_f}
\mathbb{E}\left[\prod\limits_{S\in\mathcal{S}}z_S^{Q_S}\right]  =\sum\limits_{M=0}^{\infty} \sum\limits_{(c_1,\dots,c_M)\in \mathcal{S}^M} \pi_{\text{c.o.c.}}(c_1,\dots,c_M) \prod\limits_{S\in\mathcal{S}} z_S^{q_S^{\boldsymbol{c}}},
\end{equation}
with $\boldsymbol{z}$ an $|\mathcal{S}|$-dimensional vector with entries $|z_S|\le1$ and $q_S^{\boldsymbol{c}}$ the total number of type-$S$ jobs in state $\boldsymbol{c} = (c_1,\dots,c_M)$.
The summation on the right-hand side in \eqref{eq:general_gen_f} over all possible states can be conducted in three steps:
\begin{enumerate}
\item[(i)] Fix the number of different job types $m$ that occur in the state~$\boldsymbol{c}$, $m=1,\dots,|\mathcal{S}|$.
\item[(ii)] Fix $m$ distinct job types and the order in which they occur, $\boldsymbol{S} = [S_1,\dots,S_m]$. This implies that the oldest job in the system is of type $S_1$, the possible following jobs are of the same type, the first time a different type is observed it will be a job of type $S_2$, etc. The set containing all these vectors of length $m$ is denoted by $\mathcal{S}_m$.
\item[(iii)] Sum over all states with this particular order. For instance, for the vector $[S_1,S_2,S_3]\in\mathcal{S}_3$ with only three job types one has to sum over all states
\begin{equation}
\boldsymbol{c} = (S_1,\underbrace{S_1,S_1}_{k_1}, S_2, \underbrace{\times,\times,\times}_{k_2}, S_3, \underbrace{\circ,\circ,\circ}_{k_3}),
\end{equation}
where $\times$ denotes jobs of types $S_1$ and/or types $S_2$ and $\circ$ denotes jobs of types $S_1$, $S_2$ and/or $S_3$. The values $k_1$, $k_2$ and $k_3$ can be any natural number. 
\end{enumerate}
The third step might warrant some illustration. For example, the contribution of the ordered vector $[S_1,S_2,S_3]$  to \eqref{eq:general_gen_f} can be computed by first determining how many jobs there are present in total, $M\ge 3$. Then, the values of $k_1$, $k_2$ and $k_3 = M-k_1-k_2-3$ are set. Finally, the $k_2$ and $k_3$ intermediate jobs are labeled as type-$S_1$ or type-$S_2$ jobs and type-$S_1$, type-$S_2$ or type-$S_3$ jobs, respectively. Relying on~\eqref{eq:statdistr}, the total contribution to \eqref{eq:general_gen_f} is then given by
\begin{equation}
\begin{array}{c}
{C}{}\sum\limits_{M=3}^{\infty}\frac{N\lambda p_{S_1}z_{S_1}}{\mu(S_{1})} 
\sum\limits_{k_1=0}^{M-3}
\left[
\left(\frac{N\lambda p_{S_1}z_{S_1}}{\mu(S_{1})}\right)^{k_1}

\frac{N\lambda p_{S_2}z_{S_2}}{\mu(S_1,S_2)}   
\sum\limits_{k_2=0}^{M-3-k_1}
\left[
\left(\frac{N\lambda}{\mu(S_1,S_2)}\right)^{k_2} \sum\limits_{l=0}^{k_2}
\binom{k_2}{l}
\left(p_{S_1}z_{S_1}\right)^{l}\left(p_{S_2}z_{S_2}\right)^{k_2-l}
\right.
\right. \\

\left.
\left.
\times
 \frac{N\lambda p_{S_3}z_{S_3}}{\mu(S_1,S_2,S_3)} 
\left(\frac{N\lambda}{\mu(S_1,S_2,S_3)}\right)^{k_3}  
  \sum\limits_{l_1=0}^{k_3}\sum\limits_{l_2=0}^{k_3{-}l_1}
  \left[
\binom{k_3}{l_1,l_2,k_3{-}l_1{-}l_2}
\left(p_{S_1}z_{S_1}\right)^{l_1}\left(p_{S_2}z_{S_2}\right)^{l_2}\left(p_{S_3}z_{S_3}\right)^{k_3-l_1-l_2}
\right]
\right]
\right].
\end{array}
\end{equation}
Applying the multinomium of Newton leads to
\begin{equation}
\begin{array}{c}
C\prod\limits_{j=1}^3\frac{N\lambda p_{S_j}z_{S_j}}{\mu(S_{1},\dots,S_j)} 

\sum\limits_{M=3}^{\infty}
\sum\limits_{k_1=0}^{M-3}
\left[
\left(\frac{N\lambda}{\mu(S_{1})}p_{S_1}z_{S_1}\right)^{k_1}

\sum\limits_{k_2=0}^{M-3-k_1} 
\left[ 
\left(\frac{N\lambda}{\mu(S_1,S_2)}\right)^{k_2} (p_{S_1}z_{S_1}+p_{S_2}z_{S_2})^{k_2}
\right.
\right. \\

\left.\left.\times \left(\frac{N\lambda}{\mu(S_1,S_2,S_3)}\right)^{k_3} \left(p_{S_1}z_{S_1}+p_{S_2}z_{S_2}+p_{S_3}z_{S_3}\right)^{k_3}\right]\right].
\end{array}
\end{equation}
Interchanging the order of summation results in
\begin{equation}
\begin{array}{c}
C\prod\limits_{j=1}^3\frac{N\lambda p_{S_j}z_{S_j} }{\mu(S_{1},\dots,S_j)}

\left[\sum\limits_{k_1=0}^{\infty}\left(\frac{N\lambda p_{S_1}z_{S_1}}{\mu(S_{1})}\right)^{k_1}\right] 

 \left[\sum\limits_{k_2=0}^{\infty}\left(\frac{N\lambda}{\mu(S_1,S_2)}\right)^{k_2} \left(p_{S_1}z_{S_1}+p_{S_2}z_{S_2}\right)^{k_2}\right] \\

\times\sum\limits_{k_3=0}^{\infty}\left[\left(\frac{N\lambda}{\mu(S_1,S_2,S_3)}\right)^{k_3} \left(p_{S_1}z_{S_1}+p_{S_2}z_{S_2}+p_{S_3}z_{S_3}\right)^{k_3}\right].
\end{array}
\end{equation}
Due to the stability conditions in \eqref{stabcond1} and \eqref{stabcond2}, the expression for the infinite geometric sum may be applied to obtain
\begin{equation}\label{eq:proof_GF_inter}
C\prod\limits_{j=1}^3\frac{N\lambda p_{S_j}z_{S_j}}{\mu(S_{1},\dots,S_j)}  \prod\limits_{j=1}^3\left(1-\frac{N\lambda}{\mu(S_{1},\dots,S_j)} (p_{S_1}z_{S_1}+\dots+ p_{S_j}z_{S_j})\right)^{-1}.
\end{equation}
Generalizing the above reasoning and applying the above-mentioned three steps will give an expression for \eqref{eq:general_gen_f}, namely
\begin{equation}
    C\left[1+ \sum\limits_{m=1}^{|\mathcal{S}|}\sum\limits_{\boldsymbol{S}\in \mathcal{S}_m} \prod\limits_{j=1}^m \frac{N\lambda p_{S_j}z_{S_j}}{\mu(S_1,\dots,S_j)}  \prod\limits_{j=1}^m\left(1-\frac{N\lambda}{\mu(S_1,\dots,S_j)}\sum\limits_{i=1}^j p_{S_i}z_{S_i}\right)^{-1}\right].
\end{equation}
The three steps only consider states with at least one job and the contribution of the empty state can be seen in the additional `1'.
Now, substituting $z_S=1$ in~\eqref{eq:general_gen_f} for all $S\in\mathcal{S}$ should give 1 as a result, and this yields an expression for the normalization constant $C$. This concludes the derivation of the joint PGF~\eqref{eq:pgf}.
\end{proof}

A similar result is derived in~\cite{ayesta2019token} concerning the generating function of both the total and the waiting number of jobs of each type in the token-based central queue setting, which includes among others matching models and redundancy models. Due to a slightly different state description and an alternative enumeration of the states, these generating functions still consist of infinite sums including general expressions for the probability that particular servers are processing particular job types.

Circumventing the above obstacles by directly using the product-form expressions for the stationary distribution in~\eqref{eq:statdistr} allows us to study the heavy-traffic limit in Theorem~\ref{statespacecollapse_general} for redundancy c.o.c.\ by interchanging the summation and limit operator in the expression provided in Proposition~\ref{prop:pgf}.

\begin{proof}[Proof of Theorem~\ref{statespacecollapse_general} \textup{(}c.o.c.\ mechanism\textup{)}.]
 Let $\mathcal{T}^*\subseteq\mathcal{S}$ be the critical subset, $p_{\mathcal{T}^*}$, $\mu_{\mathcal{T}^*}$ and $\lambda^*=\mu_{\mathcal{T}^*}/(Np_{\mathcal{T}^*})$ as defined in Definition~\ref{def:critical_subset_arrrate}.
We will prove the following heavy-traffic behavior of the Moment Generating Function~(MGF) of the number of jobs of each type $(Q_S)_{S\in\mathcal{S}}$.
\begin{equation}\label{eq:genfuncconv}
\mathbb{E}\left[ \mathrm{exp}\left(-\left(1-\frac{\lambda}{\lambda^*}\right) \sum\limits_{S\in\mathcal{S}} t_S Q_S \right) \right] \rightarrow \left(1+\sum\limits_{S\in\mathcal{T}^*} \frac{p_S}{p_{\mathcal{T}^*}}t_S \right)^{-1},
\end{equation}
as $\lambda\uparrow\lambda^*$ and $t_S\ge 0$ for all $S\in\mathcal{S}$. 
Moreover, it can easily be seen that the MGF of the random vector $\boldsymbol{X}\coloneqq\left({\mathrm{Exp}}(1)(p_S/p_{\mathcal{T}^*})_{S\in\mathcal{T}^*},(0)_{S\notin\mathcal{T}^*}\right)$ is given by 
\begin{equation}\label{eq:mgf_exp}
\mathbb{E}\left[\prod_{S\in\mathcal{T}^*}\mathrm{exp}\left({-t_S\left(\frac{p_S}{p_{\mathcal{T}^*}}{\mathrm{Exp}}(1)\right)}\right)\right] = \mathbb{E}\left[\mathrm{exp}\left(-\left(\sum\limits_{S\in\mathcal{T}^*}\frac{p_S}{p_{\mathcal{T}^*}}t_S\right){\mathrm{Exp}}(1)\right)\right] = \left(1{+}\sum_{S\in\mathcal{T}^*}\frac{p_S}{p_{\mathcal{T}^*}}t_S\right)^{-1}.
\end{equation}
By Feller's Convergence Theorem~\cite{Feller1971}, the non-negative random vector $(1-\frac{\lambda}{\lambda^*})(Q_S)_{S\in\mathcal{S}}$ converges in distribution to the random vector $\boldsymbol{X}$ when its MGF converges pointwise to the MGF in~\eqref{eq:mgf_exp}.
Hence, it is sufficient to show that \eqref{eq:genfuncconv} holds, in order to conclude the result stated in Theorem~\ref{statespacecollapse_general}. 

To obtain the MGF of $(1-\frac{\lambda}{\lambda^*})(Q_S)_{S\in\mathcal{S}}$ define $z_S\coloneqq \mathrm{exp}(-(1-\frac{\lambda}{\lambda^*})t_S)$ and use the expression for the PGF in~\eqref{eq:pgf}. As we allow $t_S  \ge 0$, it follows that $|z_S|\le 1$. For any $m=1,\dots,|\mathcal{S}|$, for any $\boldsymbol{S} = [S_1,\dots,S_m]\in\mathcal{S}_m$ and for any $j=1,\dots,m$ it can be observed that
\begin{equation}\label{eq:obs1}
\lim_{\lambda\uparrow \lambda^*} \frac{N\lambda}{\mu(S_1,\dots,S_j)}p_{S_j}z_{S_j} = \frac{\mu_{\mathcal{T}^*}}{\mu(S_1,\dots,S_j)}\frac{p_{S_j}}{p_{\mathcal{T}^*}} \in (0,\infty).
\end{equation}
From Assumption~\ref{assumption:non_local_stab}, it can be deduced that there exists some~$\epsilon>0$ such that for all $\mathcal{T}\subsetneq\mathcal{T}^*$ there holds that $\frac{N\lambda^* p_{\mathcal{T}}}{\mu_{\mathcal{T}}} = \frac{\mu_{\mathcal{T}^*}}{\mu_{\mathcal{T}}} \frac{p_{\mathcal{T}}}{p_{\mathcal{T}^*}}<1-\epsilon$.
 Hence,
\begin{equation}\label{eq:obs2}
\lim_{\lambda\uparrow \lambda^*}\left(1-\frac{N\lambda}{\mu(S_1,\dots,S_j)}\sum\limits_{i=1}^j p_{S_i}z_{S_i} \right)^{-1} 
= 
\begin{cases}
\infty & \text{if } \{S_1,\dots,S_j\} = \mathcal{T}^*,\\
\left(1-\frac{\mu_{\mathcal{T}^*}}{\mu(S_1,\dots,S_j)}\sum\limits_{i=1}^j \frac{p_{S_i}}{p_{\mathcal{T}^*}} \right)^{-1} \in (0,\infty) & \text{otherwise}.
\end{cases}
\end{equation}
Therefore the dominating terms in both the numerator and denominator of \eqref{eq:pgf} in the heavy-traffic regime are those with $\boldsymbol{S}\in \mathcal{S}_m$ such that $m\ge|\mathcal{T}^*|$ and $\{S_1,\dots,S_{|\mathcal{T}^*|}\} = \mathcal{T}^*$. Let $\mathcal{S}^{\mathcal{T}^*}$ denote the set of all vectors $\boldsymbol{S}$ that satisfy this property. Note that if $\mathcal{T}^*$ is contained in $\boldsymbol{S}$, but not as the $|\mathcal{T}^*|$ first occurring job types, this $\boldsymbol{S}$ will only have a finite contribution to the value in the numerator and denominator due to the above observations. This leads to 
\begin{equation}
 \lim\limits_{\lambda\uparrow \lambda^*} \mathbb{E} \left[ \mathrm{exp}\left(-\left(1-\frac{\lambda}{\lambda^*}\right) \sum\limits_{S\in\mathcal{S}} t_S Q_S \right) \right] 
 =  \lim\limits_{\lambda\uparrow \lambda^*} \dfrac{\sum\limits_{\boldsymbol{S}\in \mathcal{S}^{\mathcal{T}^*}} \prod\limits_{j=1}^{|\boldsymbol{S}|}  \frac{N\lambda p_{S_j}z_{S_j}}{\mu(S_1,\dots,S_j)} \prod\limits_{j=1}^{|\boldsymbol{S}|}\left(1-\frac{N\lambda}{\mu(S_1,\dots,S_j)}\sum\limits_{i=1}^j p_{S_i}z_{S_i}\right)^{-1}}{\sum\limits_{\boldsymbol{S}\in \mathcal{S}^{\mathcal{T}^*}} \prod\limits_{j=1}^{|\boldsymbol{S}|} \frac{N\lambda p_{S_j}}{\mu(S_1,\dots,S_j)} \prod\limits_{j=1}^{|\boldsymbol{S}|}\left(1-\frac{N\lambda}{\mu(S_1,\dots,S_j)}\sum\limits_{i=1}^j p_{S_i}\right)^{-1}}.
\end{equation}
Since $\mu(S_1,\dots,S_{|\mathcal{T}^*|}) = \mu_{\mathcal{T}^*}$ and $\sum_{i=1}^{|\mathcal{T}^*|}p_{S_i} = p_{\mathcal{T}^*}$, the above fraction can be rewritten as
\begin{equation}
\lim\limits_{\lambda\uparrow \lambda^*} \dfrac{\sum\limits_{\boldsymbol{S}\in \mathcal{S}^{\mathcal{T}^*}} \prod\limits_{j=1}^{|\boldsymbol{S}|} \ \frac{N\lambda p_{S_j}z_{S_j}}{\mu(S_1,\dots,S_j)} \prod\limits_{\substack{j=1\\ j \neq |\mathcal{T}^*|}}^{|\boldsymbol{S}|}\left(1{-}\frac{N\lambda}{\mu(S_1,\dots,S_j)}\sum\limits_{i=1}^j p_{S_i}z_{S_i}\right)^{-1}}{\sum\limits_{\boldsymbol{S}\in \mathcal{S}^{\mathcal{T}^*}} \prod\limits_{j=1}^{|\boldsymbol{S}|}  \frac{N\lambda p_{S_j}}{\mu(S_1,\dots,S_j)} \prod\limits_{\substack{j=1\\ j \neq |\mathcal{T}^*|}}^{|\boldsymbol{S}|}\left(1{-}\frac{N\lambda}{\mu(S_1,\dots,S_j)}\sum\limits_{i=1}^j p_{S_i}\right)^{-1}}\cdot 
 \lim\limits_{\lambda\uparrow \lambda^*} \dfrac{1-\frac{N\lambda}{\mu_{\mathcal{T}^*}}p_{\mathcal{T}^*}}{1-\frac{N\lambda}{\mu_{\mathcal{T}^*}}\sum\limits_{S\in\mathcal{T}^*} p_{S}z_{S}}.
\end{equation}
The first limit evaluates to 1 due to the above observations and, after applying l'H\^{o}pital's rule, the second limit indeed evaluates to the right-hand side of~\eqref{eq:mgf_exp}.
This concludes the proof of Theorem~\ref{statespacecollapse}.
\end{proof}

The proofs of Corollaries~\ref{cor:R_i} and~\ref{cor:SojournWaiting} are deferred to Appendix~\ref{app:HT}.

\section{Outlook}
\label{sec:conclusion}
The broader lay of the land and paucity of results for parallel-server systems with arbitrary assignment constraints as visualized in Figure~\ref{fig:diagram_literature} suggest a few natural directions for further research.

First of all, the papers of Rutten and Mukherjee~\cite{Rutten2020} and in particular of Weng et al.~\cite{Weng2020} are the primary counterparts of the present paper for JSQ type strategies in a many-server scenario rather than redundancy scheduling in a heavy-traffic regime.
Surprisingly, the results in these two papers also entail a certain notion of universality, with similar achievable performance as in a fully flexible system under relatively stringent assignment constraints.
While this universality property manifests itself in a different form in a many-server scenario, it suggests that this paradigm may not just apply to a given policy in a given limiting regime, but in fact unifies and spans across different conditions and different policies, with JSW providing a natural bridge between JSQ and redundancy scheduling as mentioned earlier.
In particular, we strongly conjecture that similar heavy-traffic results as derived in the present paper for redundancy scheduling hold for JSQ policies (up to speed dependent weight factors), except that these would need to be established in terms of process-level limits as they are obtained in~\cite{Atar2019} for power-of-$d$ settings in the absence of any explicit stationary distribution.
Likewise, it would be interesting to investigate whether similar many-server asymptotics as obtained in~\cite{Rutten2020,Weng2020} apply for redundancy scheduling models, for which the same product-form distributions as used in the present paper would provide a natural toolset. 

A further research direction would be to use the PGFs to establish convergence rates and refined approximations with improved accuracy in pre-limit scenarios of moderate instead of high load. Although the PGFs exist in closed form, the analysis would be both numerically and analytically challenging due to the intricate dependencies on the assignment constraints. 


\section*{Acknowledgements}
We are grateful to the two anonymous reviewers and the associate editor for their insightful comments and would further like to thank C\'{e}line Comte for helpful discussions that have led to the generalized result in Theorem~\ref{statespacecollapse_general}.


\begin{spacing}{0.2}
\bibliographystyle{abbrv}
\bibliography{references_OR} 

\begin{thebibliography}{10}

\bibitem{Adan2018}
I.~Adan, I.~Kleiner, R.~Righter, and G.~Weiss.
\newblock {FCFS parallel service systems and matching models}.
\newblock {\em Perform. Eval.}, 127-128:253--272, Nov. 2018.

\bibitem{adan2014skill}
I.~Adan and G.~Weiss.
\newblock A skill based parallel service system under {FCFS-ALIS} {—} steady
  state, overloads, and abandonments.
\newblock {\em Stoch. Syst.}, 4(1):250--299, 2014.

\bibitem{Afeche2021}
P.~Af{\`{e}}che, R.~Caldentey, and V.~Gupta.
\newblock On the optimal design of a bipartite matching queueing system.
\newblock {\em Oper. Res.}, 70(1):363--401, 2021.

\bibitem{Ananthanarayanan2013}
G.~Ananthanarayanan, A.~Ghodsi, S.~Shenker, and I.~Stoica.
\newblock {Effective straggler mitigation: Attack of the clones}.
\newblock In {\em Proc. 10th USENIX Conference on Networked Systems Design and
  Implementation}, pages 185--198, Lombard, IL, 2013. {USENIX}.

\bibitem{Anton2021}
E.~Anton, U.~Ayesta, M.~Jonckheere, and I.~M. Verloop.
\newblock On the stability of redundancy models.
\newblock {\em Oper. Res.}, 69(5):1540--1565, 2021.

\bibitem{Atar2019replicate}
R.~Atar, I.~Keslassy, and G.~Mendelson.
\newblock {Replicate to the shortest queues}.
\newblock {\em Queueing Syst.}, 92:1--23, 2019.

\bibitem{Atar2019}
R.~Atar, I.~Keslassy, and G.~Mendelson.
\newblock {Subdiffusive load balancing in time-varying queueing systems}.
\newblock {\em Oper. Res.}, 67(6):1678--1698, 2019.

\bibitem{ayesta2019token}
U.~Ayesta, T.~Bodas, J.~L. Dorsman, and I.~M. Verloop.
\newblock A token-based central queue with order-independent service rates.
\newblock {\em Oper. Res.}, 70(1):545--561, 2021.

\bibitem{ayesta2018unifying}
U.~Ayesta, T.~Bodas, and I.~M. Verloop.
\newblock On a unifying product form framework for redundancy models.
\newblock {\em Perform. Eval.}, 127:93--119, 2018.

\bibitem{Banerjee2020}
S.~Banerjee, Y.~Kanoria, and P.~Qian.
\newblock Dynamic assignment control of a closed queueing network under
  complete resource pooling, 2020.
\newblock {Preprint}.

\bibitem{Bell2001}
S.~L. Bell and R.~J. Williams.
\newblock {Dynamic scheduling of a system with two parallel servers in heavy
  traffic with resource pooling: Asymptotic optimality of a threshold policy}.
\newblock {\em Ann. Appl. Probability}, 11(3):608--649, 2001.

\bibitem{Bonald2017b}
T.~Bonald and C.~Comte.
\newblock {Balanced fair resource sharing in computer clusters}.
\newblock {\em Perform. Eval.}, 116:70--83, Nov. 2017.

\bibitem{Bonald2017}
T.~Bonald, C.~Comte, and F.~Mathieu.
\newblock {Performance of balanced fairness in resource pools: A recursive
  approach}.
\newblock {\em Proc. ACM Meas. Anal. Comput. Syst.}, 1(2):1--25, Dec. 2017.

\bibitem{Bramson1998}
M.~Bramson.
\newblock {State space collapse with application to heavy traffic limits for
  multiclass queueing networks}.
\newblock {\em Queueing Sys.}, 30(1-2):89--148, Nov. 1998.

\bibitem{Bramson2011}
M.~Bramson.
\newblock {Stability of join the shortest queue networks}.
\newblock {\em Ann. Appl. Probability}, 21(4):1568--1625, Aug. 2011.

\bibitem{Budhiraja2019}
A.~Budhiraja, D.~Mukherjee, and R.~Wu.
\newblock {Supermarket model on graphs}.
\newblock {\em Ann. Appl. Probab.}, 29(3):1740--1777, 2019.

\bibitem{Cardinaels2019}
E.~Cardinaels, S.~C. Borst, and J.~S.~H. van Leeuwaarden.
\newblock {Job assignment in large-scale service systems with affinity
  relations}.
\newblock {\em Queueing Syst.}, 93:227--268, 2019.

\bibitem{Chen2012}
H.~Chen and H.~Q. Ye.
\newblock {Asymptotic optimality of balanced routing}.
\newblock {\em Oper. Res.}, 60(1):163--179, Jan. 2012.

\bibitem{Cruise2020}
J.~Cruise, M.~Jonckheere, and S.~Shneer.
\newblock {Stability of JSQ in queues with general server-job class
  compatibilities}.
\newblock {\em Queueing Syst.}, 95(3):271--279, Jun. 2020.

\bibitem{Dai2008}
J.~G. Dai and W.~Lin.
\newblock Asymptotic optimality of maximum pressure policies in stochastic
  processing networks.
\newblock {\em Ann. Appl. Probab.}, 18(6):2239 -- 2299, 2008.

\bibitem{Eschenfeldt2016}
P.~Eschenfeldt and D.~Gamarnik.
\newblock Supermarket queueing system in the heavy traffic regime. {S}hort
  queue dynamics, 2018.
\newblock {Preprint}.

\bibitem{Feller1971}
W.~Feller.
\newblock {\em An Introduction to Probability Theory and its Applications},
  volume~2 of {\em Wiley Series in Probability and Statistics}.
\newblock Wiley, 2 edition, 1971.

\bibitem{FS99}
P.~J. Fleming and B.~Simon.
\newblock Heavy-traffic approximations for a system of infinite servers with
  load balancing.
\newblock {\em Prob. Engrg. Inform. Sci.}, 13(3):251--273, 1999.

\bibitem{Foss1998}
S.~G. Foss and N.~I. Chernova.
\newblock {On the stability of a partially accessible multi-station queue with
  state-dependent routing}.
\newblock {\em Queueing Syst.}, 29:55--73, 1998.

\bibitem{Gardner2017scheduling}
K.~Gardner, M.~Harchol-Balter, E.~Hyyti{\"{a}}, and R.~Righter.
\newblock {Scheduling for efficiency and fairness in systems with redundancy}.
\newblock {\em Perform. Eval.}, 116:1--25, Nov. 2017.

\bibitem{gardner2017redundancy}
K.~Gardner, M.~Harchol-Balter, A.~Scheller-Wolf, M.~Velednitsky, and
  S.~Zbarsky.
\newblock Redundancy-$d$: The power of $d$ choices for redundancy.
\newblock {\em Oper. Res.}, 65(4):1078--1094, 2017.

\bibitem{Gardner2020}
K.~Gardner and R.~Righter.
\newblock {Product forms for FCFS queueing models with arbitrary server-job
  compatibilities: an overview}.
\newblock {\em Queueing Syst.}, 96(1):3--51, Oct. 2020.

\bibitem{Gardner2016queueing}
K.~Gardner, S.~Zbarsky, S.~Doroudi, M.~Harchol-Balter, E.~Hyyti{\"a}, and
  A.~Scheller-Wolf.
\newblock Queueing with redundant requests: Exact analysis.
\newblock {\em Queueing Syst.}, 83(3-4):227--259, 2016.

\bibitem{gast2015power}
N.~Gast.
\newblock {The power of two choices on graphs: The pair-approximation is
  accurate?}
\newblock {\em ACM SIGMETRICS Perform. Eval. Rev.}, 43(2):69--71, 2015.

\bibitem{Graves2003}
S.~C. Graves and B.~T. Tomlin.
\newblock Process flexibility in supply chains.
\newblock {\em Manag. Sci.}, 49(7):907--919, 2003.

\bibitem{Harrison1998}
J.~M. Harrison.
\newblock {Heavy traffic analysis of a system with parallel servers: Asymptotic
  optimality of discrete-review policies}.
\newblock {\em Ann. Appl. Probab.}, 8(3):822--848, 1998.

\bibitem{Harrison1999}
J.~M. Harrison and M.~J. L{\'{o}}pez.
\newblock {Heavy traffic resource pooling in parallel server systems}.
\newblock {\em Queueing Syst.}, 33(4):339--368, 1999.

\bibitem{He2008}
Y.~T. He and D.~G. Down.
\newblock {Limited choice and locality considerations for load balancing}.
\newblock {\em Perform. Eval.}, 65(9):670--687, Aug. 2008.

\bibitem{Hellemans2019}
T.~Hellemans, T.~Bodas, and B.~{Van Houdt}.
\newblock {Performance analysis of workload dependent load balancing policies}.
\newblock {\em Proc. ACM Meas. Anal. Comput. Syst.}, 3(2):1--35, Jun. 2019.

\bibitem{Hellemans2018}
T.~Hellemans and B.~{Van Houdt}.
\newblock On the power-of-$d$-choices with least loaded server selection.
\newblock {\em Proc. ACM Meas. Anal. Comput. Syst.}, 2(2):1--22, Jun. 2018.

\bibitem{Hellemans2021}
T.~Hellemans and B.~Van~Houdt.
\newblock Mean waiting time in large-scale and critically loaded power of $d$
  load balancing systems.
\newblock {\em Proc. ACM Meas. Anal. Comput. Syst.}, 5(2), 2021.

\bibitem{hurtado2020transform}
D.~Hurtado-Lange and S.~T. Maguluri.
\newblock Transform methods for heavy-traffic analysis.
\newblock {\em Stoch. Syst.}, 10(4):275--309, 2020.

\bibitem{HurtadoLange2020}
D.~Hurtado-Lange and S.~T. Maguluri.
\newblock Throughput and delay optimality of power-of-$d$ choices in
  inhomogeneous load balancing systems.
\newblock {\em Oper. Res. Lett.}, 49(4):616--622, jul 2021.

\bibitem{Jordan1995}
W.~C. Jordan and S.~C. Graves.
\newblock Principles on the benefits of manufacturing process flexibility.
\newblock {\em Manag. Sci.}, 41(4):577--594, 1995.

\bibitem{Joshi2018}
G.~Joshi.
\newblock Synergy via redundancy: Boosting service capacity with adaptive
  replication.
\newblock {\em ACM SIGMETRICS Perform. Eval. Rev.}, 45(3):21–28, Mar. 2018.

\bibitem{Joshi2016}
G.~Joshi, E.~Soljanin, and G.~Wornell.
\newblock Efficient redundancy techniques for latency reduction in cloud
  systems.
\newblock {\em ACM Trans. Model. Perform. Eval. Comput. Syst.}, 2(2), Apr.
  2017.

\bibitem{Keilson1988}
J.~Keilson and L.~D. Servi.
\newblock {A distributional form of Little's law}.
\newblock {\em Oper. Res. Lett.}, 7(5):223--227, Oct. 1988.

\bibitem{Kelly1993}
F.~P. Kelly and C.~N. Laws.
\newblock Dynamic routing in open queueing networks: Brownian models, cut
  constraints and resource pooling.
\newblock {\em Queueing Syst.}, 13:47--86, 1993.

\bibitem{Krzesinski2011}
A.~E. Krzesinski.
\newblock {Order independent queues}.
\newblock In R.~Boucherie and N.~van Dijk, editors, {\em Queueing Networks.
  International Series in Operations Research and Management Science}, volume
  154, pages 85--120. Springer, Boston, MA, 2011.

\bibitem{Laws1992}
C.~N. Laws.
\newblock Resource pooling in queueing networks with dynamic routing.
\newblock {\em Adv. in Appl. Probab.}, 24(3):699--726, 1992.

\bibitem{Lee2017}
K.~Lee, R.~Pedarsani, and K.~Ramchandran.
\newblock {On scheduling redundant requests with cancellation overheads}.
\newblock {\em IEEE/ACM Trans. Netw.}, 25(2):1279--1290, Apr. 2017.

\bibitem{Maguluri2015}
S.~T. Maguluri and R.~Srikant.
\newblock {Heavy-traffic behavior of the Maxweight algorithm in a switch with
  uniform traffic}.
\newblock {\em ACM SIGMETRICS Perform. Eval. Rev.}, 43(2):72--74, Sept. 2015.

\bibitem{Mandelbaum2004}
A.~Mandelbaum and A.~L. Stolyar.
\newblock Scheduling flexible servers with convex delay costs: Heavy-traffic
  optimality of the generalized $c\mu$-rule.
\newblock {\em Oper. Res.}, 52(6):836--855, 2004.

\bibitem{mitzenmacher2001power}
M.~Mitzenmacher.
\newblock The power of two choices in randomized load balancing.
\newblock {\em IEEE Trans. Parallel Distrib. Syst.}, 12(10):1094--1104, 2001.

\bibitem{Mukherjee2018}
D.~Mukherjee, S.~C. Borst, and J.~S.~H. van Leeuwaarden.
\newblock Asymptotically optimal load balancing topologies.
\newblock {\em Proc. ACM Meas. Anal. Comput. Syst.}, 2(1), Apr. 2018.

\bibitem{Rutten2020}
D.~Rutten and D.~Mukherjee.
\newblock Load balancing under strict compatibility constraints.
\newblock {\em Math. Oper. Res.}, 2021.
\newblock {Forthcoming}.

\bibitem{Shah2012}
D.~Shah and D.~Wischik.
\newblock {Switched networks with maximum weight policies: Fluid approximation
  and multiplicative state space collapse}.
\newblock {\em Ann. Appl. Probab.}, 22(1):70--127, Feb. 2012.

\bibitem{Shah2016}
N.~B. Shah, K.~Lee, and K.~Ramchandran.
\newblock {When do redundant requests reduce latency?}
\newblock {\em IEEE Trans. Comm.}, 64(2):715--722, 2016.

\bibitem{Sharifnassab2020}
A.~Sharifnassab, J.~N. Tsitsiklis, and S.~J. Golestani.
\newblock {Fluctuation bounds for the Max-weight policy with applications to
  state space collapse}.
\newblock {\em Stoch. Syst.}, 10(3):223--250, Apr. 2020.

\bibitem{Shi2019}
C.~Shi, Y.~Wei, and Y.~Zhong.
\newblock Process flexibility for multiperiod production systems.
\newblock {\em Oper. Res.}, 67(5):1300--1320, 2019.

\bibitem{Simchi-Levi2012}
D.~Simchi-Levi and Y.~Wei.
\newblock Understanding the performance of the long chain and sparse designs in
  process flexibility.
\newblock {\em Oper. Res.}, 60(5):1125--1141, 2012.

\bibitem{Sloothaak2019}
F.~Sloothaak, J.~Cruise, S.~Shneer, M.~Vlasiou, and B.~Zwart.
\newblock Complete resource pooling of a load-balancing policy for a network of
  battery swapping stations.
\newblock {\em Queueing Syst.}, 99:65--120, 2021.

\bibitem{Stolyar2004}
A.~L. Stolyar.
\newblock {MaxWeight scheduling in a generalized switch: State space collapse
  and workload minimization in heavy traffic}.
\newblock {\em Ann. Appl. Probab.}, 14(1):1--53, 2004.

\bibitem{Stolyar2005}
A.~L. Stolyar.
\newblock Optimal routing in output-queued flexible server systems.
\newblock {\em Prob. Engrg. Inform. Sci.}, 19:141--189, 2005.

\bibitem{TX13a}
J.~N. Tsitsiklis and K.~Xu.
\newblock On the power of (even a little) resource pooling.
\newblock {\em Stoch. Syst.}, 2(1):1--66, 2012.

\bibitem{TX13b}
J.~N. Tsitsiklis and K.~Xu.
\newblock Queueing system topologies with limited flexibility.
\newblock SIGMETRICS '13, pages 167--178, New York, NY, USA, 2013. AMC.

\bibitem{TX17}
J.~N. Tsitsiklis and K.~Xu.
\newblock Flexible queueing architectures.
\newblock {\em Oper. Res.}, 5(65):1398--1413, 2017.

\bibitem{Turner1998}
S.~R.~E. Turner.
\newblock {The effect of increasing routing choice on resource pooling}.
\newblock {\em Prob. Engrg. Inform. Sci.}, 12:109--124, 1998.

\bibitem{Varma2021}
S.~M. Varma and S.~T. Maguluri.
\newblock Transportation polytope and its applications in parallel server
  systems, 2021.
\newblock {Preprint}.

\bibitem{visschers2012product}
J.~Visschers, I.~Adan, and G.~Weiss.
\newblock A product form solution to a system with multi-type jobs and
  multi-type servers.
\newblock {\em Queueing Syst.}, 70(3):269--298, 2012.

\bibitem{vvedenskaya1996queueing}
N.~D. Vvedenskaya, R.~L. Dobrushin, and F.~I. Karpelevich.
\newblock Queueing system with selection of the shortest of two queues: An
  asymptotic approach.
\newblock {\em Problemy Peredachi Informatsii}, 32(1):20--34, 1996.

\bibitem{Weng2020}
W.~Weng, X.~Zhou, and R.~Srikant.
\newblock Optimal load balancing with locality constraints.
\newblock {\em Proc. ACM Meas. Anal. Comput. Syst.}, 4(3), Nov. 2020.

\bibitem{Zheng2021}
L.~Zheng, S.~Chandratre, A.~Ali, A.~Szabo, L.~Durham, L.~D. Joyce, and D.~L.
  Joyce.
\newblock How does multiple listing affect lung transplantation? a
  retrospective analysis.
\newblock {\em Semin. Thorac Cardiovasc. Surg.}, may 2021.

\end{thebibliography}
\end{spacing}

\appendix

\section{Local stability conditions and CRP conditions}\label{app:CRP}
In this appendix we formulate two alternative CRP conditions, involving a linear programming characterization and a geometric representation respectively, and prove that these are equivalent to the local stability conditions~\eqref{stabcond3} and~\eqref{stabcond4}.

We first introduce some useful definitions and preliminaries.
For any $\mathcal{T}\subseteq\mathcal{S}$, define
\[
f_{\mathcal{T}}(\boldsymbol{\lambda}) \coloneqq \sum_{S \in \mathcal{T}} \lambda_S - \mu_{\mathcal{T}} = \boldsymbol{w}_{\mathcal{T}} \cdot \boldsymbol{\lambda} - \mu_{\mathcal{T}},
\]
with $\boldsymbol{\lambda} = (\lambda_S)_{S \in \mathcal{S}}$ and $\boldsymbol{w}_{\mathcal{T}}$ a 0--1 vector of length $|\mathcal{S}|$ indicating if $S \notin \mathcal{T}$ or $S \in \mathcal{T}$, respectively.
The stability conditions in~\eqref{stabcond1} for both the c.o.c.\ and c.o.s.\ strategies may be equivalently written as
\[
f_{\mathcal{T}}(\boldsymbol{\lambda}) < 0 \hspace*{.4in} \text{ for all } \mathcal{T} \subseteq \mathcal{S},
\]
and the closure of the stability region may thus be defined as
\[
\mathcal{C} \coloneqq \left\{\boldsymbol{\lambda} \in \mathbb{R}_+^{|\mathcal{S}|} \colon f_{\mathcal{T}}(\boldsymbol{\lambda}) \le 0 \text{ for all } \mathcal{T} \subseteq \mathcal{S}\right\},
\]
with Pareto boundary
\[
\partial\mathcal{C} \coloneqq \left\{\boldsymbol{\lambda} \in \mathbb{R}_+^{|\mathcal{S}|} \colon f_{\mathcal{T}}(\boldsymbol{\lambda}) \le 0 \mbox{ for all } \mathcal{T} \subseteq \mathcal{S} \text{~and~} f_{\mathcal{T}'}(\boldsymbol{\lambda}) = 0 \text{ for some } \mathcal{T}' \subseteq \mathcal{S}\right\}.
\]
Each `face' of the boundary thus corresponds to a hyperplane of the form $\boldsymbol{w}_{\mathcal{T}} \cdot \boldsymbol{\lambda} - \mu_{\mathcal{T}} = 0$ with variables $\boldsymbol{\lambda}$, offset $- \mu_{\mathcal{T}}$ and normal vector $\boldsymbol{w}_{\mathcal{T}}$.

We further define the `rate region' of the system as
\begin{equation}
\mathcal{R} \coloneqq \{\boldsymbol{r} \in \mathbb{R}_+^{|\mathcal{S}|}\colon r_S = \sum_{n \in S} \mu_n x_{nS} \mbox{ for some } x \in \mathcal{X}\},
\label{rate1}
\end{equation}
with
\[
\mathcal{X} \coloneqq
\{(x_{nS})^{\mathcal{I}}\colon \sum_{S : n\in S} x_{nS} \leq 1 \mbox{ for all } n = 1, \dots, N \mbox{ and } x_{nS} \geq 0 \mbox{ for all } (n, S) \in \mathcal{I}\},
\]
and
\[
\mathcal{I} = \bigcup_{S \in \mathcal{S}} \bigcup_{n \in S} \{(n, S)\}
\]
denoting the set of all compatible pairs of server and job types.
Interpreting $x_{nS}$ as the fraction of time that server~$n$ spends serving type-$S$ jobs, it is seen that $\mathcal{R}$ indeed represents the set of all achievable service rates for the various job types.

It is easily verified that the rate region may equivalently be expressed as
\begin{equation}
\mathcal{R} = \{\boldsymbol{\lambda} \in \mathbb{R}_+^{|\mathcal{S}|}\colon \sum_{S \in \mathcal{S}} \lambda_S u_{nS} \leq \mu_n \mbox{ for some } u \in \mathcal{U}\},
\label{rate2}
\end{equation}
with
\[
\mathcal{U} \coloneqq
\{(u_{nS})^{\mathcal{I}}\colon \sum_{n \in S} u_{nS} = 1 \mbox{ for all } S \in \mathcal{S} \mbox{ and } u_{nS} \geq 0 \mbox{ for all } (n, S) \in \mathcal{I}\}.
\]
Interpreting $u_{nS}$ as the fraction of type-$S$ jobs handled by server~$n$, we note that $\mathcal{R}$ as defined in~\eqref{rate2} represents the set of all arrival rates of the various job types that can be supported without overloading any of the servers.

Using similar arguments as in~\cite{Cardinaels2019}, it can also be shown that the set $\mathcal{R}$ coincides with the closure of the stability region~$\mathcal{C}$ .
Henceforth, we will use these two representations interchangeably and refer to both $\mathcal{C}$ and $\mathcal{R}$ as the `capacity' region.
Let $\boldsymbol{r}^*$ be the point on the boundary of the capacity region that is approached by the vector $N \lambda (p_S)_{S \in \mathcal{S}}$ in heavy traffic as $\lambda \uparrow \lambda^*$.

\mathtoolsset{showonlyrefs=false}
Consider the linear program
\begin{eqnarray}
\min & & u \\
\mbox{sub} & & \sum_{n \in S} \mu_n x_{nS} \geq N \lambda p_S
\hspace*{.4in} \mbox{ for all } S \in \mathcal{S} \\
& & \sum_{S \ni n} x_{nS} \leq u \hspace*{.4in} \mbox{ for all } n = 1, \dots, N \\
& & x_{nS} \geq 0 \hspace*{.8in} \mbox{ for all } (n, S) \in \mathcal{I}.
\end{eqnarray}
\mathtoolsset{showonlyrefs=true}
Let $u^*$ denote the optimal value of the objective function. It is easily seen that $u^* \leq 1$ iff $N \lambda (p_S)_{S \in \mathcal{S}}$ belongs to the capacity region, and in particular $u^* = 1$ corresponds to a `heavy-traffic' regime.

\mathtoolsset{showonlyrefs=false}
The dual version of the above linear program is
\begin{eqnarray}
\max & & N \lambda \sum_{s \in \mathcal{S}} p_S y_S \label{dualobje1} \\
\mbox{sub} & & \mu_n y_S \leq z_n \hspace*{.4in} \mbox{ for all } (n, S) \in \mathcal{I} \label{dualcons1} \\
& & \sum_{n = 1}^{N} z_n \leq 1 \label{dualcons2} \\
& & z_n \geq 0 \hspace*{.4in} \mbox{ for all } n = 1, \dots, N. \label{dualcons3}
\end{eqnarray}
\mathtoolsset{showonlyrefs=true}

\begin{theorem}[Equivalence among various CRP conditions] \label{th:equiv_CRP}
The following statements are equivalent:
\begin{itemize}
\item[(i)] For all (non-empty) $\mathcal{T} \subsetneq \mathcal{S}$ there holds that $p_{\mathcal{T}} < \frac{1}{N \mu} \mu_{\mathcal{T}}$, i.e., the local stability conditions~\eqref{stabcond3} and~\eqref{stabcond4} are satisfied.
\item[(ii)] For all (non-empty) $\mathcal{T} \subsetneq \mathcal{S}$ there holds that $f_{\mathcal{T}}(\boldsymbol{r}^*) < 0$ and $f_{\mathcal{S}}(\boldsymbol{r}^*) = 0$.
\item[(iii)] There holds that $\lambda^* = \mu$, i.e., $\boldsymbol{r}^* = N \mu (p_S)_{S \in \mathcal{S}}$, and the normal vector to the boundary of the capacity region at~$\boldsymbol{r}^*$ is unique up to a scalar coefficient, and all its components are strictly positive, i.e., is given by $\boldsymbol{w}_{\mathcal{S}} = (1, 1, \dots, 1)$ up to a scalar coefficient.
\item[(iv)] The optimal solution of the dual linear program \eqref{dualobje1}--\eqref{dualcons3} for $\lambda = \mu$ is unique and given by
\[
\boldsymbol{y}^* = \frac{1}{N \mu} (1, 1, \dots, 1) \mbox{ and } \boldsymbol{z}^* = \frac{1}{N \mu} (\mu_1, \mu_2, \dots, \mu_N).
\]
\end{itemize}
\end{theorem}

\begin{remark} \normalfont
It is worth observing that in~\cite{Laws1992} it is mentioned that saturation of a single `generalized cut constraint' as a more general notion of statement~(i) implies uniqueness of the solution of a dual version of a linear program associated with a multi-commodity flow problem, closely mirroring statement~(iv).
Also, Stolyar~\cite{Stolyar2004} uses a Lagrangian formulation to prove equivalence between statements~(iii) and~(iv) in a more general setting, albeit for a somewhat different incarnation of the dual program.
\end{remark}

\begin{proof}[Proof of Theorem~\ref{th:equiv_CRP}]
We will successively prove the following implications: (i) $\Rightarrow$ (iii); (iii) $\Rightarrow$ (ii); (ii) $\Rightarrow$ (i); (iv) $\Rightarrow$ (i); and (i) $\Rightarrow$ (iv). \\

(i) $\Rightarrow$ (iii):
Since $\boldsymbol{r}^*$ belongs to~$\mathcal{C}$, it must in particular satisfy the global constraint $f_{\mathcal{S}}(\boldsymbol{r}^*) \le 0$, hence $N \lambda^* = N \lambda^* p_{\mathcal{S}} \le \mu_{\mathcal{S}} = N \mu$, i.e., $\lambda^* \le \mu$.
Thus, by virtue of~(i), for all (non-empty) ${\mathcal T} \subsetneq {\mathcal S}$,
\[
f_{\mathcal{T}}(\boldsymbol{r}^*) =
\boldsymbol{w}_{\mathcal{T}} \cdot \boldsymbol{r}^* - \mu_{\mathcal{T}} =
N \lambda^* p_{\mathcal{T}} - \mu_{\mathcal{T}} \le
N \mu p_{\mathcal{T}} - \mu_{\mathcal{T}} =
N \mu \left(p_{\mathcal{T}} - \frac{\mu_{\mathcal{T}}}{N \mu}\right) < 0,
\]
i.e., $\boldsymbol{r}^*$ does not belong to a face corresponding to any of the local inequalities.
Because $\boldsymbol{r}^*$ is an element of the boundary $\partial\mathcal{C}$, it follows that it must satisfy the global constraint with equality, $f_{\mathcal{S}}(\boldsymbol{r}^*) = 0$, implying that in fact $\lambda^* = \mu$,
and the normal vector to the boundary $\partial\mathcal{C}$ at $\boldsymbol{r}^* = 0$ is unique and given by $\boldsymbol{w}_{\mathcal{S}} = (1, 1, \dots, 1)$ up to a scalar coefficient. \\

(iii) $\Rightarrow$ (ii):
In view of~(iii), since $\boldsymbol{w}_{\mathcal{S}}$ is the normal vector to the face corresponding to the global constraint, $\boldsymbol{r}^*$ must satisfy this constraint with equality, i.e., $f_{\mathcal{S}}(\boldsymbol{r}^*) = 0$.
Also, because $\boldsymbol{r}^*$ belongs to~$\mathcal{C}$ by definition, it must satisfy all the local constraints, but the uniqueness of the normal vector rules out that any of these are satisfied with equality, i.e., $f_{\mathcal{T}}(\boldsymbol{r}^*) < 0$ for all $\mathcal{T} \subsetneq \mathcal{S}$. \\

(ii) $\Rightarrow$ (i):
Noting that
\[
f_{\mathcal{S}}(\boldsymbol{r}^*) = \boldsymbol{w}_{\mathcal{S}} \cdot \boldsymbol{r}^* - \mu_{\mathcal{S}} = N \lambda^* - N \mu = 0,
\]
we deduce that $\lambda^* = \mu$.
Consequently, for any $\mathcal{T} \subsetneq \mathcal{S}$, we have
\[
0 > f_{\mathcal{T}}(\boldsymbol{r}^*) =
\boldsymbol{w}_{\mathcal{T}} \cdot \boldsymbol{r}^* - \mu_{\mathcal{T}} =
N \lambda^* p_{\mathcal{T}} - \mu_{\mathcal{T}} =
N \mu p_{\mathcal{T}} - \mu_{\mathcal{T}},
\]
i.e., $p_{\mathcal{T}} < \frac{1}{N \mu} \mu_{\mathcal{T}}$. \\

(iv) $\Longrightarrow$ (i):
Suppose that (iv) holds, but (i) does not, so there exists a (non-empty) $\mathcal{T} \subsetneq \mathcal{S}$ for which $p_{\mathcal{T}} \geq \frac{1}{N \mu} \mu_{\mathcal{T}}$, i.e.,
\begin{equation}
N \mu \sum_{S \in \mathcal{T}} p_S \geq \sum_{n \in \mathcal{N}(\mathcal{T})} \mu_n,
\label{eq1}
\end{equation}
with $\mathcal{N}(\mathcal{T}) = \bigcup_{S \in \mathcal{T}} S$.

We will construct a feasible solution $(\hat{\boldsymbol{y}}, \hat{\boldsymbol{z}})$ of the dual linear program with the same objective value as $(\boldsymbol{y}^*, \boldsymbol{z}^*)$, which yields a contradiction with (iv) and hence shows that (i) must hold.
In order to do so, first of all observe that the optimal value of the dual objective function for $\lambda = \mu$ equals
$N \lambda \sum_{S \in \mathcal{S}} p_S y_S^* = N \mu \sum_{S \in \mathcal{S}} p_S \frac{1}{N \mu} = 1$.
Hence the optimal value of the primal objective function also equals $u^* = 1$.
Thus, the optimal solution $({u}^*, \boldsymbol{x}^*)$ of the primal linear program for $\lambda = \mu$ satisfies
\begin{equation}
\sum_{n \in S} \mu_n x_{nS}^* \geq N \mu p_S
\hspace*{.4in} \mbox{ for all } S \in \mathcal{S},
\label{ineq1}
\end{equation}
and
\begin{equation}
\sum_{S: n\in S} x_{nS}^* \leq u^* = 1 \hspace*{.4in} \mbox{ for all } n = 1, \dots, N.
\label{ineq2}
\end{equation}
Summing the inequalities in~\eqref{ineq1} over all job types $S \in \mathcal{T}$, we obtain.
\begin{equation}
\sum_{S \in \mathcal{T}} \sum_{n \in S} \mu_n x_{nS}^* \geq N \mu \sum_{S \in \mathcal{T}} p_S.
\label{ineq3}
\end{equation}
For convenience, define
\[
\mathcal{I}^0 = \{(n, S) \in \mathcal{I}\colon n \in \mathcal{N}(\mathcal{T}), S \notin \mathcal{T}\}
\]
 as the set of all pairs of servers in $\mathcal{N}(\mathcal{T})$ and compatible job types that do not belong to~$\mathcal{T}$.
Then, using the inequalities in~\eqref{ineq2} and invoking~\eqref{eq1}, we derive
\mathtoolsset{showonlyrefs=false}
\begin{eqnarray}
\sum_{S \in \mathcal{T}} \sum_{n \in S} \mu_n x_{nS}^*
&=&
\sum_{n \in \mathcal{N}(\mathcal{T})} \mu_n \sum_{S : n\in S\in\mathcal{T}} x_{nS}^* \nonumber \\
&=&
\sum_{n \in \mathcal{N}(\mathcal{T})} \mu_n \left[\sum_{S: n\in S} x_{nS}^* - \sum_{S : n\in S\notin\mathcal{T}} x_{nS}^*\right] \nonumber \\
&\leq&
\sum_{n \in \mathcal{N}(\mathcal{T})} \mu_n \left[1 - \sum_{S : n\in S\notin\mathcal{T}} x_{nS}^*\right] \nonumber \\
&=&
\sum_{n \in \mathcal{N}(\mathcal{T})} \mu_n - \sum_{n \in \mathcal{N}(\mathcal{T})} \mu_n \sum_{S : n\in S\notin\mathcal{T}} x_{nS}^* \nonumber \\
&\leq&
N \mu \sum_{S \in \mathcal{T}} p_{\mathcal{S}} - \sum_{(n, S) \in \mathcal{I}^0} \mu_n x_{nS}^*.
\label{ineq4}
\end{eqnarray}
\mathtoolsset{showonlyrefs=true}
Combining the inequalities in~\eqref{ineq3} and~\eqref{ineq4}, we deduce
\begin{equation}
x_{nS}^* = 0 \hspace*{.4in} \mbox{ for all } (n, S) \in \mathcal{I}^0.
\label{zero}
\end{equation}

Now define $(\hat{\boldsymbol{y}}$, $\hat{\boldsymbol{z}})$ by
\[
\hat{y}_S = \left\{\begin{array}{ll} y_S^* (1 + C^+) & \text{~if~}S \in \mathcal{T} \\ y_S^* (1 - C^-) & \text{~if~}S \notin \mathcal{T} \end{array}\right.
\]
and
\[
\hat{z}_n = \left\{\begin{array}{ll} z_n^* (1 + C^+) & \text{~if~}n \in \mathcal{N}(\mathcal{T}) \\ z_n^* (1 - C^-) & \text{~if~}n \notin \mathcal{N}(\mathcal{T}) \end{array}\right.,
\]
with
\begin{equation}
C^+ = \frac{\epsilon}{\mu_{\mathcal{T}}}, \hspace{.4in}
C^- = \frac{\epsilon}{N \mu - \mu_{\mathcal{T}}},
\end{equation}
and $\epsilon \in [0, \epsilon_{\max}]$, with $\epsilon_{\max} = N \mu - \mu_{\mathcal{T}} > 0$.

Note that
\begin{eqnarray*}
\sum_{n = 1}^{N} \hat{z}_n
&=&
\sum_{n \in \mathcal{N}(\mathcal{T})} \hat{z}_n + \sum_{n \notin \mathcal{N}(\mathcal{T})} \hat{z}_n \\
&=&
\sum_{n \in \mathcal{N}(\mathcal{T})} z_n^* (1 + C^+) + \sum_{n \notin \mathcal{N}(\mathcal{T})} z_n^* (1 - C^-) \\
&=&
\sum_{n \in \mathcal{N}(\mathcal{T})} z_n^* + \sum_{n \notin \mathcal{N}(\mathcal{T})} z_n^* + C^+ \sum_{n \in \mathcal{N}(\mathcal{T})} \mu_n - C^- \sum_{n \notin \mathcal{N}(\mathcal{T})} \mu_n \\
&=&
\sum_{n = 1}^{N} z_n^* = 1,
\end{eqnarray*}
implying that the dual constraint~\eqref{dualcons1} is satisfied.

Further observe that the index set $\mathcal{I}$ may be partitioned into three disjoint subsets $\mathcal{I}^+$, $\mathcal{I}^0$ and $\mathcal{I}^-$, with
$\mathcal{I}^+ = \{(n, S) \in \mathcal{I}\colon n \in \mathcal{N}(\mathcal{T}),S \in \mathcal{T}\}$, $\mathcal{I}^0$ as defined above and $\mathcal{I}^- = \{(n, S) \in \mathcal{I}\colon n \notin \mathcal{N}(\mathcal{T}), S\notin\mathcal{T}\}$.
By construction,
\[
(\hat{y}_S, \hat{z}_n) = \left\{\begin{array}{ll}
( y_S^*(1+C^+), z_n^*(1+C^+)) & \mbox{ for all } (n, S) \in \mathcal{I}^+ \\
( y_S^*(1-C^-),  z_n^*(1+C^+)) & \mbox{ for all } (n, S) \in \mathcal{I}^0 \\
( y_S^*(1-C^-),  z_n^*(1-C^-)) & \mbox{ for all } (n, S) \in \mathcal{I}^-
\end{array} \right. .
\]
Since the dual constraints~\eqref{dualcons1} are all satisfied with equality for $(\boldsymbol{y}^*, \boldsymbol{z}^*)$, it follows that
\[
\mu_n \hat{y}_S = \hat{z}_n \hspace*{.4in} \mbox{ for all } (n, S) \in \mathcal{I}^+ \cup \mathcal{I}^-,
\]
while
\[
\mu_n \hat{y}_S \leq \hat{z}_n \hspace*{.4in} \mbox{ for all } (n, S) \in \mathcal{I}^0
\]
for any $\epsilon \geq 0$.
Thus, $(\hat{\boldsymbol{y}}, \hat{\boldsymbol{z}})$ is a feasible solution of the dual linear program, and $(u^*, \boldsymbol{x}^*)$ and $(\hat{\boldsymbol{y}}, \hat{\boldsymbol{z}})$ satisfy the complementary slackness conditions in view of~\eqref{zero}.
This implies that $(\hat{\boldsymbol{y}}, \hat{\boldsymbol{z}}) \neq (\boldsymbol{y}^*, \boldsymbol{z}^*)$ is also an optimal solution of the dual linear program for any $\epsilon \in (0, \epsilon_{\max}]$, yielding a contradiction with the initial supposition. \\

(i) $\Longrightarrow$ (iv):

Note that $(\boldsymbol{y}^*, \boldsymbol{z}^*)$ is a feasible solution of the dual linear program for $\lambda = \mu$, with objective value $N \lambda \sum_{s \in \mathcal{S}} p_S y_S^* = N \mu \sum_{S \in \mathcal{S}} p_S \frac{1}{N \mu} = 1$.
Let $(u^*, \boldsymbol{x}^*)$ be an optimal solution of the primal linear program for $\lambda = \mu$.
Since (i) implies (iii), we know that $\boldsymbol{r}^* = N \mu (p_S)_{S \in \mathcal{S}}$ belongs to the capacity region, yielding $u^* \leq 1$.
Because the maximum objective value of the dual problem is bounded from above by the minimum value of the primal problem by weak duality, it follows that $(\boldsymbol{y}^*, \boldsymbol{z}^*)$ is in fact an optimal solution of the dual linear program.

In order to prove uniqueness, let $(\tilde{\boldsymbol{y}}, \tilde{\boldsymbol{z}})$ be an optimal solution of the dual linear program, and
\[
\mathcal{T} = \arg\max_{S \in \mathcal{S}} \tilde{y}_S = \{S \in \mathcal{S}: \tilde{y}_S = y_{\max}\},
\]
with $y_{\max} = \max_{S \in \mathcal{S}} \tilde{y}_S$.
We only need to consider the case that $\mathcal{T}$ is a strict subset of $\mathcal{S}$ because otherwise all the components of $\tilde{\boldsymbol{y}}$ are equal, and $\tilde{\boldsymbol{y}}$ must be identical to $\boldsymbol{y}^*$ in order to be optimal.

Now define
\[
\hat{y}_S = \left\{\begin{array}{ll} \tilde{y}_S - C^- = y_{\max} - C^- & \text{~if~}S \in \mathcal{T} \\
\tilde{y}_S + C^+ & \text{~if~}S \notin \mathcal{T} \end{array} \right.
\]
and
\[
\hat{z}_n = \left\{\begin{array}{ll} \tilde{z}_n - C^- \mu_n & \text{~if~}n \in \mathcal{N}(\mathcal{T}) \\
\tilde{z}_n + C^+ \mu_n & \text{~if~}n \notin \mathcal{N}(\mathcal{T}) \end{array} \right.,
\]
with $C^-$, $C^+$ and $\mathcal{N}(\mathcal{T}) = \bigcup_{S \in \mathcal{T}} S$ as defined earlier.

Note that
\begin{eqnarray*}
\sum_{n = 1}^{N} \hat{z}_n
&=&
\sum_{n \in \mathcal{N}(\mathcal{T})} \hat{z}_n + \sum_{n \notin \mathcal{N}(\mathcal{T})} \hat{z}_n \\
&=&
\sum_{n \in \mathcal{N}(\mathcal{T})} (\tilde{z}_n - C^- \mu_n) + \sum_{n \notin \mathcal{N}(\mathcal{T})} (\tilde{z}_n + C^+ \mu_n) \\
&=&
\sum_{n \in \mathcal{N}(\mathcal{T})} \tilde{z}_n + \sum_{n \notin \mathcal{N}(\mathcal{T})} \tilde{z}_n - C^- \sum_{n \in \mathcal{N}(\mathcal{T})} \mu_n + C^+ \sum_{n \notin \mathcal{N}(\mathcal{T})} \mu_n \\
&=&
\sum_{n = 1}^{N} \tilde{z}_n \leq 1,
\end{eqnarray*}
implying that the dual constraint~\eqref{dualcons2} is satisfied.

As before, the index set $\mathcal{I}$ may be partitioned into the three disjoint subsets $\mathcal{I}^+$, $\mathcal{I}^0$ and $\mathcal{I}^-$.

By construction,
\[
(\hat{y}_S, \hat{z}_n) = \left\{\begin{array}{ll}
(\tilde{y}_S - C^-, \tilde{z}_n - C^- \mu_n) = (y_{\max} - C^-, \tilde{z}_n - C^- \mu_n) & \mbox{ for all } (n, S) \in \mathcal{I}^+ \\
(\tilde{y}_S + C^-, \tilde{z}_n - C^+ \mu_n) & \mbox{ for all } (n, S) \in \mathcal{I}^0 \\
(\tilde{y}_S + C^+, \tilde{z}_n + C^+ \mu_n) & \mbox{ for all } (n, S) \in \mathcal{I}^-
\end{array} \right. .
\]

Since the dual constraints~\eqref{dualcons1} are all satisfied for $(\tilde{\boldsymbol{y}}, \tilde{\boldsymbol{z}})$, it follows that
\[
\mu_n \hat{y}_S \leq \hat{z}_n \hspace*{.4in} \mbox{ for all } (n, S) \in \mathcal{I}^+ \cup \mathcal{I}^-.
\]

Also, for all $n \in \mathcal{N}(\mathcal{T})$,
\[
\tilde{z}_n \geq \max_{S \ni n} \mu_n \tilde{y}_S = \mu_n y_{\max},
\]
and thus, defining $\Delta = y_{\max} - \max_{S \notin \mathcal{T}} \tilde{y}_S > 0$,
\begin{eqnarray*}
\mu_n \hat{y}_S
&=&
\mu_n (\tilde{y}_S + C^+) \\
&=&
\mu_n y_{\max} + \mu_n (\tilde{y}_S - y_{\max} + C^+) \\
&\leq&
\tilde{z}_n + \mu_n (\tilde{y}_S - y_{\max} + C^+) \\
&=&
\hat{z}_n + \mu_n (\tilde{y}_S - y_{\max} + C^+ + C^-) \\
&=&
\hat{z}_n + \mu_n (C^+ + C^- - \Delta) \leq \hat{z}_n
\end{eqnarray*}
for all $\Delta \geq C^+ + C^-$, i.e., $\epsilon \leq \Delta \left[\frac{1}{\mu_{\mathcal{T}}} + \frac{1}{N \mu - \mu_{\mathcal{T}}}\right]^{- 1}$.

Thus, the dual constraints~\eqref{dualcons1} are satisfied for $(\hat{\boldsymbol{y}}, \hat{\boldsymbol{z}})$ as well, so that $(\hat{\boldsymbol{y}}, \hat{\boldsymbol{z}})$ is a feasible solution of the dual linear program.

Also, \mathtoolsset{showonlyrefs=false}
\begin{eqnarray}
N \mu \sum_{S \in \mathcal{S}} p_S \hat{y}_S
&=&
N \mu \sum_{S \in \mathcal{T}} p_S \hat{y}_S + N \mu \sum_{S \notin \mathcal{T}} p_S \hat{y}_S \nonumber \\
&=&
N \mu \sum_{S \in \mathcal{T}} p_S (\tilde{y}_S - C^-) + N \mu \sum_{S \notin \mathcal{T}} p_S (\tilde{y}_S + C^+) \nonumber \\
&=&
N \mu \sum_{S \in \mathcal{T}} p_S \tilde{y}_S + N \mu \sum_{S \notin \mathcal{T}} p_S \tilde{y}_S - N \mu \sum_{S \in \mathcal{T}} p_S C^- + N \mu \sum_{S \notin \mathcal{T}} p_S C^+ \nonumber \\
&=&
N \mu \sum_{S \in \mathcal{S}} p_S \tilde{y}_S - N \mu \left[C^- \sum_{S \in \mathcal{T}} p_S - C^+ \sum_{S \notin \mathcal{T}} p_S \right]. \label{nega1}
\end{eqnarray}
\mathtoolsset{showonlyrefs=true}

Now observe that $p_{\mathcal{T}} < \frac{1}{N \mu} \mu_{\mathcal{T}}$ by virtue of~(i), so
\[
C^- \sum_{S \in \mathcal{T}} p_S = \frac{p_{\mathcal{T}}}{\mu_{\mathcal{T}}} < \frac{1}{N \mu},
\]
and
\[
\sum_{S \notin \mathcal{T}} p_S = 1 - p_{\mathcal{T}} > 1 - \frac{1}{N \mu} \mu_{\mathcal{T}} = \frac{N \mu - \mu_{\mathcal{T}}}{N \mu},
\]
so
\[
C^+ \sum_{S \notin \mathcal{T}} p_S > \frac{1}{N \mu - \mu_{\mathcal{T}}} \frac{N \mu - \mu_{\mathcal{T}}}{N \mu} = \frac{1}{N \mu}.
\]

We conclude that the term in square brackets in~\eqref{nega1} is strictly negative, which means that the objective value is strictly larger for $(\hat{\boldsymbol{y}}, \hat{\boldsymbol{z}})$ than for $\tilde{\boldsymbol{y}}$, $\tilde{\boldsymbol{z}}$, and hence $(\tilde{\boldsymbol{y}}, \tilde{\boldsymbol{z}})$ cannot be an optimal solution different from $(\boldsymbol{y}^*, \boldsymbol{z}^*)$.
\end{proof}

\section{Interpretation of Proposition~\ref{prop:pgf}}\label{app:interpretation_pgf}

The derivation of the PGF in Proposition~\ref{prop:pgf} induces a probabilistic interpretation in terms of the ordered vectors $\boldsymbol{S}$ of job types and geometrically distributed random variables.

More precisely, let us focus on an ordered vector $\boldsymbol{S} = [S_1,\dots,S_m] \in \mathcal{S}_m$, for some $m = 1,\dots,|\mathcal{S}|$. Then, the corresponding states in the state space are of the form 
\[
\boldsymbol{c}  = (S_1,\times_1,S_2,\times_2,S_3,\times_3\dots,S_m,\times_m),
\] 
where $\times_i \in \{S_1,\dots,S_i\}^n$ for any $n\ge 0$. Because of the same argument that gives rise to~\eqref{eq:proof_GF_inter}, the total contribution of these states to the generating function of $(Q_S)_{S\in\mathcal{S}}$ is given by 
\begin{equation}\label{eq:interpretation}
\begin{array}{rcl}
\mathbb{E}\left[ \prod\limits_{S\in\mathcal{S}} z_S^{Q_S}\mathds{1}\{\boldsymbol{S}(\boldsymbol{c}) = \boldsymbol{S}\}\right] &=& C  \prod\limits_{j=1}^m \frac{N\lambda p_{S_j}z_{S_j}}{\mu(S_1,\dots,S_j)}  \left(1-\frac{N\lambda}{\mu(S_1,\dots,S_j)}\sum\limits_{i=1}^j p_{S_i}z_{S_i}\right)^{-1}\\
& =& C  \prod\limits_{j=1}^m \frac{\frac{N\lambda p_{S_j}}{\mu(S_1,\dots,S_j)}}{1-\frac{N\lambda\sum_{i=1}^j p_{S_i}}{\mu(S_1,\dots,S_j)}}z_{S_j}\frac{1-\frac{N\lambda\sum_{i=1}^j p_{S_i}}{\mu(S_1,\dots,S_j)}}{1-\frac{N\lambda\sum_{i=1}^j p_{S_i}}{\mu(S_1,\dots,S_j)}\frac{\sum_{i=1}^j p_{S_i}z_{S_i}}{\sum_{i=1}^j p_{S_i}}},
\end{array}
\end{equation}
where $\boldsymbol{S}(\boldsymbol{c}) = \boldsymbol{S}$ denotes the event that the ordered vector of job types of the state $\boldsymbol{c}$ coincides with the vector $\boldsymbol{S}$.
Define $\hat{p}_{j,S_k}\coloneqq p_{S_k} / \sum_{i=1}^j p_{S_i}$ with $k=1,\dots,j$ and $j=1,\dots,m$, the probability that a job is of type $S_k$ given that it must be one of the job types in $\{S_1,\dots,S_j\}$, and rewrite the last factor of the above expression as
\begin{equation}
\begin{array}{rcl}
\dfrac{1-\frac{N\lambda\sum_{i=1}^j p_{S_i}}{\mu(S_1,\dots,S_j)}}{1-\frac{N\lambda\sum_{i=1}^j p_{S_i}}{\mu(S_1,\dots,S_j)}\frac{\sum_{i=1}^j p_{S_i}z_{S_i}}{\sum_{i=1}^j p_{S_i}}} &=& \left(1-\dfrac{N\lambda\sum_{i=1}^j p_{S_i}}{\mu(S_1,\dots,S_j)}\right) \sum\limits_{q = 0}^{\infty} \left(\dfrac{N\lambda\sum_{i=1}^j p_{S_i}}{\mu(S_1,\dots,S_j)}\dfrac{\sum_{i=1}^j p_{S_i}z_{S_i}}{\sum_{i=1}^j p_{S_i}}\right)^q\\
&=& \sum\limits_{q = 0}^{\infty}\left(\sum\limits_{i=1}^j \hat{p}_{j,S_i}z_{S_i}\right)^q \left(1-\dfrac{N\lambda\sum_{i=1}^j p_{S_i}}{\mu(S_1,\dots,S_j)}\right)\left(\dfrac{N\lambda\sum_{i=1}^j p_{S_i}}{\mu(S_1,\dots,S_j)}\right)^q\\
& = & \mathbb{E}[(H_j(z_{S_1},\dots,z_{S_j}))^{Q_j}]. 
\end{array}
\end{equation}
The random variable $Q_j$ is geometrically distributed with parameter $\frac{N\lambda\sum_{i=1}^j p_{S_i}}{\mu(S_1,\dots,S_j)}<1$ and 
\begin{equation}
H_j(z_{S_1},\dots,z_{S_j}) = \mathbb{E}\left[ \prod\limits_{i=1}^j z_{S_i}^{Y_{S_i}}\right],
\end{equation}
with $(Y_{S_1},\dots,Y_{S_j})$ a $j$-dimensional random vector where the $i$th component is equal to 1 and all other components equal to 0 with probability $\hat{p}_{j,S_i}$ with $i=1,\dots,j$. Observe that
\begin{equation}
\mathbb{E}\left[ \prod\limits_{S\in\mathcal{S}} z_S^{Q_S}\mathds{1}\{\boldsymbol{S}(\boldsymbol{c}) = \boldsymbol{S}\}\right] = \mathbb{E}\left[ \prod\limits_{S\in\mathcal{S}} z_S^{Q_S}\mid \boldsymbol{S}(\boldsymbol{c}) = \boldsymbol{S}\right] \mathbb{P}(\boldsymbol{S}(\boldsymbol{c}) = \boldsymbol{S}).
\end{equation}
Now, via substitution of $z_S\equiv 1$ in~\eqref{eq:interpretation} we see that
\begin{equation}
\mathbb{P}(\boldsymbol{S}(\boldsymbol{c}) = \boldsymbol{S}) = C \prod\limits_{j=1}^m \frac{\frac{N\lambda p_{S_j}}{\mu(S_1,\dots,S_j)}}{1-\frac{N\lambda\sum_{i=1}^j p_{S_i}}{\mu(S_1,\dots,S_j)}},
\end{equation}
hence
\begin{equation}\label{eq:ht_interpret_defHj}
\mathbb{E}\left[ \prod\limits_{S\in\mathcal{S}} z_S^{Q_S}\mid \boldsymbol{S}(\boldsymbol{c}) = \boldsymbol{S}\right] = \prod\limits_{j=1}^m z_{S_j} \mathbb{E}[(H_j(z_{S_1},\dots,z_{S_j}))^{Q_j}].
\end{equation}
Define $Q_{j,S_k}$ as the number of type-$S_k$ jobs that occur between the occurrences of the first $S_j$ and $S_{j+1}$ jobs in a state, with $j=1,\dots,m-1$ and $k=1,\dots,j$. Let $Q_{m,S_k}$ denote the number of type-$S_k$ jobs occurring after the first type-$S_m$ job, with $k=1,\dots,m$. For instance, in this state $\boldsymbol{c}= (S_1,S_1,S_2,S_1,S_2,S_2,S_3,S_1,S_1,S_1,S_3)$ we have that $Q_{2,S_1} = 1$ and $Q_{3,S_1} = 3$. Hence, by definition,
\begin{equation}
\mathbb{E}\left[ \prod\limits_{S\in\mathcal{S}} z_S^{Q_S}\mid \boldsymbol{S}(\boldsymbol{c}) = \boldsymbol{S}\right]  = \prod\limits_{j=1}^m z_{S_j}\cdot \mathbb{E}\left[ \prod\limits_{j=1}^m \prod\limits_{i=1}^jz_{S_j}^{Q_{i,S_j}}\right]  = \prod\limits_{j=1}^m z_{S_j}\cdot \mathbb{E}\left[ \prod\limits_{j=1}^mz_{S_j}^{ \sum_{i=1}^j Q_{i,S_j}}\right],
\end{equation}
which suggests that 
\begin{equation}
\mathbb{E}\left[\prod\limits_{i=1}^jz_{S_j}^{Q_{i,S_j}}\right]  =  \mathbb{E}[(H_j(z_{S_1},\dots,z_{S_j}))^{Q_j}].
\end{equation}
In particular, for $z_S \equiv z $
\begin{equation}
\mathbb{E}\left[z^{\sum_{i=1}^j Q_{i,S_j}}\right]  =  \mathbb{E}[z^{Q_j}],
\end{equation}
which implies that the number of jobs between the first occurrences of a type-$S_j$ and a type-$S_{j+1}$ is geometrically distributed with parameter $\frac{N\lambda\sum_{i=1}^j p_{S_i}}{\mu(S_1,\dots,S_j)}<1$, for $j=1,\dots,m-1$. Moreover, the identities of the various intermediate jobs are independent and identically distributed with probabilities $\hat{p}_{j,S_1},\dots,\hat{p}_{j,S_j}$. Likewise, the number of jobs in a state after the first occurrence of a type-$S_m$ job is geometrically distributed with parameter $\frac{N\lambda\sum_{i=1}^m p_{S_i}}{\mu(S_1,\dots,S_m)}<1$, and the labels of these jobs are independent and identically distributed with probabilities $\hat{p}_{m,S_1},\dots,\hat{p}_{m,S_m}$. \\

When the system satisfies the local stability condition~\eqref{stabcond3} and~\eqref{stabcond4}, as described in Subsection~\ref{subsec:max_stable}, using similar arguments as in the proof of Theorem~\ref{statespacecollapse_general} it can be shown that
\begin{equation}\label{eq:ht_allJobtypes}
\lim_{\lambda\uparrow\mu}\mathbb{P}(\boldsymbol{S}(\boldsymbol{c}) = \boldsymbol{S}) = \dfrac{\prod_{j=1}^{|\mathcal{S}|-1}\left(\mu(S_1,\dots,S_j) -N\mu \sum_{i=1}^jp_{S_i}\right)^{-1}}{\sum_{T\in\mathcal{S}_{|\mathcal{S}|}}\prod_{j=1}^{|\mathcal{S}|-1}\left(\mu(T_1,\dots,T_j) -N\mu \sum_{i=1}^jp_{T_i}\right)^{-1}}.
\end{equation}
Hence, $\mathbb{P}(\boldsymbol{S}(\boldsymbol{c}) \in \mathcal{S}_{|\mathcal{S}|})$ tends to one if $\lambda$ approaches $\mu$, 
i.e., all different job types will occur at any state of the system in the heavy-traffic regime. Or in other words, $\mathbb{P}(\boldsymbol{S}(\boldsymbol{c}) =\boldsymbol{S})\approx 0$ when $\boldsymbol{S}$ does not cover all job types and $\lambda$ is close to $\mu$. Observe that the expression
\begin{equation}
\mathbb{E}[(H_{|\mathcal{S}|}(z_{S_1},\dots,z_{S_{|\mathcal{S}|}}))^{Q_{|\mathcal{S}|}}]  = \dfrac{1-\frac{\lambda}{\mu}}{1-\frac{\lambda}{\mu}\sum_{S\in\mathcal{S}}p_Sz_s}
\end{equation}
corresponds to the generating function of the joint stationary distribution of a multi-class M/M/1 queue with class probabilities $(p_S)_{S\in\mathcal{S}}$. 
In what follows we will argue that this term will indeed be dominant in the heavy-traffic limit. More precisely, when $z_S = \exp(-(1-\frac{\lambda}{\mu})t_S)$ for all $S\in\mathcal{S}$ and $\lambda$ approaches $\mu$,
\begin{equation}
\mathbb{E}\left[ \prod\limits_{S\in\mathcal{S}} z_S^{Q_S}\right] \approx \mathbb{E}[(H_{|\mathcal{S}|}(z_{S_1},\dots,z_{S_{|\mathcal{S}|}}))^{Q_{|\mathcal{S}|}}],
\end{equation}
where the right-hand side converges to $\left(1+\sum_{S\in\mathcal{S}}p_St_S\right)^{-1}$. By Feller's Convergence theorem this implies that $\left(1-\lambda/\mu\right)\left(Q_S\right)_{S\in\mathcal{S}}$ converges in distribution to $\mathrm{Exp}(1)(p_S)_{S\in\mathcal{S}}$, as stated and proved rigorously in Theorems~\ref{statespacecollapse} and~\ref{statespacecollapse_general}.

First we note that the parameter $\frac{N\lambda\sum_{i=1}^j p_{S_i}}{\mu(S_1,\dots,S_j)}$ of the geometrically distributed random variable $Q_j$, for $j=1,\dots,|\mathcal{S}|-1$, remains strictly smaller than 1 in the heavy-traffic regime.
Then, by definition of $z_S$ and $\hat{p}_{j,S_i}$, $z_{S} \approx 1$ and $\sum_{i=1}^j \hat{p}_{j,S_i}z_{S_i} \approx 1$ when $\lambda$ is close to $\mu$.
This results in 
\begin{equation}\label{eq:ht_approx1}
\mathbb{E}[(H_{j}(z_{S_1},\dots,z_{S_{j}}))^{Q_{j}}] = \mathbb{E}\left[\left(\sum\limits_{i=1}^j \hat{p}_{j,S_i}z_{S_i}\right)^{Q_{j}}\right] \approx  \mathbb{E}\left[1^{Q_{j}}\right] = 1.
\end{equation}
This means that the (on average finite) number of jobs in between first occurrences of different job types in a state will have a negligible contribution to the (scaled) total number of jobs in the system. However, the (infinite number of) jobs that occur once all job types are encountered in a state will significantly contribute to the joint stationary distribution of the number of jobs of each type in the limiting regime.
So, from~\eqref{eq:ht_interpret_defHj} and~\eqref{eq:ht_allJobtypes} we know that 
\begin{equation}
\begin{array}{rcl}
\mathbb{E}\left[ \prod\limits_{S\in\mathcal{S}} z_S^{Q_S}\right] &\approx & \sum\limits_{\boldsymbol{S}\in\mathcal{S}_{|\mathcal{S}|}} \mathbb{E}\left[ \prod\limits_{S\in\mathcal{S}} z_S^{Q_S}\mid \boldsymbol{S}(\boldsymbol{c})=\boldsymbol{S}\right]\mathbb{P}(\boldsymbol{S}(\boldsymbol{c})=\boldsymbol{S})\\
 &= & \sum\limits_{\boldsymbol{S}\in\mathcal{S}_{|\mathcal{S}|}} \prod\limits_{j=1}^{|\mathcal{S}|} z_{S_j}\mathbb{E}[(H_{j}(z_{S_1},\dots,z_{S_{j}}))^{Q_{j}}] \mathbb{P}(\boldsymbol{S}(\boldsymbol{c})=\boldsymbol{S}).\\
\end{array}
\end{equation}
Due to the observation in~\eqref{eq:ht_approx1}, we see that 
\begin{equation}
\prod\limits_{j=1}^{|\mathcal{S}|} z_{S_j}\mathbb{E}[(H_{j}(z_{S_1},\dots,z_{S_{j}}))^{Q_{j}}] \approx \left(\prod\limits_{j=1}^{|\mathcal{S}|-1} 1 \right) \mathbb{E}[(H_{|\mathcal{S}|}(z_{S_1},\dots,z_{S_{|\mathcal{S}|}}))^{Q_{|\mathcal{S}|}}]  = \mathbb{E}[(H_{|\mathcal{S}|}(z_{S_1},\dots,z_{S_{|\mathcal{S}|}}))^{Q_{|\mathcal{S}|}}],
\end{equation}
where the last factor is independent of the order of the vector $\boldsymbol{S}$ of all $|\mathcal{S}|$ different job types. Hence,
\begin{equation}
\begin{array}{rcl}
\mathbb{E}\left[ \prod\limits_{S\in\mathcal{S}} z_S^{Q_S}\right] 
  &\approx & \mathbb{E}[(H_{|\mathcal{S}|}(z_{S_1},\dots,z_{S_{|\mathcal{S}|}}))^{Q_{|\mathcal{S}|}}]\sum\limits_{\boldsymbol{S}\in\mathcal{S}_{|\mathcal{S}|}} \mathbb{P}(\boldsymbol{S}(\boldsymbol{c})=\boldsymbol{S})\\
  &\approx & \mathbb{E}[(H_{|\mathcal{S}|}(z_{S_1},\dots,z_{S_{|\mathcal{S}|}}))^{Q_{|\mathcal{S}|}}]
\end{array}
\end{equation}
when $\lambda$ approaches $\mu$.

\section{Proof of Theorem~\ref{statespacecollapse_general} for the c.o.s.\ mechanism}\label{app:cos}

In Appendix~\ref{app:cos_indirect} we prove Theorem~\ref{statespacecollapse_general} for the c.o.s.\ mechanism based on the relation between the product-form expressions~\eqref{eq:statdistr} and~\eqref{eq:product_form_cos}; we refer to this proof method as the indirect method. The result of Theorem~\ref{statespacecollapse_general} for the c.o.s.\ mechanism can also be derived in a similar way as outlined in Section~\ref{sec:HT} for the c.o.c.\ mechanism, relying on its joint~PGF. This proof method is referred to as the direct method, and is described in Appendix~\ref{app:cos_direct}.

\subsection{Indirect method}\label{app:cos_indirect}

The indirect proof of Theorem~\ref{statespacecollapse_general} for the redundancy c.o.s. policy relies on the results derived in Section~\ref{sec:HT}.
\begin{proof}[Proof of Theorem~\ref{statespacecollapse_general} \textup{(}c.o.s.\ mechanism\textup{)}.]
Let $\mathcal{T}^*\subseteq\mathcal{S}$, $p_{\mathcal{T}^*}$, $\mu_{\mathcal{T}^*}$ and $\lambda^*=\mu_{\mathcal{T}^*}/(Np_{\mathcal{T}^*})$ be as defined in Section~\ref{sec:main_results}.
In order for Theorem~\ref{statespacecollapse_general} to hold for redundancy c.o.s., it is sufficient to show that
\begin{equation}\label{eq:convergence_waiting}
\left(1 - \frac{\lambda}{\lambda^*}\right) (\tilde{Q}_S)_{S \in {\mathcal S}} \xrightarrow{d}
\left(\mathrm{Exp}(1) \left(\frac{p_S}{p_{\mathcal{T}^*}}\right)_{S \in {\mathcal T}^*}, (0)_{S\notin\mathcal{T}^*}\right),
\end{equation}
as $\lambda \uparrow \lambda^*$ with $\tilde{Q}_S$ the number of {\em waiting} type-$S$ jobs. 
It is clear that the number of type-$S$ jobs in service can never exceed the total number of servers $N$, and a finite number of jobs will have a negligible contribution in the heavy-traffic limit. More precisely, the MGFs of $\left(1 - \frac{\lambda}{\lambda^*}\right) (\tilde{Q}_S)_{S \in {\mathcal S}}$ and $\left(1 - \frac{\lambda}{\lambda^*}\right) ({Q}_S)_{S \in {\mathcal S}}$ coincide in the heavy-traffic limit as
\begin{equation}
\prod_{S\in\mathcal{S}}\mathrm{e}^{-\left(1-\frac{\lambda}{\lambda^*}\right)t_S N}\mathbb{E}\left[\prod\limits_{S\in\mathcal{S}}\mathrm{e}^{-\left(1-\frac{\lambda}{\lambda^*}\right)t_S\tilde{Q}_S}\right] \le
 \mathbb{E}\left[\prod\limits_{S\in\mathcal{S}}\mathrm{e}^{-\left(1-\frac{\lambda}{\lambda^*}\right)t_SQ_S}\right] \le \mathbb{E}\left[\prod\limits_{S\in\mathcal{S}}\mathrm{e}^{-\left(1-\frac{\lambda}{\lambda^*}\right)t_S\tilde{Q}_S}\right]  
\end{equation}
with $\boldsymbol{t}\ge 0$ and the deterministic product on the left-hand side converging to 1 as $\lambda\uparrow\lambda^*$.
Hence, proving~\eqref{eq:convergence_waiting} would yield Theorem~\ref{statespacecollapse_general} for the c.o.s.\ mechanism.

We will show the following heavy-traffic result for the MGF of the number of waiting jobs of each type $(\tilde{Q}_S)_{S\in\mathcal{S}}$.
\begin{equation}\label{eq:genfuncconv_cos}
\mathbb{E}\left[ \prod\limits_{S\in\mathcal{S}} z_S^{\tilde{Q}_S} \right] = \sum\limits_{\boldsymbol{c},\boldsymbol{u}} \pi_{\text{c.o.s.}}(\boldsymbol{c},\boldsymbol{u}) \prod\limits_{S\in\mathcal{S}} z_S^{q_S^{\boldsymbol{c}}} \rightarrow \left(1+\sum\limits_{S\in\mathcal{T}^*} \frac{p_S}{p_{\mathcal{T}^*}}t_S \right)^{-1},
\end{equation}
as $\lambda\uparrow\lambda^*$ and $z_S \coloneqq \exp(-(1-\frac{\lambda}{\lambda^*})t_S)$ with $t_S  \ge 0$ for all $S\in\mathcal{S}$. The number of waiting type-$S$ jobs in state $\boldsymbol{c}$ is denoted by $q_S^{\boldsymbol{c}}$. This establishes~\eqref{eq:convergence_waiting} due to Feller's convergence theorem~\citep{Feller1971}.
We will split the state space into two disjoint parts while computing the generating function: states where all servers that can process jobs of types $\mathcal{T}^*$ are busy and states where at least one of these servers is idle. With a slight abuse of notation, we refer to these compatible servers as the set $\mathcal{T}^*$ as well. 
\begin{equation}\label{eq:proofHT_parts12}
\sum\limits_{\boldsymbol{c},\boldsymbol{u}} \pi_{\text{c.o.s.}}(\boldsymbol{c},\boldsymbol{u}) \prod\limits_{S\in\mathcal{S}} z_S^{q_S^{\boldsymbol{c}}} =
\underbrace{\sum\limits_{\substack{\boldsymbol{c},\boldsymbol{u}: \\ \boldsymbol{u}\cap\mathcal{T}^*=\emptyset}} \pi_{\text{c.o.s.}}(\boldsymbol{c},\boldsymbol{u}) \prod\limits_{S\in\mathcal{S}} z_S^{q_S^{\boldsymbol{c}}}}_{\text{Part 1}} 
+ \underbrace{\sum\limits_{\substack{\boldsymbol{c},\boldsymbol{u}: \\ \boldsymbol{u}\cap\mathcal{T}^*\neq\emptyset}} \pi_{\text{c.o.s.}}(\boldsymbol{c},\boldsymbol{u}) \prod\limits_{S\in\mathcal{S}} z_S^{q_S^{\boldsymbol{c}}}}_{\text{Part 2}} 
\end{equation}
The limits of Parts 1 and 2 will be evaluated separately.\\

\textbf{Part~1:} To establish the limit of Part~1, we rely on the stationary relation between the c.o.s.\ and c.o.c.\ mechanism. It is known that $\mathbb{P}_{\text{c.o.s.}}\{\boldsymbol{c} \mid \boldsymbol{u}=\emptyset\} = \pi_{\text{c.o.c.}}(\boldsymbol{c})$ for any job-type sequence~$\boldsymbol{c}$ \citep{Adan2018}.
\cite{Gardner2020} observed that this relation can even be extended further by conditioning on the set of idle servers $\boldsymbol{u}$, without assuming that $\boldsymbol{u}=\emptyset$.
So, 
\begin{equation}\label{eq:relation2}
    \mathbb{P}_{\text{c.o.s.}}\{\boldsymbol{c} \mid \boldsymbol{u}\} = \pi^{\boldsymbol{u}}_{\text{c.o.c.}}(\boldsymbol{c}),
\end{equation}
where $\pi^{\boldsymbol{u}}_{\text{c.o.c.}}(\boldsymbol{c})$ denotes the stationary distribution of a truncated system operating under the redundancy c.o.c.\ policy. Only the servers $\{1,\dots,N\} \setminus \boldsymbol{u}$ are present and only those types of jobs in $\mathcal{S}$ that have no compatible servers in $\boldsymbol{u}$ are considered; the set of these job types is denoted by $\mathcal{S}'\subseteq\mathcal{S}$. Note that $\mathcal{T}^*\subseteq\mathcal{S}'$ as all servers compatible to the job types in $\mathcal{T}^*$ are busy by assumption. The relation in \eqref{eq:relation2} allows us to rewrite Part~1 as
\begin{equation}
    \sum\limits_{\boldsymbol{u}:\boldsymbol{u}\cap\mathcal{T}^* = \emptyset}\mathbb{P}_{\text{c.o.s.}} (\boldsymbol{u})\left( \sum\limits_{\boldsymbol{c}} \pi^{\boldsymbol{u}}_{\text{c.o.c.}}(\boldsymbol{c}) \prod\limits_{S\in\mathcal{S}} z_S^{q_S^{\boldsymbol{c}}}\right).
\end{equation}
Focusing on the latter summation, we observe that the jobs in state $\boldsymbol{c}$ can only be compatible with the servers $\{1,\dots,N\}\setminus\boldsymbol{u}$, hence they are of types in~$\mathcal{S}'$. So, this summation is the joint~PGF of the number of jobs of each type in~$\mathcal{S}'$ of the truncated system under the c.o.c.\ mechanism.
 (As a comparison, \eqref{eq:general_gen_f} is the expression for the joint generating function of a system with job types in $\mathcal{S}$, servers $\{1,\dots,N\}$ and operating according to the redundancy c.o.c.\ policy.) Moreover, $\mathcal{T}^*\subseteq\mathcal{S}$ is also the critical subset of types of the truncated system, otherwise this would violate the fact that we have a unique critical set $\mathcal{T}^*$ of the original (larger) system. This allows us to apply the result of Theorem~\ref{statespacecollapse_general} for the truncated system operating according to the redundancy c.o.c.\ policy, hence the sum converges to $(1+\sum_{S\in\mathcal{T}^*} \frac{p_s}{p_{\mathcal{T}^*}}t_S)^{-1}$. Furthermore, note that 
 \begin{equation}\label{eq:argument_cos_sum}
  \sum\limits_{\boldsymbol{u}:\boldsymbol{u}\cap\mathcal{T}^* = \emptyset}\mathbb{P}_{\text{c.o.s.}} (\boldsymbol{u}) 
 \end{equation}
 equals the probability that all servers in $\mathcal{T}^*$ are non-idle as $\mathbb{P}_{\text{c.o.s.}} (\boldsymbol{u})$ denotes the probability that precisely the servers in $\boldsymbol{u}$ are idle for any $\boldsymbol{u}$ with $\boldsymbol{u}\cap\mathcal{T}^* = \emptyset$. Intuitively it is clear that the above expression will converge to 1 as $\lambda\uparrow\lambda^*$ since the aggregate arrival rate to the servers in $\mathcal{T}^*$ reaches the aggregate service rate. The formal proof of this statement is deferred to Appendix~\ref{app:cos_direct}.
 From these observations we conclude that Part~1 converges to the right-hand side of~\eqref{eq:genfuncconv_cos} when $\lambda\uparrow\lambda^*$.\\

\textbf{Part~2:} It remains to be shown that Part~2 converges to~0 as $\lambda\uparrow\lambda^*$. Intuitively this is clear as all states with idle servers in $\mathcal{T}^*$ are considered, while we know that the arrival rate to these servers reached the boundary of the stability region.
This reasoning involves an interchange of limit and infinite sum and is justified using Lebesgue's dominated convergence theorem.

First, it is shown that the point-wise limit exists for any fixed job sequence~$\boldsymbol{c}$ and vector of idle servers~$\boldsymbol{u}$. Observe that
\begin{equation}
0\le \lim_{\lambda\uparrow\lambda^*}\pi_{\text{c.o.s.}}(\boldsymbol{c},\boldsymbol{u})\prod_{S\in\mathcal{S}} \exp\left(-\left(1-\frac{\lambda}{\lambda^*}\right)t_Sq_S^{\boldsymbol{c}}\right) \le \lim_{\lambda\uparrow\lambda^*}\pi_{\text{c.o.s.}}(\boldsymbol{c},\boldsymbol{u}) \le \lim_{\lambda\uparrow\lambda^*} \pi_{\text{c.o.s.}}(\emptyset,\emptyset)\prod\limits_{j=1}^{|\boldsymbol{u}|} \frac{\mu_{u_j}}{\lambda_{\mathcal{C}(u_1,\dots,u_j)}},
\end{equation}
where the last inequality is due to the stability conditions~\eqref{stabcond1} and~\eqref{stabcond2}. The normalization constant $\pi_{\text{c.o.s.}}(\emptyset,\emptyset)$ tends to 0 as the servers in $\mathcal{T}^*$ will be busy with probability 1 and jobs of the corresponding types will start to accumulate.
Now,
\begin{equation}
0\le\lim_{\lambda\uparrow\lambda^*}\prod\limits_{j=1}^{|\boldsymbol{u}|} \frac{\mu_{u_j}}{\lambda_{\mathcal{C}(u_1,\dots,u_j)}} \le \lim_{\lambda\uparrow\lambda^*}\prod\limits_{j=1}^{|\boldsymbol{u}|} \frac{\mu^+}{\lambda N p^- } \le \max\left\{1, \left( \frac{\mu^+p_{\mathcal{T}^*}}{\mu_{\mathcal{T}^*} p^- }\right)^{|\boldsymbol{u}|}\right\}<\infty,
\end{equation}
with $\mu^+\coloneqq\max\{\mu_n \colon n=1,\dots,N\}$ and $p^- \coloneqq \min\{p_S\colon S\in \mathcal{S}\}$, respectively.
Hence, the point-wise limit is 0.

Second, define the function $h(\boldsymbol{c},\boldsymbol{u})\coloneqq\pi_{\text{c.o.s.}}(\boldsymbol{c},\boldsymbol{u})$; obviously $h$ dominates $\pi_{\text{c.o.s.}}(\boldsymbol{c},\boldsymbol{u})\prod_{S\in\mathcal{S}} z_S^{q_S^{\boldsymbol{c}}}$. Moreover, it is trivial that $h$ is finitely summable over all states under consideration.

Hence, Lebesgue's dominated convergence theorem can be applied to conclude convergence of Part~2 to 0 as $\lambda\uparrow\lambda^*$. This concludes the proof. 
\end{proof}

\begin{remark}\normalfont
The stationary distribution under the assignment rules according to the assignment condition and ALIS coincide when aggregating all states with the same central queue of jobs~$\boldsymbol{c}$ but a different order of the idle servers $\boldsymbol{u}$ \cite[Theorem~2.2]{adan2014skill}. Hence, Theorem~\ref{statespacecollapse_general} is also true when using the assignment condition to assign jobs to idle servers. 
\end{remark}

\subsection{Direct method}\label{app:cos_direct}
 
\begin{proposition}\label{prop:pgfcos}
Assuming that the stability conditions~\eqref{stabcond1} and~\eqref{stabcond2} are satisfied, the joint~PGF of the number of waiting jobs of each type for the redundancy c.o.s.\ policy is given by 
\begin{equation}\label{eq:pgfcos}
\mathbb{E}\left[\prod\limits_{S \in {\mathcal S}} z_S^{\tilde{Q}_S}\right] = \frac{g(\boldsymbol{z})}{g(\boldsymbol{1})},
\end{equation}
where $\boldsymbol{z}$ and $\boldsymbol{1}$  are $|\mathcal{S}|$-dimensional vectors with entries $|z_S| \le 1$ and
\[
g(\boldsymbol{z})=\sum\limits_{L=0}^{N}\sum\limits_{\boldsymbol{u}\in \mathcal{N}_L}\prod\limits_{l=1}^L \frac{\mu_{u_l}}{\lambda_{\mathcal{C}(u_1,\dots,u_l)}} 
 \sum\limits_{m=0}^{|\mathcal{S}|}\sum\limits_{\boldsymbol{S}\in \mathcal{S}_m^{\boldsymbol{u}}} \prod\limits_{j=1}^m \frac{N\lambda p_{S_j}z_{S_j}}{\mu(S_1,\dots,S_j)}\left(1{-}\frac{N\lambda}{\mu(S_1,\dots,S_j)}\sum\limits_{i=1}^j p_{S_i}z_{S_i}\right)^{-1}.
\]
The set consisting of all ordered vectors of $L$ idle servers is denoted by $\mathcal{N}_L$.
The $m$-dimensional vector $\boldsymbol{S}$ consists of $m$ different waiting job types that are not compatible with the idle servers $\boldsymbol{u}$, and the set consisting of all these vectors is denoted by $\mathcal{S}_m^{\boldsymbol{u}}$.
\end{proposition}

\begin{proof}[Proof of Proposition~\ref{prop:pgfcos}.]
Using the stationary distribution given in~\eqref{eq:product_form_cos}, the joint~PGF of the number of waiting jobs of each type may be written as
\begin{equation}
    \mathbb{E}\left[\prod\limits_{S\in\mathcal{S}}z_{S}^{\tilde{Q}_S}\right] = \sum\limits_{\boldsymbol{c},\boldsymbol{u}} \pi_{\text{c.o.s.}}(\boldsymbol{c},\boldsymbol{u})\prod\limits_{S\in\mathcal{S}}z_{S}^{q_S^{\boldsymbol{c}}},
\end{equation}
with $\boldsymbol{z}$ an $|\mathcal{S}|$-dimensional vector with entries $|z_S|\le 1$ and $q_S^{\boldsymbol{c}}$ the total number of (waiting) type-$S$ jobs in state $\boldsymbol{c}$. The summation over all possible states will be conducted in a few steps, just as we did in the proof of Proposition~\ref{prop:pgf}. First, the number of idle servers $L$ is fixed, with $L=0,\dots,N$. Then the ordered vector $\boldsymbol{u}$ containing $L$ idle servers in order of longest idleness is set, the collection of all these vectors is denoted by $\mathcal{N}_L$. From then the three steps in the proof of Proposition~\ref{prop:pgf} are applied to sum over all job states $\boldsymbol{c}$ that have no compatible server in $\boldsymbol{u}$, yielding the desired expression.
\end{proof}

Along the same lines as the proof of Proposition~\ref{prop:pgf}, we observe that the value $g(\boldsymbol{1})$ coincides with the normalization constant $C'$ in~\eqref{eq:product_form_cos}. Using the notation outlined in Section~\ref{sec:main_results}, Proposition~\ref{prop:pgfcos} can be used to give an alternative proof of Theorem~\ref{statespacecollapse_general} and to formalize the convergence argument of~\eqref{eq:argument_cos_sum} in Appendix~\ref{app:cos_indirect}.

\begin{proof}[Proof of Theorem~\ref{statespacecollapse_general} \textup{(}c.o.s.\ mechanism\textup{)}.]
As argued in Subsection~\ref{subsec:cos}, it is sufficient to show that~\eqref{eq:convergence_waiting} holds in order to establish the result in Theorem~\ref{statespacecollapse_general}. To obtain the moment generating function of $(1-\frac{\lambda}{\lambda^*})(\tilde{Q}_S)_S$ define $z_S\coloneqq \exp(-(1-\frac{\lambda}{\lambda^*})t_S)$ and use the expression obtained in Proposition~\ref{prop:pgfcos}. Observe that for any $L=0,\dots,N$, for any $\boldsymbol{u}\in \mathcal{N}_L$ and for any $l=1,\dots,L$
\begin{equation}
    \lim_{\lambda\uparrow\lambda^*} \frac{\mu_{u_l}}{\lambda_{\mathcal{C}(u_1,\dots,u_l)}}= \dfrac{\mu_{u_l}}{N\lambda^*\sum\limits_{S:S\cap\{u_1,\dots,u_l\}\neq\emptyset}p_s} \in (0,\infty).
\end{equation}
As $\mathcal{N}_L$ is a finite-size set for all $L=0,\dots,N$, together with the observations in~\eqref{eq:obs1} and~\eqref{eq:obs2}, we can apply similar steps as in the proof of Theorem~\ref{statespacecollapse_general} for the c.o.c.\ mechanism as outlined in Section~\ref{sec:HT}. This concludes the proof.
\end{proof}

\begin{lemma}\label{lem:argument_cos}
Let $\mathbb{P}_{\text{c.o.s.}}(\boldsymbol{u})$ denote the stationary probability that precisely the servers in the vector $\boldsymbol{u}$ are idle, and $\boldsymbol{u}_l$ is the $l$th longest idle server. Then
\begin{equation}\label{eq:argument_cos}
\lim_{\lambda\uparrow\lambda^*}\sum\limits_{\boldsymbol{u}:\boldsymbol{u}\cap \mathcal{T}^* = \emptyset} \mathbb{P}_{\text{c.o.s.}}(\boldsymbol{u}) = 1.
\end{equation} 
In other words, the servers in set $\mathcal{T^*}$ will be occupied with probability 1 in the limit.
\end{lemma}

\begin{proof}
For any fixed $\boldsymbol{u}$, let 
\[
B_{\boldsymbol{u}} \coloneqq \prod\limits_{l=1}^{|\boldsymbol{u}|}\frac{\mu_{\boldsymbol{u}_l}}{\lambda_{\mathcal{C}(\boldsymbol{u}_1,\dots,\boldsymbol{u}_l)}}.
\]
Then, \eqref{eq:argument_cos} can be rewritten as follows, using the stationary distribution in~\eqref{eq:product_form_cos}, 
\begin{equation}
\sum\limits_{\boldsymbol{u}: \boldsymbol{u}\cap \mathcal{T}^* = \emptyset}\sum\limits_{\boldsymbol{c}}\pi_{\text{c.o.s.}}(\boldsymbol{c},\boldsymbol{u})
 = C'\sum\limits_{\boldsymbol{u}: \boldsymbol{u}\cap \mathcal{T}^* = \emptyset}B_{\boldsymbol{u}}\sum\limits_{\boldsymbol{c}}\prod\limits_{i=1}^{|\boldsymbol{c}|}\frac{N\lambda p_{\boldsymbol{c}_i}}{\mu(\boldsymbol{c}_1,\dots,\boldsymbol{c}_i)}.
\end{equation}
Now, 
\begin{equation}\label{eq:D_u}
D_{\boldsymbol{u}}\coloneqq  \sum\limits_{\boldsymbol{c}}\prod\limits_{i=1}^{|\boldsymbol{c}|}\frac{N\lambda p_{\boldsymbol{c}_i}}{\mu(\boldsymbol{c}_1,\dots,\boldsymbol{c}_i)} = \sum\limits_{m=0}^{|\mathcal{S}|}\sum\limits_{\boldsymbol{S}\in \mathcal{S}_m^{\boldsymbol{u}}} \prod\limits_{j=1}^m \frac{N\lambda p_{S_j}z_{S_j}}{\mu(S_1,\dots,S_j)}\left(1{-}\frac{N\lambda}{\mu(S_1,\dots,S_j)}\sum\limits_{i=1}^j p_{S_i}z_{S_i}\right)^{-1},
\end{equation}
with $\mathcal{S}_m^{\boldsymbol{u}}$ the set of all ordered vectors $\boldsymbol{S}$ consisting of $m$ distinct  job types that are not compatible with the servers in $\boldsymbol{u}$. The second equality in the above expression can be established using the three-step approach outlined in the proof of Proposition~\ref{prop:pgf}. The normalization constant $C'$ is given by the value $g(\boldsymbol{1})$ obtained in Proposition~\ref{prop:pgfcos}. It can be seen that the first two summations, concerning the vectors of idle servers, can be split into two separate summations conditional on having an empty or a non-empty intersection with the set $\mathcal{T}^*$. So, \eqref{eq:argument_cos} is given by
\begin{equation}
C' \sum\limits_{\boldsymbol{u}:\boldsymbol{u}\cap \mathcal{T}^* = \emptyset}B_{\boldsymbol{u}} D_{\boldsymbol{u}} = \frac{\sum\limits_{\boldsymbol{u}: \boldsymbol{u}\cap \mathcal{T}^* = \emptyset}B_{\boldsymbol{u}} D_{\boldsymbol{u}}}{\sum\limits_{\boldsymbol{u}: \boldsymbol{u}\cap \mathcal{T}^* = \emptyset}B_{\boldsymbol{u}} D_{\boldsymbol{u}} + \sum\limits_{\boldsymbol{u}: \boldsymbol{u}\cap \mathcal{T}^* \neq \emptyset}B_{\boldsymbol{u}} D_{\boldsymbol{u}}}= \left(1+\frac{\sum\limits_{\boldsymbol{u}: \boldsymbol{u}\cap \mathcal{T}^* \neq \emptyset}B_{\boldsymbol{u}} D_{\boldsymbol{u}}}{\sum\limits_{\boldsymbol{u}: \boldsymbol{u}\cap \mathcal{T}^* = \emptyset}B_{\boldsymbol{u}} D_{\boldsymbol{u}}}\right)^{-1}.
\end{equation}
In order to obtain the limit of~\eqref{eq:argument_cos} as $\lambda\uparrow\lambda^*$, we first use the observation made in the proof of Proposition~\ref{prop:pgfcos} that the limit of $B_{\boldsymbol{u}}$ exists and is strictly positive for any vector of idle server $\boldsymbol{u}$. Second, we investigate the limit of $D_{\boldsymbol{u}}$ for any $\boldsymbol{u}$ based on~\eqref{eq:obs1} and~\eqref{eq:obs2} in the proof of Theorem~\ref{statespacecollapse_general} for the c.o.c.\ mechanism. In case $\boldsymbol{u}\cap \mathcal{T}^* \neq \emptyset$, no vector $\boldsymbol{S}$ in the summation of~\eqref{eq:D_u} will consist of all servers in $\mathcal{T}^*$. Hence,
\begin{equation}\label{eq:lim_component_cos}
\lim\limits_{\lambda\uparrow\lambda^*} \sum\limits_{\boldsymbol{u}: \boldsymbol{u}\cap \mathcal{T}^* \neq \emptyset}B_{\boldsymbol{u}} D_{\boldsymbol{u}} \in (0,\infty).
\end{equation}
In case $\boldsymbol{u}\cap \mathcal{T}^* = \emptyset$, some vectors $\boldsymbol{S}$ in the summation of~\eqref{eq:D_u} will consist of all servers in $\mathcal{T}^*$. Thus, $D_{\boldsymbol{u}}$ diverges if $\lambda\uparrow\lambda^*$. We conclude that~\eqref{eq:argument_cos} tends to~1 as $\lambda\uparrow\lambda^*$.
\end{proof}

\section{Convergence of the moments}\label{app:conv_moments}
\subsection{Preliminaries and notation}
In order to establish the result in Theorem~\ref{th:conv_in_mean} when the local stability conditions~\eqref{stabcond3} and~\eqref{stabcond4} are satisfied, we rely on the joint~PGFs derived in Proposition~\ref{prop:pgf} and~\ref{prop:pgfcos} for the c.o.c.\ and c.o.s.\ mechanisms, respectively.
Moreover, we will make use of the two following results.

\begin{lemma}[General Leibniz rule] Let $f_1,f_2\dots,f_m$ be $m$ functions that are at least $n$ times differentiable. Then, for any $n\ge 1$,
\begin{equation}
\left(\prod_{i=1}^m f_i\right)^{(n)} = \sum\limits_{\boldsymbol{k}\colon k_1+\dots + k_m=n}\binom{n}{k_1,\dots,k_m}\prod_{i=1}^m f_i^{(k_i)}.
\end{equation}
\end{lemma}

\begin{lemma}[Fa\`a di Bruno's formula]
Let $f$ and $g$ be functions for which the necessary derivatives are well defined. Then, for any $n\ge 1$, 
\begin{equation}
\frac{\mathrm{d}^n}{\mathrm{d}t^n} f\left(g(t)\right) = \sum\limits_{\boldsymbol{m}\in R(n)}  \frac{n!}{m_1!m_2!\dots m_n!} f^{(|\boldsymbol{m}|)}\left(g(t)\right)\prod\limits_{j=1}^n \left(\frac{g^{(j)}(t)}{j!}\right)^{m_j},
\end{equation}
with
\begin{equation}\label{eq:R(n)}
R(n) \coloneqq \left\{ \boldsymbol{m} \in\mathbb{N}^n \colon 1 m_1+2m_2+\dots + nm_n = n \right\}
\text{~~~and~~~}
|\boldsymbol{m}| \coloneqq m_1+m_2+\dots+m_n.
\end{equation}
\end{lemma}

In the remainder of this section we will use both notations $f^{n}(t)$ and $\frac{\mathrm{d}^n}{\mathrm{d}t^n}f(t)$ to denote the $n$th derivative of a function $f$ with respect to $t$. Moreover, we set $\rho\coloneqq\frac{\lambda}{\mu}$.

\subsection{The moments of the total number of jobs}

\begin{theorem}\label{th:conv_moments_tot}
If the local stability conditions~\eqref{stabcond3} and~\eqref{stabcond4} hold, then for both the c.o.c.\ and c.o.s. mechanisms
\begin{equation}
\lim_{\lambda\uparrow\mu}\mathbb{E}\left[ \left((1-\rho)Q\right)^n \right] = n! 
\end{equation}
for any $n\ge 1$.
\end{theorem}
The following corollary could already be deduced from Theorem~\ref{statespacecollapse}; it can also be deduced from Theorem~\ref{th:conv_moments_tot} via the method of moments.

\begin{corollary} 
If the local stability conditions~\eqref{stabcond3} and~\eqref{stabcond4} hold, then for both the c.o.c.\ and c.o.s. mechanisms
\begin{equation}
\left(1-\frac{\lambda}{\mu}\right)Q   \xrightarrow{d} \mathrm{Exp}(1)
\end{equation}
as $\lambda\uparrow\mu$. 
\end{corollary}

To prove Theorem~\ref{th:conv_moments_tot} we first derive expressions for the $n$th moments of $Q$ and $\tilde{Q}$ for any fixed value of $\lambda$ within the stability region for the c.o.c.\ and c.o.s. mechanisms, respectively.

\begin{lemma}\label{lem:nth_moment_total}
Assuming that the stability conditions~\eqref{stabcond1} and~\eqref{stabcond2} are satisfied, the $n$th moment of the total number of jobs in the system, $Q$, under the c.o.c.\ mechanism is given by
\begin{equation}
\mathbb{E}\left[ \left((1-\rho)Q\right)^n \right] = \frac{n!(1-\rho)^n}{f(1)} \sum\limits_{m=1}^{|\mathcal{S}|}\sum\limits_{\boldsymbol{S}\in \mathcal{S}_m} \gamma(\boldsymbol{S}) \prod\limits_{j=1}^m A(\boldsymbol{S},j)\prod\limits_{j=1}^m (1-B(\boldsymbol{S},j))^{-1}.
\end{equation}
The $n$th moment of the total number of waiting jobs in the system, $\tilde{Q}$, under the c.o.s.\ mechanism is given by
\begin{equation}
\mathbb{E}\left[ \left((1-\rho)\tilde{Q}\right)^n \right] = \frac{n!(1-\rho)^n}{g(1)}\sum\limits_{L=0}^N\sum\limits_{\boldsymbol{u}\in\mathcal{N}_L}\prod\limits_{l=1}^L \frac{\mu_{u_l}}{\lambda_{\mathcal{C}}(u_1,\dots,u_l)} \sum\limits_{m=1}^{|\mathcal{S}|}\sum\limits_{\boldsymbol{S}\in \mathcal{S}_m^{\boldsymbol{u}}} \gamma(\boldsymbol{S}) \prod\limits_{j=1}^m A(\boldsymbol{S},j)\prod\limits_{j=1}^m (1-B(\boldsymbol{S},j))^{-1}.
\end{equation}
For any  $\boldsymbol{S}\in\mathcal{S}_m$ and $j=1,\dots,m$ we define
\begin{equation}
\gamma(\boldsymbol{S}) = \sum\limits_{\boldsymbol{k}\colon k_0+\dots + k_m=n} \frac{m^{k_0}}{k_0!}\prod\limits_{i=1}^m \sum\limits_{\boldsymbol{m}\in R(k_i)}\binom{m_1+\dots+m_{k_i}}{m_1,\dots,m_{k_i}} \left(1-B(\boldsymbol{S},i)\right)^{-|\boldsymbol{m}|}
\prod\limits_{j=1}^{k_i}\frac{B(\boldsymbol{S},i)^{m_j}}{j!^{m_j}},
\end{equation}
\begin{equation}\label{eq:A_S_j}
A(\boldsymbol{S},j) \coloneqq  \frac{N\lambda p_{S_j}}{\mu(S_1,\dots,S_j)},
\end{equation}
and
\begin{equation}\label{eq:B_S_j}
B(\boldsymbol{S},j) \coloneqq  \frac{N\lambda\sum_{i=1}^j p_{S_i}}{\mu(S_1,\dots,S_j)}.
\end{equation}
The functions $f$ and $g$ are defined in Proposition~\ref{prop:pgf} and~\ref{prop:pgfcos}, respectively.
\end{lemma}

\begin{proof}[Proof of Lemma~\ref{lem:nth_moment_total}.]
\textbf{C.o.c.\ mechanism:}\\
\textit{Step~1: Find a general expression for $\mathbb{E}\left[ \left((1-\rho)Q\right)^n \right]$.}
From Proposition~\ref{prop:pgf} we have that 
\begin{equation}
\mathbb{E}\left[ \mathrm{e}^{-(1-\rho)t Q}\right] = \frac{f(t)}{f(1)},
\end{equation}
with
\begin{equation}
f(t) \coloneqq 1+\sum\limits_{m=1}^{|\mathcal{S}|}\sum\limits_{\boldsymbol{S}\in \mathcal{S}_m} \left(\prod\limits_{j=1}^m \frac{N\lambda p_{S_j}}{\mu(S_1,\dots,S_j)}\right) \mathrm{e}^{-(1-\rho)mt}\left(\prod\limits_{j=1}^m\left(1{-}\frac{N\lambda\sum_{i=1}^j p_{S_i}}{\mu(S_1,\dots,S_j)} \mathrm{e}^{-(1-\rho)t}\right)^{-1}\right),
\end{equation}
with $t\ge 0$.
Hence, taking the appropriate derivatives we will obtain for $n\ge 1$,
\begin{equation}\label{proof:general derivative}
\mathbb{E}\left[ \left((1-\rho)Q\right)^n \right] = \frac{(-1)^n}{f(1)} \sum\limits_{m=1}^{|\mathcal{S}|}\sum\limits_{\boldsymbol{S}\in \mathcal{S}_m} \left(\prod\limits_{j=1}^m A(\boldsymbol{S},j)\right) \frac{\mathrm{d}^n}{\mathrm{d}t^n} h(\boldsymbol{S},t)\big{\vert}_{t=0},
\end{equation}
with
\begin{equation}\label{eq:h_general}
h(\boldsymbol{S},t) \coloneqq  \mathrm{e}^{-(1-\rho)mt}\prod\limits_{j=1}^m\left(1-B(\boldsymbol{S},j)\mathrm{e}^{-(1-\rho)t}\right)^{-1}.
\end{equation}
Hence, we need to determine the $n$th derivative of $h(\boldsymbol{S},t)$ with respect to $t$. Note that $A(\boldsymbol{S},j)$, $B(\boldsymbol{S},j)$ and $h(\boldsymbol{S},t)$ also depend on $\lambda$, but this is omitted to not make the notation too complicated and for now $\lambda$ is also assumed to be fixed. \\
\textit{Step 2: Rewrite $h(\boldsymbol{S},t)$.}\\
Let 
\begin{equation}
h(\boldsymbol{S},t) = \prod\limits_{i=0}^m h_i(\boldsymbol{S},t),
\end{equation}
with 
\begin{equation}
h_i(\boldsymbol{S},t) = \left\{
\begin{array}{lcl}
\mathrm{e}^{-(1-\rho)mt}, & & \text{if $i=0$}\\
\left(1-B(\boldsymbol{S},i)\mathrm{e}^{-(1-\rho)t}\right)^{-1}, & & \text{if $i=1,\dots,m$.}
\end{array}
\right.
\end{equation}
\textit{Step 3: Apply the general Leibniz rule.}
\begin{equation}
\frac{\mathrm{d}^n}{\mathrm{d}t^n} h(\boldsymbol{S},t) = \sum\limits_{\boldsymbol{k}\colon k_0+\dots + k_m=n} \binom{n}{k_0,\dots , k_m} \prod_{i=0}^m h_i^{(k_i)}(\boldsymbol{S},t).
\end{equation}
\textit{Step 4: Determine $h_i^{(k_i)}(\boldsymbol{S},t)$.}\\
\begin{equation}
h_0^{(k_0)}(\boldsymbol{S},t) = \left(-(1-\rho)m\right)^{k_0}\mathrm{e}^{-(1-\rho)mt}.
\end{equation}
For $i=1,\dots,k$ we define the function $v(t) \coloneqq t^{-1}$ and
\begin{equation}
g_i(\boldsymbol{S},t) \coloneqq 1-B(\boldsymbol{S},i)\mathrm{e}^{-(1-\rho)t}
\end{equation}
such that $h_i(\boldsymbol{S},t) = v\left(g_i(\boldsymbol{S},t) \right)$. Now we apply Fa\`a di Bruno's formula,
\begin{equation}
 h_i^{(k_i)}(\boldsymbol{S},t) = \sum\limits_{\boldsymbol{m}\in R(k_i)}   \frac{k_i!}{m_1!m_2!\dots m_{k_i}!} v^{(|\boldsymbol{m}|)}\left(g_i(\boldsymbol{S},t)\right)\prod\limits_{j=1}^{k_i} \left(\frac{g^{(j)}_i(\boldsymbol{S},t)}{j!}\right)^{m_j},
\end{equation}
with $|\boldsymbol{m}| = m_1+\dots + m_{k_i}$.
We observe that
\begin{equation}
\begin{array}{rcl}
v^{(|\boldsymbol{m}|)}(t) & =& (-1)^{|\boldsymbol{m}|}\cdot |\boldsymbol{m}|! \cdot t^{-|\boldsymbol{m}|-1}\\
g^{(j)}_i(\boldsymbol{S},t)&=& (-1)^{j+1}B(\boldsymbol{S},i) (1-\rho)^j\mathrm{e}^{-(1-\rho)t}.
\end{array}
\end{equation}
\textit{Step 5: Combine the above steps.}
\begin{equation}
\begin{array}{rcl}
h^{(n)}(\boldsymbol{S},t)  &=& \sum\limits_{\boldsymbol{k}\colon k_0+\dots + k_m=n} \binom{n}{k_0,\dots,k_m}\prod\limits_{i=0}^m h_i^{(k_i)}(\boldsymbol{S},t)\\

&=& \sum\limits_{\boldsymbol{k}\colon k_0+\dots + k_m=n} \frac{n!}{k_0!k_1!\dots k_m!}\left[(-1)(1-\rho)m\right]^{k_0}\mathrm{e}^{-(1-\rho)mt}  \\
& & \times\prod\limits_{i=1}^{m} \left\{\sum\limits_{\boldsymbol{m}\in R(k_i)} \frac{k_i!}{m_1!m_2!\dots m_{k_i}!} v^{(|\boldsymbol{m}|)}\left(g_i(\boldsymbol{S},t)\right) \prod\limits_{j=1}^{k_i} \left(\frac{g_i^{(j)}(\boldsymbol{S},t)}{j!}\right)^{m_j} \right\} \\

&=& \sum\limits_{\boldsymbol{k}\colon k_0+\dots + k_m=n} \frac{n!}{k_0!}\left[(-1)(1-\rho)m\right]^{k_0}\mathrm{e}^{-(1-\rho)mt} \\
& &  \times\prod\limits_{i=1}^{m} \left\{\sum\limits_{\boldsymbol{m}\in R(k_i)} \frac{|\boldsymbol{m}|!}{m_1!m_2!\dots m_{k_i}!}(-1)^{|\boldsymbol{m}|} \left( 1- B(\boldsymbol{S},i)\mathrm{e}^{-(1-\rho)t}\right)^{-|\boldsymbol{m}|-1} \right.\\
& & \times \left. \prod\limits_{j=1}^{k_i} \left(\frac{(-1)^{j+1}B(\boldsymbol{S},i)(1-\rho)^j\mathrm{e}^{-(1-\rho)t}}{j!}\right)^{m_j} \right\} \\

& =& n!(1-\rho)^n(-1)^n \sum\limits_{\boldsymbol{k}\colon k_0+\dots + k_m=n} \frac{m^{k_0}}{k_0!}\mathrm{e}^{-(1-\rho)mt}\\

& & \times  \prod\limits_{i=1}^m \left\{  \sum\limits_{\boldsymbol{m}\in R(k_i)}   \binom{m_1+\dots+m_{k_i}}{m_1,\dots,m_{k_i}}\left(1-B(\boldsymbol{S},i) \mathrm{e}^{-(1-\rho)t}\right) ^{-|\boldsymbol{m}|-1} \right.\\

& & \times \left.  \prod\limits_{j=1}^{k_i}\frac{B(\boldsymbol{S},i)^{m_j}\mathrm{e}^{-(1-\rho)m_jt}}{j!^{m_j}} \right\}.
\end{array}
\end{equation}
Substitute $t=0$ in the above expression and plug it into~\eqref{proof:general derivative} to obtain the desired result.\\
\textbf{C.o.s.\ mechanism:} Considering Proposition~\ref{prop:pgfcos}, we note that
\begin{equation}\label{proof:general derivative cos}
\mathbb{E}\left[ \left((1-\rho)\tilde{Q}\right)^n \right] = \frac{(-1)^n n!(1-\rho)^n}{g(1)}\sum\limits_{L=0}^N\sum\limits_{\boldsymbol{u}\in\mathcal{N}_L}\prod\limits_{l=1}^L \frac{\mu_{u_l}}{\lambda_{\mathcal{C}}(u_1,\dots,u_l)}
 \sum\limits_{m=1}^{|\mathcal{S}|}\sum\limits_{\boldsymbol{S}\in \mathcal{S}_m^{\boldsymbol{u}}} \left(\prod\limits_{j=1}^m A(\boldsymbol{S},j)\right) \frac{\mathrm{d}^n}{\mathrm{d}t^n} h(\boldsymbol{S},t)\big{\vert}_{t=0},
\end{equation}
with $h(\boldsymbol{S},t)$ as defined in~\eqref{eq:h_general}. Hence, following the steps applied for the c.o.c.\ mechanism and combining all results will lead to the desired expression.
\end{proof}

\begin{proof}[Proof of Theorem~\ref{th:conv_moments_tot}.] 
We will study the heavy-traffic limit of the $n$th moments as $\lambda\uparrow\mu$ obtained in Lemma~\ref{lem:nth_moment_total} using similar arguments as the ones in the proof of Theorem~\ref{statespacecollapse_general}. In the proof below we will focus on the c.o.c.\ mechanism; the proof for the c.o.s.\ version follows the same lines.\\ \textbf{C.o.c.\ mechanism:}  Observe that for the c.o.c.\ version we have
\begin{equation}
\lim_{\lambda\uparrow\mu} \mathbb{E}\left[((1-\rho)Q)^n \right] = \lim_{\lambda\uparrow\mu} \mathbb{E}\left[((1-\rho)Q)^n \right].
\end{equation}
Before computing the overall limit we first focus on the numerator and denominator (multiplied by $(1-\rho)$) of the above expression separately. \\
\textit{Denominator:}
\begin{equation}
\begin{array}{rl}
& (1-\rho)f(1)\\

= & (1-\rho) + (1-\rho)\sum\limits_{m=1}^{|\mathcal{S}|-1}\sum\limits_{\boldsymbol{S}\in \mathcal{S}_m}\prod\limits_{j=1}^m A(\boldsymbol{S},j)\prod\limits_{j=1}^m\left(1-B(\boldsymbol{S},j) \right)^{-1}\\
& + \sum\limits_{\boldsymbol{S}\in \mathcal{S}_{|\mathcal{S}|}}\prod\limits_{j=1}^{|\mathcal{S}|} A(\boldsymbol{S},j)\prod\limits_{j=1}^{|\mathcal{S}|-1}\left(1-B(\boldsymbol{S},j) \right)^{-1}\\

= & \sum\limits_{\boldsymbol{S}\in \mathcal{S}_{|\mathcal{S}|}}\prod\limits_{j=1}^{|\mathcal{S}|}A(\boldsymbol{S},j)\prod\limits_{j=1}^{|\mathcal{S}|-1}\left(1-B(\boldsymbol{S},j)\right)^{-1} + o(1-\rho).
\end{array}
\end{equation}
The last equality is due to the local stability conditions~\eqref{stabcond3} and~\eqref{stabcond4} as 
\begin{equation}
\lim_{\lambda\uparrow\mu}\frac{N\lambda\sum_{i=1}^j p_{S_i}}{\mu(S_1,\dots,S_j)} < 1
\end{equation}
for any $\boldsymbol{S}\in \mathcal{S}_m$, $j=1,\dots,m$ and $m=1,\dots,|\mathcal{S}|-1$.\\
\textit{Numerator:}
\begin{equation}
\begin{array}{rl}
& n!(1-\rho)^{n+1} \sum\limits_{m=1}^{|\mathcal{S}|}\sum\limits_{\boldsymbol{S}\in \mathcal{S}_m} \prod\limits_{j=1}^m A(\boldsymbol{S},j)\prod\limits_{j=1}^m (1-B(\boldsymbol{S},j))^{-1}  \\
& \times \sum\limits_{\boldsymbol{k}\colon k_0+\dots + k_m=n} \frac{m^{k_0}}{k_0!}\left\{ \prod\limits_{i=1}^m \sum\limits_{\boldsymbol{m}\in R(k_i)}\frac{|\boldsymbol{m}|!}{m_1!\dots m_{k_i}!} \left(1-B(\boldsymbol{S},i)\right)^{-|\boldsymbol{m}|}
\prod\limits_{j=1}^{k_i}\frac{B(\boldsymbol{S},i)^{m_j}}{j!^{m_j}} \right\}\\

=& n!(1-\rho)^{n+1} \sum\limits_{m=1}^{|\mathcal{S}|-1}\sum\limits_{\boldsymbol{S}\in \mathcal{S}_m} \prod\limits_{j=1}^m A(\boldsymbol{S},j)\prod\limits_{j=1}^m (1-B(\boldsymbol{S},j))^{-1} \\
&  \times \sum\limits_{\boldsymbol{k}\colon k_0+\dots + k_m=n} \frac{m^{k_0}}{k_0!} \left\{ \prod\limits_{i=1}^m \sum\limits_{\boldsymbol{m}\in R(k_i)}\frac{|\boldsymbol{m}|!}{m_1!\dots m_{k_i}!} \left(1-B(\boldsymbol{S},i)\right)^{-|\boldsymbol{m}|}
\prod\limits_{j=1}^{k_i}\frac{B(\boldsymbol{S},i)^{m_j}}{j!^{m_j}} \right\}\\
+& n!(1-\rho)^n \sum\limits_{\boldsymbol{S}\in \mathcal{S}_{|\mathcal{S}|}} \prod\limits_{j=1}^{|\mathcal{S}|} A(\boldsymbol{S},j)\prod\limits_{j=1}^{|\mathcal{S}|-1} (1-B(\boldsymbol{S},j))^{-1}  \\
& \times \sum\limits_{\boldsymbol{k}\colon k_0+\dots + k_{|\mathcal{S}|}=n} \frac{{|\mathcal{S}|}^{k_0}}{k_0!}\prod\limits_{i=1}^{|\mathcal{S}|}\left\{  \sum\limits_{\boldsymbol{m}\in R(k_i)}\frac{|\boldsymbol{m}|!}{m_1!\dots m_{k_i}!} \left(1-B(\boldsymbol{S},i)\right)^{-|\boldsymbol{m}|}
\prod\limits_{j=1}^{k_i}\frac{B(\boldsymbol{S},i)^{m_j}}{j!^{m_j}} \right\}.
\end{array}
\end{equation}
The first summation vanishes in the limit due to the prefactor $(1-\rho)^{n+1}$ and the fact that all its terms are finite as the local stability conditions~\eqref{stabcond3} and~\eqref{stabcond4} are satisfied. We get
\begin{equation}
\begin{array}{l}
 n!(1-\rho)^n \sum\limits_{\boldsymbol{S}\in \mathcal{S}_{|\mathcal{S}|}} \prod\limits_{j=1}^{|\mathcal{S}|} A(\boldsymbol{S},j)\prod\limits_{j=1}^{|\mathcal{S}|-1} (1-B(\boldsymbol{S},j))^{-1}  \\
 
\times \sum\limits_{\boldsymbol{k}\colon k_0+\dots + k_{|\mathcal{S}|}=n} \frac{{|\mathcal{S}|}^{k_0}}{k_0!}\prod\limits_{i=1}^{|\mathcal{S}|-1}\left\{  \sum\limits_{\boldsymbol{m}\in R(k_i)}\frac{|\boldsymbol{m}|!}{m_1!\dots m_{k_i}!} \left(1-B(\boldsymbol{S},i)\right)^{-|\boldsymbol{m}|}
\prod\limits_{j=1}^{k_i}\frac{B(\boldsymbol{S},i)^{m_j}}{j!^{m_j}} \right\} \\

\times\left\{  \sum\limits_{\boldsymbol{m}\in R(k_{|\mathcal{S}|})}\frac{|\boldsymbol{m}|!}{m_1!\dots m_{k_{|\mathcal{S}|}}!} \left(1-\rho\right)^{-|\boldsymbol{m}|}
\prod\limits_{j=1}^{k_{|\mathcal{S}|}}\frac{B(\boldsymbol{S},i)^{m_j}}{j!^{m_j}} \right\} + o(1-\rho).
\end{array}
\end{equation}
Note that if $\boldsymbol{m}\in R(k_{|\mathcal{S}|})$ such that $n-|\boldsymbol{m}| > 0$, the corresponding term vanishes as $\lambda\uparrow\mu$ due to the presence of the factor $(1-\rho)^{n-|\boldsymbol{m}|}$. So, we only need to focus on those $\boldsymbol{m}\in R(k_{|\mathcal{S}|})$ such that $n=|\boldsymbol{m}|$ and $1m_1+2m_2+\dots, k_{|\mathcal{S}|}m_{k_{|\mathcal{S}|}} = k_{|\mathcal{S}|}$. Moreover, we know that $k_0+k_1+\dots k_{|\mathcal{S}|} = n$. This implies that
\begin{equation}
n = |\boldsymbol{m}| \le 1m_1+2m_2+\dots, k_{|\mathcal{S}|}m_{k_{|\mathcal{S}|}} = k_{|\mathcal{S}|} \le n,
\end{equation} 
or $k_0 = k_1 = \dots  = k_{|\mathcal{S}|-1} = 0$ and $k_{|\mathcal{S}|} = n$. From this we deduce that the non-vanishing terms in the numerator are such that $\boldsymbol{m}\in R(k_i)$ for $i=0,\dots,|\mathcal{S}|-1$ are the all-zeros vectors and that $\boldsymbol{m}\in R(k_{|\mathcal{S}|})$ is such that $m_1 = n$ and $m_2 = \dots = m_{k_{|\mathcal{S}|}} = 0$. This results in the following simplification for the numerator:
\begin{equation}
\begin{array}{rl}
&  n!(1-\rho)^n \sum\limits_{\boldsymbol{S}\in \mathcal{S}_{|\mathcal{S}|}} \prod\limits_{j=1}^{|\mathcal{S}|} A(\boldsymbol{S},j)\prod\limits_{j=1}^{|\mathcal{S}|-1} (1-B(\boldsymbol{S},j))^{-1}  \left(1-\rho\right)^{-n}
\frac{\rho^{n}}{1!^{n}}  + o(1-\rho) \\

= & n!\rho^n \sum\limits_{\boldsymbol{S}\in \mathcal{S}_{|\mathcal{S}|}} \prod\limits_{j=1}^{|\mathcal{S}|} A(\boldsymbol{S},j)\prod\limits_{j=1}^{|\mathcal{S}|-1} (1-B(\boldsymbol{S},j))^{-1}  + o(1-\rho). 
\end{array}
\end{equation}
Combined this gives us
\begin{equation}
\lim\limits_{\lambda\uparrow\mu} \mathbb{E}\left[((1-\rho)Q)^n \right]  = 
\lim\limits_{\lambda\uparrow\mu}\frac{n!\rho^n \sum\limits_{\boldsymbol{S}\in \mathcal{S}_{|\mathcal{S}|}} \prod\limits_{j=1}^{|\mathcal{S}|} A(\boldsymbol{S},j)\prod\limits_{j=1}^{|\mathcal{S}|-1} (1-B(\boldsymbol{S},j))^{-1}  + o(1-\rho)}{ \sum\limits_{\boldsymbol{S}\in \mathcal{S}_{|\mathcal{S}|}} \prod\limits_{j=1}^{|\mathcal{S}|} A(\boldsymbol{S},j)\prod\limits_{j=1}^{|\mathcal{S}|-1} (1-B(\boldsymbol{S},j))^{-1}  + o(1-\rho)} =  n!.
\end{equation}
\textbf{C.o.s.\ mechanism:}  Via a similar reasoning, where the dominant terms in the limit are those including all job types, we can show that 
\begin{equation}
\lim_{\lambda\uparrow\mu} \mathbb{E}\left[\left((1-\rho)\tilde{Q}\right)^n \right] = n!.
\end{equation}
In order to prove the convergence result in Theorem~\ref{th:conv_moments_tot} in terms of the total number of jobs in the system, we use that $\tilde{Q}\le Q\le \tilde{Q}+N$. Hence, 
\begin{equation}
\mathbb{E}\left[((1-\rho)\tilde{Q})^n \right] \le \mathbb{E}\left[((1-\rho)Q)^n \right] \le  \mathbb{E}\left[\left((1-\rho)(\tilde{Q}+N)\right)^n \right],
\end{equation}
where the last term can be rewritten as
\begin{equation}
\mathbb{E}\left[((1-\rho)\tilde{Q})^n \right] + \sum\limits_{k=0}^{n-1} \binom{n}{k} N^{n-k}(1-\rho)^{n-k} \mathbb{E}\left[((1-\rho)\tilde{Q})^k \right] = \mathbb{E}\left[((1-\rho)\tilde{Q})^n \right] + o(1-\rho)
\end{equation}
using an inductive argument that the $k$th moment of $(1-\rho)\tilde{Q}$ converges to $k!$, for $k<n$. This concludes the proof.
\end{proof}

\subsection{The moments of the total number of jobs of a fixed type}
\begin{theorem}\label{th:conv_moments_type}
If the local stability condition~\eqref{stabcond3} and~\eqref{stabcond4} hold, then for both the c.o.c.\ and c.o.s. mechanisms
\begin{equation}
\lim_{\lambda\uparrow\mu}\mathbb{E}\left[ \left((1-\rho)Q_T\right)^n \right] = n!p_T^n,
\end{equation} 
for any $n\ge 1$ and any job type $T\in\mathcal{S}$.
\end{theorem}

It can again be deduced, via the method of moments, that $(1-\rho)Q_T$ converges in distribution to $\mathrm{Exp}(1)p_T$. This result is weaker than the result in Theorem~\ref{statespacecollapse} where we prove convergence in distribution of the full random vector $(1-\rho)(Q_S)_{S\in\mathcal{S}}$.\\

To prove Theorem~\ref{th:conv_moments_tot} we first derive expressions for the $n$th moments of $Q$ and $\tilde{Q}$ for any fixed value of $\lambda$ within the stability region for the c.o.c.\ and c.o.s.\ mechanisms, respectively.

\begin{lemma}\label{lem:nth_moment_type}
Assuming that the stability conditions~\eqref{stabcond1} and~\eqref{stabcond2} are satisfied, the $n$th moment of the number of type-$T$ jobs in the system, $Q_T$ under the c.o.c.\ mechanism is given by
\begin{equation}
\mathbb{E}\left[ \left((1-\rho)Q_T\right)^n \right]  = \frac{n!(1-\rho)^n}{f(1)} \sum\limits_{m=1}^{|\mathcal{S}|}\sum\limits_{\boldsymbol{S}\in \mathcal{S}_m(T)} \gamma(\boldsymbol{S})\prod\limits_{j=1}^m A(\boldsymbol{S},j)\prod\limits_{j=1}^m (1-B(\boldsymbol{S},j))^{-1}.
\end{equation}
The $n$th moment of the number of waiting type-$T$ jobs in the system, $\tilde{Q}_T$ under the c.o.s.\ mechanism is given by
\begin{equation}
\mathbb{E}\left[ \left((1-\rho)\tilde{Q}_T\right)^n \right] = \frac{n!(1-\rho)^n}{g(1)} \sum\limits_{L=0}^N \sum\limits_{\boldsymbol{u}\in\mathcal{N}_L}\prod\limits \frac{\mu_{u_l}}{\lambda_{\mathcal{C}}(u_1,\dots,u_l)}\sum\limits_{m=1}^{|\mathcal{S}|}\sum\limits_{\boldsymbol{S}\in \mathcal{S}_m(T)} \gamma(\boldsymbol{S}) \prod\limits_{j=1}^m A(\boldsymbol{S},j) \prod\limits_{j=1}^m (1-B(\boldsymbol{S},i))^{-1}, 
\end{equation}
with 
\begin{equation}
\gamma(\boldsymbol{S})\coloneqq  \sum\limits_{\boldsymbol{k}\in\tilde{R}(n)}  
\frac{1}{k_0!} \prod\limits_{i=l}^m \left\{
\sum\limits_{\boldsymbol{m}\in R(k_i)}   \frac{(|\boldsymbol{m}|)!}{m_1!m_2!\dots m_{k_i}!} \left(1-B(\boldsymbol{S},i)\right)^{-|\boldsymbol{m}|}\prod\limits_{j=1}^{k_i} \left(\frac{ N\lambda p_{S_l}}{j!\mu(S_1,\dots,S_i)}\right)^{m_j}
\right\}.
\end{equation}
We define $\mathcal{S}_m(T) \subseteq \mathcal{S}_m$ as the set of all $m$-dimensional ordered vectors $\boldsymbol{S}$ with unique elements from $\mathcal{S}$ and containing job type~$T$. We write $|\boldsymbol{m}|\coloneqq m_0+m_1+\dots+m_{k_i}$ for  $\boldsymbol{m}\in R(k_i)$, with $R(k_i)$ as in~\eqref{eq:R(n)}. Moreover, $\tilde{R}(n) \subset R(n)$ such that the entries $k_1,\dots,k_{l-1}=0$ for $k\in\tilde{R}(n)$, with $l$ the position of type $T$ in the ordered vector $\boldsymbol{S}$.
Define $A(\boldsymbol{S},i)$ and $B(\boldsymbol{S},i)$ as in~\eqref{eq:A_S_j} and~\eqref{eq:B_S_j}, respectively.
\end{lemma}

\begin{proof}[Proof of Lemma~\ref{lem:nth_moment_type}.]
The proof follows  similar steps as the proof of Lemma~\ref{lem:nth_moment_total}.\\
\textbf{C.o.c.\ mechanism:}\\
\textit{Step 1: Find a general expression for $\mathbb{E}\left[ \left((1-\rho)Q_T\right)^n \right]$.}
From Proposition~\ref{prop:pgf} we have that 
\begin{equation}
\mathbb{E}\left[ \mathrm{e}^{-(1-\rho)t Q_T}  \right] = \mathbb{E}\left[ \prod\limits_{S\in\mathcal{S}}z_S^{Q_S}\right] = \frac{f(\boldsymbol{z})}{f(\boldsymbol{1})},
\end{equation}
with $t\ge 0$ and $z_S =  \mathrm{e}^{-(1-\rho)t}$ if $S= T$ and 1 otherwise. The function~$f$ is defined in~\eqref{eq:pgf_f}. Now, let 
\begin{equation}\label{eq:h(S,t)_type}
h(\boldsymbol{S},t) \coloneqq \prod\limits_{j=1}^m z_{S_j}
\prod\limits \left(1-\frac{N\lambda\sum_{i=1}^jp_{S_i}z_{S_i}}{\mu(S_1,\dots,S_j)}\right)^{-1}.
\end{equation}
Hence, taking the appropriate derivatives we will obtain for $n\ge 1$,
\begin{equation}\label{proof:general derivative_type}
\mathbb{E}\left[ \left((1-\rho)Q_T\right)^n \right] = \frac{(-1)^n}{f(\boldsymbol{1})} \sum\limits_{m=1}^{|\mathcal{S}|}\sum\limits_{\boldsymbol{S}\in \mathcal{S}_m} \left(\prod\limits_{j=1}^m A(\boldsymbol{S},j)\right) \frac{\mathrm{d}^n}{\mathrm{d}t^n} h(\boldsymbol{S},t)\big{\vert}_{t=0}.
\end{equation}
Hence, we need to determine the $n$th derivative if $h(\boldsymbol{S},t)$ with respect to $t$.\\
\textit{Step 2: Rewrite $h(\boldsymbol{S},t)$.}
Whenever $T\notin\boldsymbol{S}$, the $n$th derivative of $h$ will always be 0. Hence, we assume in the derivations below that $T\in\boldsymbol{S}$, more precisely at position $l\in\{1,\dots,m\}$. So, $z_{S_l} = \mathrm{e}^{-(1-\rho)t}$ and $z_{S_j} =1$ for $j\neq l$.
Let 
\begin{equation}
h(\boldsymbol{S},t) = \prod\limits_{i=0}^m h_i(\boldsymbol{S},t),
\end{equation}
with 
\begin{equation}
h_i(\boldsymbol{S},t) = \left\{
\begin{array}{lcl}
\prod\limits_{j=1}^m z_{S_j} = z_{S_l} = \mathrm{e}^{-(1-\rho)t}, & & \text{if $i=0$}\\
\left(1-\frac{N\lambda\sum_{i=1}^jp_{S_i}z_{S_i}}{\mu(S_1,\dots,S_j)}\right)^{-1}, & & \text{if $i=1,\dots,m$.}
\end{array}
\right.
\end{equation}
\textit{Step 3: Apply the general Leibniz rule.}
\begin{equation}
\frac{\mathrm{d}^n}{\mathrm{d}t^n} h(\boldsymbol{S},t) = \sum\limits_{\boldsymbol{k}\colon k_0,\dots ,k_m=n} \binom{n}{k_0+\dots + k_m} \prod_{i=0}^m h_i^{(k_i)}(\boldsymbol{S},t).
\end{equation}
\textit{Step 4: Determine $h_i^{(k_i)}(\boldsymbol{S},t)$.}\\
\begin{equation}
h_0^{(k_0)}(\boldsymbol{S},t) = \left(-(1-\rho)\right)^{k_0}\mathrm{e}^{-(1-\rho)t}.
\end{equation}
For $i=1,\dots,k$ we define the function $v(t) \coloneqq t^{-1}$ and
\begin{equation}
g_i(\boldsymbol{S},t) \coloneqq 1-\frac{N\lambda\sum_{i=1}^jp_{S_i}z_{S_i}}{\mu(S_1,\dots,S_j)}
\end{equation}
such that $h_i(\boldsymbol{S},t) = v\left(g_i(\boldsymbol{S},t) \right)$. Now we apply Fa\`a di Bruno's formula,
\begin{equation}
 h_i^{(k_i)}(\boldsymbol{S},t) = \sum\limits_{\boldsymbol{m}\in R(k_i)}   \frac{k_i!}{m_1!m_2!\dots m_{k_i}!} v^{(|\boldsymbol{m}|)}\left(g_i(\boldsymbol{S},t)\right)\prod\limits_{j=1}^{k_i} \left(\frac{g^{(j)}_i(\boldsymbol{S},t)}{j!}\right)^{m_j},
\end{equation}
with $|\boldsymbol{m}| = m_1+\dots + m_{k_i}$. We observe that $g^{(j)}_i(\boldsymbol{S},t) = 0$ if $i<l$ as $g_i$ is in this case independent of $t$. If $i\ge l$, we obtain
\begin{equation}
g^{(j)}_i(\boldsymbol{S},t) = \frac{(-1)^{j+1}(1-\rho)^j N\lambda p_{S_l}}{\mu(S_1,\dots,S_i)}\mathrm{e}^{-(1-\rho)t}.
\end{equation}
\textit{Step 5: Combine the above steps.}
\begin{equation}
\begin{array}{rcl}
\frac{\mathrm{d}^n}{\mathrm{d}t^n} h(\boldsymbol{S},t)\big{\vert}_{t=0}& = &  \sum\limits_{\boldsymbol{k}\colon k_0,\dots ,k_m=n} \binom{n}{k_0+\dots + k_m} 
\left(-(1-\rho)\right)^{k_0} \\

& & \times\prod\limits_{i=1}^m \left\{
\sum\limits_{\boldsymbol{m}\in R(k_i)}   \frac{k_i!(-1)^{|\boldsymbol{m}|}(|\boldsymbol{m}|)!}{m_1!m_2!\dots m_{k_i}!} \left(1-B(\boldsymbol{S},i)\right)^{-1-|\boldsymbol{m}|}\prod\limits_{j=1}^{k_i} \left(\frac{g^{(j)}_i(\boldsymbol{S},0)}{j!}\right)^{m_j}
\right\}\\

& = & n! \sum\limits_{\boldsymbol{k}\colon k_0,\dots ,k_m=n}  
\frac{1}{k_0!}\left(-(1-\rho)\right)^{k_0} \\

& &\times \prod\limits_{i=1}^m \left\{
\sum\limits_{\boldsymbol{m}\in R(k_i)}   \frac{(-1)^{|\boldsymbol{m}|}(|\boldsymbol{m}|)!}{m_1!m_2!\dots m_{k_i}!} \left(1-B(\boldsymbol{S},i)\right)^{-1-|\boldsymbol{m}|}\prod\limits_{j=1}^{k_i} \left(\frac{g^{(j)}_i(\boldsymbol{S},0)}{j!}\right)^{m_j}
\right\}.
\end{array}
\end{equation}
Observe that for any $\boldsymbol{k}$ for which there exists an $i\in\{1,\dots,l-1\}$ such that $k_i >0$ the whole corresponding term is equal to 0. This is due to the fact that $g_i^{(j)}(\boldsymbol{S},0) = 0$ for any $j\ge 1$. Hence, we only need to consider the natural vectors $\boldsymbol{k}$ such that $k_1=\dots = k_{l-1} = 0$ and $k_0+k_l+\dots+k_m = n$. We depict this by $\boldsymbol{k}\in\tilde{R}(\boldsymbol{S},n)$. Hence,
\begin{equation}
\begin{array}{rcl}
\frac{\mathrm{d}^n}{\mathrm{d}t^n} h(\boldsymbol{S},t)\big{\vert}_{t=0}& = &  n! \sum\limits_{\boldsymbol{k}\in\tilde{R}(\boldsymbol{S},n)}  
\frac{1}{k_0!}\left(-(1-\rho)\right)^{k_0} \\

& &\times \prod\limits_{i=l}^m \left\{
\sum\limits_{\boldsymbol{m}\in R(k_i)}   \frac{(-1)^{|\boldsymbol{m}|}(|\boldsymbol{m}|)!}{m_1!m_2!\dots m_{k_i}!} \left(1-B(\boldsymbol{S},i)\right)^{-1-|\boldsymbol{m}|}\prod\limits_{j=1}^{k_i} \left(\frac{(-1)^{j+1}(1-\rho)^j N\lambda p_{S_l}}{j!\mu(S_1,\dots,S_i)}\right)^{m_j}
\right\}\\

& = &  n!(-1)^n (1-\rho)^n \sum\limits_{\boldsymbol{k}\in\tilde{R}(\boldsymbol{S},n)}  
\frac{1}{k_0!} \\
& & \times\prod\limits_{i=l}^m \left\{
\sum\limits_{\boldsymbol{m}\in R(k_i)}   \frac{(|\boldsymbol{m}|)!}{m_1!m_2!\dots m_{k_i}!} \left(1-B(\boldsymbol{S},i)\right)^{-1-|\boldsymbol{m}|}\prod\limits_{j=1}^{k_i} \left(\frac{ N\lambda p_{S_l}}{j!\mu(S_1,\dots,S_i)}\right)^{m_j}
\right\}.
\end{array}
\end{equation}
Plugging the above expression into~\eqref{proof:general derivative_type} concludes the proof for the c.o.c.\ mechanism.\\
\textbf{C.o.s.\ mechanism:} We observe that~\eqref{proof:general derivative cos} still holds with $h(\boldsymbol{S},t)$ as defined in~\eqref{eq:h(S,t)_type}. Following the steps outlined above will result in an expression for the $n$th moment under de c.o.s.\ mechanism. This concludes the proof.
\end{proof}

\begin{proof}[Proof of Theorem~\ref{th:conv_moments_type}.]
We will study the heavy-traffic limit of the $n$th moments as $\lambda\uparrow\mu$ obtained in Lemma~\ref{lem:nth_moment_type} using similar arguments as the ones in the proof of Theorem~\ref{statespacecollapse_general}. \\
\textbf{C.o.c. mechanism:} 
Write
\begin{equation}
\lim_{\lambda\uparrow\mu}\mathbb{E}\left[ \left((1-\rho)Q_T\right)^n \right] = \lim_{\lambda\uparrow\mu}\mathbb{E}\left[ \left((1-\rho)Q_T\right)^n \right]\frac{(1-\rho)}{(1-\rho)}.
\end{equation} 
and consider first the numerator and denominator (multiplied by $(1-\rho)$) separately before computing the actual limit.\\
\textit{Denominator:}
From the proof of Theorem~\ref{th:conv_moments_tot}, we know via the local stability conditions~\eqref{stabcond3} and~\eqref{stabcond4} that
\begin{equation}
(1-\rho)f(\boldsymbol{1})= \sum\limits_{\boldsymbol{S}\in \mathcal{S}_{|\mathcal{S}|}}\prod\limits_{j=1}^{|\mathcal{S}|} \frac{N\lambda p_{S_j}}{\mu(S_1,\dots,S_j)}\prod\limits_{j=1}^{|\mathcal{S}|-1}\left(1-\frac{N\lambda\sum_{i=1}^j p_{S_i}}{\mu(S_1,\dots,S_j)} \right)^{-1} + o(1-\rho).
\end{equation}

\textit{Numerator:}
\begin{equation}
\begin{array}{rl}
& n!(1-\rho)^{n+1} \sum\limits_{m=1}^{|\mathcal{S}|}\sum\limits_{\boldsymbol{S}\in \mathcal{S}_m(T)} \left(\prod\limits_{j=1}^m \frac{N\lambda p_{S_j}}{\mu(S_1,\dots,S_j)}\right) \\
&\times \sum\limits_{\boldsymbol{k}\in\tilde{R}(\boldsymbol{S},n)}  
\frac{1}{k_0!} \prod\limits_{i=l}^m \left\{
\sum\limits_{\boldsymbol{m}\in R(k_i)}   \frac{(|\boldsymbol{m}|)!}{m_1!m_2!\dots m_{k_i}!} \left(1-B(\boldsymbol{S},i)\right)^{-1-|\boldsymbol{m}|}\prod\limits_{j=1}^{k_i} \left(\frac{ N\lambda p_{S_l}}{j!\mu(S_1,\dots,S_i)}\right)^{m_j}\right\}\\

= & n!(1-\rho)^{n+1} \sum\limits_{m=1}^{|\mathcal{S}|-1}\sum\limits_{\boldsymbol{S}\in \mathcal{S}_m(T)} \left(\prod\limits_{j=1}^m \frac{N\lambda p_{S_j}}{\mu(S_1,\dots,S_j)}\right) \\
& \times\sum\limits_{\boldsymbol{k}\in\tilde{R}(\boldsymbol{S},n)} \frac{1}{k_0!} \prod\limits_{i=l}^m \left\{\sum\limits_{\boldsymbol{m}\in R(k_i)}   \frac{(|\boldsymbol{m}|)!}{m_1!m_2!\dots m_{k_i}!} \left(1-B(\boldsymbol{S},i)\right)^{-1-|\boldsymbol{m}|}\prod\limits_{j=1}^{k_i} \left(\frac{ N\lambda p_{S_l}}{j!\mu(S_1,\dots,S_i)}\right)^{m_j} \right\}\\
&+ n!(1-\rho)^{n}\sum\limits_{\boldsymbol{S}\in \mathcal{S}_{|\mathcal{S}|}(T)} \left(\prod\limits_{j=1}^{|\mathcal{S}|} \frac{N\lambda p_{S_j}}{\mu(S_1,\dots,S_j)}\right) \\
& \times \sum\limits_{\boldsymbol{k}\in\tilde{R}(\boldsymbol{S},n)} \frac{1}{k_0!} \prod\limits_{i=l}^{|\mathcal{S}|-1} \left\{\sum\limits_{\boldsymbol{m}\in R(k_i)}   \frac{(|\boldsymbol{m}|)!}{m_1!m_2!\dots m_{k_i}!} \left(1-B(\boldsymbol{S},i)\right)^{-1-|\boldsymbol{m}|}\prod\limits_{j=1}^{k_i} \left(\frac{ N\lambda p_{S_l}}{j!\mu(S_1,\dots,S_i)}\right)^{m_j} \right\} \\
&\times  \left\{\sum\limits_{\boldsymbol{m}\in R(k_{|\mathcal{S}|})}   \frac{(|\boldsymbol{m}|)!}{m_1!m_2!\dots m_{k_i}!} \left(1-\rho\right)^{-|\boldsymbol{m}|}\prod\limits_{j=1}^{k_{|\mathcal{S}|}} \left(\frac{ N\lambda p_{S_l}}{j!N\mu}\right)^{m_j} \right\}\\

=  & n!(1-\rho)^{n}\sum\limits_{\boldsymbol{S}\in \mathcal{S}_{|\mathcal{S}|}} \left(\prod\limits_{j=1}^{|\mathcal{S}|} \frac{N\lambda p_{S_j}}{\mu(S_1,\dots,S_j)}\right) \\
& \times\sum\limits_{\boldsymbol{k}\in\tilde{R}(\boldsymbol{S},n)} \frac{1}{k_0!} \prod\limits_{i=l}^{|\mathcal{S}|-1} \left\{\sum\limits_{\boldsymbol{m}\in R(k_i)}   \frac{(|\boldsymbol{m}|)!}{m_1!m_2!\dots m_{k_i}!} \left(1-B(\boldsymbol{S},i)\right)^{-1-|\boldsymbol{m}|}\prod\limits_{j=1}^{k_i} \left(\frac{ N\lambda p_{S_l}}{j!\mu(S_1,\dots,S_i)}\right)^{m_j} \right\}\\
&  \times \left\{\sum\limits_{\boldsymbol{m}\in R(k_{|\mathcal{S}|})}   \frac{(|\boldsymbol{m}|)!}{m_1!m_2!\dots m_{k_i}!} \left(1-\rho\right)^{-|\boldsymbol{m}|}\prod\limits_{j=1}^{k_{|\mathcal{S}|}} \left(\frac{ \lambda p_{S_l}}{j!\mu}\right)^{m_j} \right\} +o(1-\rho).
\end{array}
\end{equation}
The last equality follows again from the local stability conditions~\eqref{stabcond3} and~\eqref{stabcond4}, and the observation that $\mathcal{S}_{|\mathcal{S}|}(T) = \mathcal{S}_{|\mathcal{S}|}$. Note that if $\boldsymbol{m}\in R(k_{|\mathcal{S}|})$ such that $n-|\boldsymbol{m}| > 0$, the corresponding term vanishes as $\lambda\uparrow\mu$. So, we only need to focus on those $\boldsymbol{m}\in R(k_{|\mathcal{S}|})$ such that $n=|\boldsymbol{m}|$ and $1m_1+2m_2+\dots, k_{|\mathcal{S}|}m_{k_{|\mathcal{S}|}} = k_{|\mathcal{S}|}$. Moreover, we know that $k_0+k_1+\dots k_{|\mathcal{S}|} = n$. This implies that
\begin{equation}
n = |\boldsymbol{m}| \le 1m_1+2m_2+\dots, k_{|\mathcal{S}|}m_{k_{|\mathcal{S}|}} = k_{|\mathcal{S}|} \le n,
\end{equation} 
or $k_0 = k_1 = \dots  = k_{|\mathcal{S}|-1} = 0$ and $k_{|\mathcal{S}|} = n$. From this we deduce that the non-vanishing terms in the numerator are such that $\boldsymbol{m}\in R(k_i)$ for $i=0,\dots,|\mathcal{S}|-1$ are the all-zeros vectors and that $\boldsymbol{m}\in R(k_{|\mathcal{S}|})$ is such that $m_1 = n$ and $m_2 = \dots = m_{k_{|\mathcal{S}|}} = 0$. This results in the following simplification for the numerator:
\begin{equation}
\begin{array}{rl}
  & n!(1-\rho)^{n}\sum\limits_{\boldsymbol{S}\in \mathcal{S}_{|\mathcal{S}|}} \prod\limits_{j=1}^{|\mathcal{S}|} \frac{N\lambda p_{S_j}}{\mu(S_1,\dots,S_j)}
 \prod\limits_{j=l}^{|\mathcal{S}|-1} \left( 1-\frac{N\lambda \sum_{i=1}^j p_{S_i}}{\mu(S_1,\dots,S_j)}\right)^{-1}
 \left(1-\rho\right)^{-n}\prod\limits_{j=1}^{n} \left(\frac{ \lambda p_{S_l}}{j!\mu}\right)^{m_j} +o(1-\rho)\\
 
 = & n!\left(\frac{\lambda p_T}{\mu}\right)^{n}\sum\limits_{\boldsymbol{S}\in \mathcal{S}_{|\mathcal{S}|}} \prod\limits_{j=1}^{|\mathcal{S}|} \frac{N\lambda p_{S_j}}{\mu(S_1,\dots,S_j)}
 \prod\limits_{i=l}^{|\mathcal{S}|-1} \left( 1-\frac{N\lambda \sum_{i=1}^j p_{S_i}}{\mu(S_1,\dots,S_j)}\right)^{-1}
  +o(1-\rho).\\
\end{array}
\end{equation}
Combined this gives
\begin{equation}
\begin{array}{rcl}
\lim\limits_{\lambda\uparrow\mu}\mathbb{E}\left[ \left((1-\rho)Q_T\right)^n \right] &= & \lim\limits_{\lambda\uparrow\mu}
\dfrac{n!p_T^{n}\rho^n \sum\limits_{\boldsymbol{S}\in \mathcal{S}_{|\mathcal{S}|}} \prod\limits_{j=1}^{|\mathcal{S}|} \frac{N\lambda p_{S_j}}{\mu(S_1,\dots,S_j)}
 \prod\limits_{i=l}^{|\mathcal{S}|-1} \left( 1-\frac{N\lambda \sum_{i=1}^j p_{S_i}}{\mu(S_1,\dots,S_j)}\right)^{-1}
  +o(1-\rho)}{\sum\limits_{\boldsymbol{S}\in \mathcal{S}_{|\mathcal{S}|}} \prod\limits_{j=1}^{|\mathcal{S}|} \frac{N\lambda p_{S_j}}{\mu(S_1,\dots,S_j)}
 \prod\limits_{i=l}^{|\mathcal{S}|-1} \left( 1-\frac{N\lambda \sum_{i=1}^j p_{S_i}}{\mu(S_1,\dots,S_j)}\right)^{-1}
  +o(1-\rho)}\\
  
  &=& n!p_T^{n}.
\end{array}
\end{equation}
\textbf{C.o.s.\ mechanism:} We apply the steps as described in the proof of Theorem~\ref{th:conv_moments_tot} for the c.o.s.\ policy in terms of $Q_T$ rather than $Q$ to obtain the desired result.
This concludes the proof.
\end{proof}

\section{Non-critical job types} \label{sec:non_critial_jobtypes}
Subsection~\ref{subsec:general} allows for scenarios where the system does not satisfy the local stability conditions~\eqref{stabcond3} and~\eqref{stabcond4}.
While Theorem~\ref{statespacecollapse_general} focuses on the queue lengths of the job types for which the aggregate arrival rate reaches the aggregate service rate, we will now turn our attention to the remaining job types in $\mathcal{S}\setminus \mathcal{T}^*$.\\

For instance, consider the system as depicted in Figure~\ref{fig:example}, which is also often referred to as the M-model due to the shape of the compatibility graph. Pick an initial arrival rate vector $3\lambda(p_{\{1,2\}},p_{\{2,3\}})$ which could be indicated by the 5-pointed star in capacity region in Figure~\ref{fig:example4}. When $\lambda$ increases to $\lambda^* = (\mu_2+\mu_3)/(3p_{\{2,3\}})$ the boundary of the capacity region will be reached for the constraint $3\lambda p_{\{2,3\}}<\mu_2+\mu_3$ (the stability constraint imposed by job type ${\{2,3\}}$), but not for the constraint $3\lambda p_{\{1,2\}}<\mu_1+\mu_2$ (the stability constraint imposed by job type ${\{1,2\}}$), or constraint $3\lambda<3\mu = \mu_1+\mu_2+\mu_3$ (the constraint imposed by both job types).
 In fact, servers~2 and~3 will be mainly processing type-${\{2,3\}}$ jobs when $\lambda\approx\lambda^*$. In order to maintain stability, this implies that almost all type-${\{1,2\}}$ jobs must be served at server~1 alone.
 Hence $3\lambda^* p_{\{1,2\}}<\mu_1<\mu_1+\mu_2$ must hold, indeed $3\lambda^*p_{\{1,2\}} = 3\lambda^*(1-p_{\{2,3\}}) = 3\lambda^* -\mu_2-\mu_3 < 3\mu-\mu_2-\mu_3 = \mu_1$. From Theorem~\ref{statespacecollapse_general} it can be concluded that the number of type-${\{1,2\}}$ jobs becomes negligible as $\lambda\uparrow\lambda^*$ after scaling with $(1-\lambda/\lambda^*)$, and one could wonder how this quantity behaves without scaling factor. As there is a part of the system that is not critically loaded when $\lambda\uparrow\lambda^*$, this part behaves as an isolated system without job types and servers that do experience critical load. 
 
We will formalize these observations below by first defining what we mean by a \textit{truncated system}.\\

The \textit{truncated system} of a full system consists of all job types in $\mathcal{S}\setminus\mathcal{T}^*$ and all servers that are not compatible with any of the job types in $\mathcal{T}^*$, with slight abuse of notation, denoted by $\{1,\dots, N \} \setminus \mathcal{T}^*$. The arrival rate of a type-$S$ job is still given by $N\lambda p_S$, for $S\in\mathcal{S}\setminus\mathcal{T}^*$. The compatible servers of this arriving job in the truncated system are $S\setminus\mathcal{T}^*$, i.e., the compatible servers in the full system that do not experience a critical load (due to other job types). Note that in the truncated system, the probability that an arriving job is of type $S$ is now given by $p_S/\sum_{T\in\mathcal{S}\setminus\mathcal{T}^*} p_T = p_S/p_{\mathcal{S}\setminus\mathcal{T}^*}$.

\begin{lemma}\label{lem:truncated_stab}
The truncated system operating at load $\lambda^*$ is stable for both the c.o.c.\ and c.o.s.\ policies. 
\end{lemma}
\begin{proof}
From~\eqref{stabcond1} and~\eqref{stabcond2}, we know that stability of the truncated system is guaranteed if for all non-empty $\mathcal{T}\subseteq\mathcal{S}\setminus\mathcal{T}^*$
\begin{equation}
N\lambda^*p_{\mathcal{T}} = N\lambda^*\sum\limits_{S\in \mathcal{T}}p_S < \sum\limits_{n\in\left(\cup_{S\in\mathcal{T}}S\right)\setminus \mathcal{T}^*} \mu_n.
\end{equation}
By definition of the critical subset $\mathcal{T}^*$, we have $N\lambda^*p_{\mathcal{T}^*}  = \mu_{\mathcal{T}^*}$. Moreover, relying on the fact that sets $\mathcal{T}$ and $\mathcal{T}^*$ share no common job types, we have
\begin{equation}
N\lambda^*(p_{\mathcal{T}}+p_{\mathcal{T}^*}) = N\lambda^* \sum\limits_{S\in \mathcal{T}\cup\mathcal{T}^*} p_S < \sum\limits_{n\in\left(\cup_{S\in\mathcal{T}\cup\mathcal{T}^*}S\right)} \mu_n.
\end{equation}
Rewriting the above inequality results in 
\begin{equation}
N\lambda^*p_{\mathcal{T}} < \sum\limits_{n\in\left(\cup_{S\in\mathcal{T}\cup\mathcal{T}^*}S\right)} \mu_n -  \mu_{\mathcal{T}^*} = \sum\limits_{n\in\left(\cup_{S\in\mathcal{T}}S\right)\setminus \mathcal{T}^*} \mu_n.
\end{equation}
This concludes the proof.
\end{proof}

The above stability result allows us to define the random vectors $(Q^*_S)_{S\notin\mathcal{T}^*}$ and $(\tilde{Q}^*_S)_{S\notin\mathcal{T}^*}$ with the same stationary distribution as the number of jobs of each type under the c.o.c.\ policy and the number of waiting jobs of each type under the c.o.s.\ policy, respectively, in the truncated system.
As stability is guaranteed, it is meaningful to compare the queue length distribution of job types in $\mathcal{S}\setminus\mathcal{T}^*$ in the full system with the queue length distribution of the job types in the truncated system operating at load $\lambda^*$.

\begin{theorem}\label{th:ssc_noncritical }
For the c.o.c.\ policy 
\begin{equation}
(Q_S)_{S\notin \mathcal{T}^*} \xrightarrow{d} (Q_S^*)_{S\notin \mathcal{T}^*}
\end{equation} as $\lambda\uparrow\lambda^*$, and for the c.o.s.\ policy 
\begin{equation}
(\tilde{Q}_S)_{S\notin \mathcal{T}^*} \xrightarrow{d} (\tilde{Q}_S^*)_{S\notin \mathcal{T}^*}
\end{equation}
as $\lambda\uparrow\lambda^*$. Here $Q_S^*$ and $\tilde{Q}_S^*$ denote the number of type-$S$ jobs and the number of waiting type-$S$ jobs, respectively, in a truncated system as described above.
\end{theorem}

For instance, for the c.o.s.\ policy, which is essentially the JSW policy, this is a very intuitive result. It is evident that for an arriving job of type $S\in\mathcal{S}\setminus\mathcal{T}^*$ in the full system, the compatible server with the smallest workload will not be compatible with any of the jobs types in $\mathcal{T}^*$, as the queue lengths of the latter servers grow unbounded.

\begin{proof}
\textbf{C.o.c.\ mechanism:}
In order to show the weak convergence result in Theorem~\ref{th:ssc_noncritical } for the c.o.c.\ policy, it is sufficient to show convergence of the generating function of $(Q_S)_{S\notin \mathcal{T}^*}$ to the generating function $(Q_S^*)_{S\notin \mathcal{T}^*}$ by Feller's convergence theorem~\cite{Feller1971}. From Proposition~\ref{prop:pgf}, and the stability guarantees in Lemma~\ref{lem:truncated_stab}, the MGF of $(Q_S^*)_{S\notin \mathcal{T}^*}$ is given by 
\begin{equation}
\mathbb{E}\left[\prod\limits_{S\notin \mathcal{T}^*} z_S^{Q^*_S}\right] = \frac{f^*(\boldsymbol{z})}{f^*(\boldsymbol{1})},
\end{equation}
with $|z_S|\le 1$ and
\begin{equation}\label{eq:truncated_f}
f^*(\boldsymbol{z}) = 1+\sum\limits_{m=1}^{|\mathcal{S}\setminus \mathcal{T}^*|} \sum\limits_{\boldsymbol{S}\in (\mathcal{S}\setminus \mathcal{T}^*)_m} \prod\limits_{j=1}^m \frac{N\lambda^*p_{S_j}z_{S_j}}{\mu^*(S_1,\dots,S_j)}\left(1-\frac{N\lambda^*\sum_{i=1}^jp_{S_i}z_{S_i}}{\mu^*(S_1,\dots,S_j)}\right)^{-1}.
\end{equation}
Let $(\mathcal{S}\setminus \mathcal{T}^*)_m$ denote the set of ordered vectors consisting of $m$ distinct job types, not belonging to the critical set $\mathcal{T}^*$. Define
\begin{equation}
\mu^*(S_1,\dots,S_j) = \sum\limits_{n\in\left(\cup_{i=1}^j S_i\right)\setminus\mathcal{T}^*} \mu_n,
\end{equation}
i.e., the cumulative service rate of the compatible servers in the truncated system. Let $z_S \equiv 1$ if $S\in\mathcal{T}^*$, then the MGF of $(Q_S)_{S\notin \mathcal{T}^*}$ is given by
\begin{equation}
\mathbb{E}\left[\prod\limits_{S\notin \mathcal{T}^*} z_S^{Q_S}\right] =\mathbb{E}\left[\prod\limits_{S\in \mathcal{S}} z_S^{Q_S}\right] = \frac{f(\boldsymbol{z})}{f(\boldsymbol{1})},
\end{equation}
with $f$ stated in Proposition~\ref{prop:pgf}. We will multiply both the numerator and denominator with $(1-N\lambda p_{\mathcal{T}^*}/\mu _{\mathcal{T}^*})$, which is by definition strictly positive as long as $\lambda<\lambda^*$. We aim to investigate the following limit
\begin{equation}
\lim_{\lambda\uparrow\lambda^*} \frac{f(\boldsymbol{z})}{f(\boldsymbol{1})} = \lim_{\lambda\uparrow\lambda^*} \frac{\left(1-\frac{N\lambda p_{\mathcal{T}^*}}{\mu_\mathcal{T}^*}\right)f(\boldsymbol{z})}{\left(1-\frac{N\lambda p_{\mathcal{T}^*}}{\mu_\mathcal{T}^*}\right)f(\boldsymbol{1})}.
\end{equation}
Now we focus on the different components. For any $|z_S|\le 1$, we observe that 
\begin{equation}
\left(1-\frac{N\lambda p_{\mathcal{T}^*}}{\mu_\mathcal{T}^*}\right)f(\boldsymbol{z}) = \left(1-\frac{N\lambda p_{\mathcal{T}^*}}{\mu_\mathcal{T}^*}\right) + \sum\limits_{m=1}^{|\mathcal{S}|} \sum\limits_{\boldsymbol{S}\in \mathcal{S}_m}\prod\limits_{j=1}^m \frac{N\lambda p_{S_j}z_{S_j}}{\mu(S_1,\dots,S_j)}
\left(1-\frac{N\lambda \sum_{i=1}^jp_{S_i}z_{S_i}}{\mu^*(S_1,\dots,S_j)}\right)^{-1} \left(1-\frac{N\lambda p_{\mathcal{T}^*}}{\mu_\mathcal{T}^*}\right).
\end{equation}
By definition $\lambda^* = \mu_{\mathcal{T}^*}/(N p_{\mathcal{T}^*})$, so
\begin{equation}
\lim_{\lambda\uparrow\lambda^*} \left(1-\frac{N\lambda p_{\mathcal{T}^*}}{\mu_{\mathcal{T}^*}}\right) = 0.
\end{equation}
Moreover, for any $m$ and $\boldsymbol{S}$, 
\begin{equation}\label{eq:intermediate_limit}
\lim_{\lambda\uparrow\lambda^*}  \prod\limits_{j=1}^m \frac{N\lambda p_{S_j}z_{S_j}}{\mu(S_1,\dots,S_j)} = \prod\limits_{j=1}^m \frac{\mu_{\mathcal{T}^*}p_{S_j}z_{S_j}}{p_{\mathcal{T}^*}\mu(S_1,\dots,S_j)} \in (0,\infty).
\end{equation}
Computing 
\begin{equation}\label{eq:limit_1}
\lim_{\lambda\uparrow\lambda^*} \left(1-\frac{N\lambda p_{\mathcal{T}^*}}{\mu_\mathcal{T}^*}\right)  \prod\limits_{j=1}^m\left(1-\frac{N\lambda\sum_{i=1}^jp_{S_i}z_{S_i}}{\mu^*(S_1,\dots,S_j)}\right)^{-1}
\end{equation}
requires a little more attention. On the one hand, when $\{S_1,\dots,S_j\}  \neq \mathcal{T}^*$, we have
\begin{equation}
\lim_{\lambda\uparrow\lambda^*} \left(1-\frac{N\lambda\sum_{i=1}^jp_{S_i}z_{S_i}}{\mu^*(S_1,\dots,S_j)}\right)^{-1} = \left(1-\frac{\mu_{\mathcal{T}^*}\sum_{i=1}^jp_{S_i}z_{S_i}}{p_{\mathcal{T}^*}\mu^*(S_1,\dots,S_j)}\right)^{-1} \in (0,\infty).
\end{equation}
On the other hand, when $\{S_1,\dots,S_{|\mathcal{T}^*|}\}   =\mathcal{T}^*$, we observe that 
\begin{equation}
\left(1-\frac{N\lambda p_{\mathcal{T}^*}}{\mu_\mathcal{T}^*}\right) \left(1-\frac{N\lambda\sum_{i=1}^{{|\mathcal{T}^*|}}p_{S_i}z_{S_i}}{\mu^*(S_1,\dots,S_{|\mathcal{T}^*|})}\right)^{-1} = 1
\end{equation}
as $z_S\equiv 1$. 

Hence, the limit in~\eqref{eq:limit_1} is only non-zero when the first $|\mathcal{T}^*|$ elements of the ordered vector of job types $\boldsymbol{S}$ are precisely the job types in the critical set $\mathcal{T}^*$. For all other vectors~$\boldsymbol{S}$, the limit in~\eqref{eq:limit_1} is zero due to the vanishing prefactor. It allow us to focus on the ordered vectors~$\boldsymbol{S}$ which are of the form $[\boldsymbol{T},S_1,\dots,S_j ]$. The vector~$\boldsymbol{T}$ is an element of the set of all ordered permutations of the elements in the critical set~$\mathcal{T}^*$, denoted by $\mathcal{P}(\mathcal{T}^*)$. The job types $S_1,\dots,S_j $ are distinct job types belonging to $\mathcal{S}\setminus\mathcal{T}^*$. So for any $|z_S|\le 1$,
\begin{equation}\label{eq:derivation_limit}
\begin{array}{rcl}
\lim_{\lambda\uparrow\lambda^*} \left(1-\frac{N\lambda p_{\mathcal{T}^*}}{\mu_\mathcal{T}^*}\right)f(\boldsymbol{z})
&=&
\sum\limits_{m=0}^{|\mathcal{S}\setminus\mathcal{T}^*|} \sum\limits_{\boldsymbol{T}\in \mathcal{P}(\mathcal{T}^*)} \sum\limits_{\boldsymbol{S} \in (\mathcal{S}\setminus\mathcal{T}^*)_m} \prod\limits_{j=1}^{|\mathcal{T}^*|}\frac{N\lambda^*p_{T_j}}{\mu(T_1,\dots,T_j)}
\prod\limits_{j=1}^{m}\frac{N\lambda^*p_{S_j}z_{S_j}}{\mu(\boldsymbol{T},S_1,\dots,S_j)}\\

& &\times \prod\limits_{j=1}^{|\mathcal{T}^*|-1}\left(1-\frac{N\lambda^*\sum_{i=1}^jp_{T_i}}{\mu(T_1,\dots,T_j)}\right)^{-1}
\prod\limits_{j=1}^{m}\left(1-\frac{N\lambda^*p_{\mathcal{T}^*}+N\lambda^*\sum_{i=1}^j p_{S_i}z_{S_i}}{\mu(\boldsymbol{T},S_1,\dots,S_j)}\right)^{-1}\\

& =& \left\{ \sum\limits_{\boldsymbol{T}\in \mathcal{P}(\mathcal{T}^*)}  \prod\limits_{j=1}^{|\mathcal{T}^*|}\frac{N\lambda^*p_{T_j}}{\mu(T_1,\dots,T_j)} \prod\limits_{j=1}^{|\mathcal{T}^*|-1}\left(1-\frac{N\lambda^*\sum_{i=1}^jp_{T_i}}{\mu(T_1,\dots,T_j)}\right)^{-1} \right\}  \\
& &\times \left\{1+ \sum\limits_{m=1}^{|\mathcal{S}\setminus\mathcal{T}^*|}  \sum\limits_{\boldsymbol{S} \in (\mathcal{S}\setminus\mathcal{T}^*)_m} \prod\limits_{j=1}^{m}\frac{N\lambda^*p_{S_j}z_{S_j}}{\mu(\mathcal{T}^*,S_1,\dots,S_j)}\left(1-\frac{\mu_{\mathcal{T}^*}+N\lambda^*\sum_{i=1}^j p_{S_i}z_{S_i}}{\mu(\mathcal{T}^*,S_1,\dots,S_j)}\right)^{-1}\right\}\\

& =& \left\{ \sum\limits_{\boldsymbol{T}\in \mathcal{P}(\mathcal{T}^*)}  \prod\limits_{j=1}^{|\mathcal{T}^*|}\frac{N\lambda^*p_{T_j}}{\mu(T_1,\dots,T_j)} \prod\limits_{j=1}^{|\mathcal{T}^*|-1}\left(1-\frac{N\lambda^*\sum_{i=1}^jp_{T_i}}{\mu(T_1,\dots,T_j)}\right)^{-1} \right\}  \\
& & \times \left\{1+ \sum\limits_{m=1}^{|\mathcal{S}\setminus\mathcal{T}^*|}  \sum\limits_{\boldsymbol{S} \in (\mathcal{S}\setminus\mathcal{T}^*)_m} \prod\limits_{j=1}^{m}\frac{N\lambda^*p_{S_j}z_{S_j}}{\mu(S_1,\dots,S_j)-\mu_{\mathcal{T}^*}}\left(1-\frac{N\lambda^*\sum_{i=1}^j p_{S_i}z_{S_i}}{\mu(\mathcal{T}^*,S_1,\dots,S_j)-\mu_{\mathcal{T}^*}}\right)^{-1}\right\}\\

& =& \left\{ \sum\limits_{\boldsymbol{T}\in \mathcal{P}(\mathcal{T}^*)}  \prod\limits_{j=1}^{|\mathcal{T}^*|}\frac{N\lambda^*p_{T_j}}{\mu(T_1,\dots,T_j)} \prod\limits_{j=1}^{|\mathcal{T}^*|-1}\left(1-\frac{N\lambda^*\sum_{i=1}^jp_{T_i}}{\mu(T_1,\dots,T_j)}\right)^{-1} \right\}f^*(\boldsymbol{z}),\\
\end{array}
\end{equation}
the limit is finite.
The prefactor
\begin{equation}\label{eq:A}
A \coloneqq  \sum\limits_{\boldsymbol{T}\in \mathcal{P}(\mathcal{T}^*)}  \prod\limits_{j=1}^{|\mathcal{T}^*|}\frac{N\lambda^*p_{T_j}}{\mu(T_1,\dots,T_j)} \prod\limits_{j=1}^{|\mathcal{T}^*|-1}\left(1-\frac{N\lambda^*\sum_{i=1}^jp_{T_i}}{\mu(T_1,\dots,T_j)}\right)^{-1}
\end{equation} of the above expression does not depend on $\boldsymbol{z}$.
We can conclude that 
\begin{equation}
\lim_{\lambda\uparrow\lambda^*} \frac{f(\boldsymbol{z})}{f(\boldsymbol{1})} =  \frac{f^*(\boldsymbol{z})}{f^*(\boldsymbol{1})}.
\end{equation}
This concludes the proof for the c.o.c.\ policy.\\

\textbf{C.o.s.\ mechanism:} The proof of the convergence result for the c.o.s.\ policy mimics the proof of that for the c.o.c.\ policy. Namely, the aim is to show convergence of the generating function of $(\tilde{Q}_S)_{S\notin \mathcal{T}^*}$ to 
\begin{equation}
\mathbb{E}\left[\prod\limits_{S\notin\mathcal{T}^*} z_S^{\tilde{Q}^*_S} \right]= \frac{g^*(\boldsymbol{z})}{g^*(\boldsymbol{1})},
\end{equation}
with $|z_S|\le 1$ and
\begin{equation}\label{eq:truncated_g}
g^*(\boldsymbol{z}) = \sum\limits_{L=0}^{N^*}\sum\limits_{\boldsymbol{u}\in \mathcal{N}^*_L}\prod\limits_{l=1}^L \frac{\mu_{u_l}}{\lambda_{\mathcal{C}(u_1,\dots,u_l)}} 
 \sum\limits_{m=0}^{|\mathcal{S}\setminus\mathcal{T}^*|}\sum\limits_{\boldsymbol{S}\in (\mathcal{S}\setminus\mathcal{T}^*)_m^{\boldsymbol{u}}} \prod\limits_{j=1}^m \frac{N\lambda^* p_{S_j}z_{S_j}}{\mu^*(S_1,\dots,S_j)}\left(1{-}\frac{N\lambda^*}{\mu^*(S_1,\dots,S_j)}\sum\limits_{i=1}^j p_{S_i}z_{S_i}\right)^{-1}.
\end{equation}
In the above expression, based on the results in Lemma~\ref{lem:truncated_stab} and Proposition~\ref{prop:pgfcos}, $N^* = |\{1,\dots,N\}\setminus\mathcal{T}^*|$, and $\mathcal{N}^*_L$ denotes the set of ordered vectors of length $L$ presenting the idle servers in the truncated system. Again we are interested in computing the limit for $\lambda\uparrow\lambda^*$ of 
\begin{equation}
\mathbb{E}\left[\prod\limits_{S\notin \mathcal{T}^*} z_S^{\tilde{Q}_S}\right] =\mathbb{E}\left[\prod\limits_{S\in \mathcal{S}} z_S^{\tilde{Q}_S}\right] = \frac{g(\boldsymbol{z})}{g(\boldsymbol{1})} = \frac{\left(1-\frac{N\lambda p_{\mathcal{T}^*}}{\mu^{\mathcal{T}^*}}\right)g(\boldsymbol{z})}{\left(1-\frac{N\lambda p_{\mathcal{T}^*}}{\mu^{\mathcal{T}^*}}\right)g(\boldsymbol{1})} ,
\end{equation}
with $z_S\equiv 1$ if $S\in\mathcal{T}^*$ and $g$ as given in Proposition~\ref{prop:pgfcos}. The above observations for the c.o.c.\ policy, together with the observation in~\eqref{eq:lim_component_cos}, will result in
\begin{equation}
\lim_{\lambda\uparrow\lambda^*} \left(1-\frac{N\lambda p_{\mathcal{T}^*}}{\mu_{\mathcal{T}^*}}\right)g(\boldsymbol{z}) = A g^*(\boldsymbol{z})
\end{equation}
for any $|z_S|\le 1$, $S\notin\mathcal{T}^*$. The prefactor $A$ is defined in~\eqref{eq:A}. This concludes the proof of the result for the c.o.s.\ policy.
\end{proof}

We can now use the above proof method to combine the results in Theorems~\ref{statespacecollapse_general} and~\ref{th:ssc_noncritical }. More precisely, as formalized in Theorem~\ref{th:ssc_nonscaled}, when $\lambda\uparrow\lambda^*$ the full system will behave like two independent subsystems and the (scaled) queue length of the critically loaded subsystem will exhibit state space collapse, while the non-critically loaded subsystem will behave as a truncated system. Hence, the complexity of the analysis of the full system operating at critical load~$\lambda^*$ is significantly reduced. Next to the state space collapse, the truncated system can be much smaller in size compared to the full system in which case it might be feasible to manipulate the detailed product-form expressions in Subsection~\ref{sec:prodform} to obtain the stationary distributions of $(Q^*_S)_{S\notin\mathcal{T}^*}$ and $(\tilde{Q}^*_S)_{S\notin\mathcal{T}^*}$. For instance, 
\begin{equation}
\mathbb{P}(Q^* = q) = \sum\limits_{\boldsymbol{c}\colon |\boldsymbol{c}| = q} \pi_{\text{c.o.c.}} (\boldsymbol{c})~~\text{and}~~\mathbb{P}(\tilde{Q}^* = q) = \sum\limits_{\boldsymbol{c},\boldsymbol{u}\colon |\boldsymbol{c}| = q} \pi_{\text{c.o.s.}} (\boldsymbol{c},\boldsymbol{u}),
\end{equation}
for $q\ge 0$ with $Q^* = \sum_{S\notin\mathcal{T}^*} Q^*_S$ and $\tilde{Q}^* = \sum_{S\notin\mathcal{T}^*} \tilde{Q}^*_S$ the total number of (waiting) jobs.



\begin{proof}[Proof of Theorem~\ref{th:ssc_nonscaled}.]
\textbf{C.o.c.\ mechanism:}
In order to show the above-mentioned result, it is sufficient to show convergence of the corresponding generating functions. More precisely, for any $z_S\coloneqq \exp\left(-\left(1-\frac{\lambda}{\lambda^*}\right)t_S\right)$ for $S\in\mathcal{T}^*$, $z_S\coloneqq \exp\left(-t_S\right)$ for $S\notin\mathcal{T}^*$ and $t_S\ge 0$,
\begin{equation}\label{eq:non_scaled_aim}
\mathbb{E}\left[\prod\limits_{S\in\mathcal{S}}z_S^{Q_S}\right] \rightarrow \left(1+\sum\limits_{S\in\mathcal{T}^*} \frac{p_S}{p_{\mathcal{T}^*}}t_S\right)^{-1}\cdot \frac{f^*(\boldsymbol{z})}{f^*(\boldsymbol{1})},
\end{equation}
as $\lambda\uparrow\lambda^*$. The first factor on the right-hand side can be recognized as the MGF of the random variable $\mathrm{Exp}(1)\left(\frac{p_S}{p_{\mathcal{T}^*}}\right)_{S\in \mathcal{T}^*}$, while the second factor is the MGF of $(Q_S^*)_{S\notin \mathcal{T}^*}$, with $f^*$ as defined in \eqref{eq:truncated_f}. Observe that 
\begin{equation}\label{eq:nonscaled_toshow}
\mathbb{E}\left[\prod\limits_{S\in\mathcal{S}}z_S^{Q_S}\right] =  \frac{f(\boldsymbol{z})}{f(\boldsymbol{1})} = \frac{1-\frac{N\lambda p_{\mathcal{T}^*}}{\mu_{\mathcal{T}^*}}}{1-\frac{N\lambda }{\mu_{\mathcal{T}^*}}\sum\limits_{S\in\mathcal{T}^*}p_Sz_S} 
\frac{f(\boldsymbol{z})\left(1-\frac{N\lambda}{\mu_{\mathcal{T}^*}}\sum\limits_{S\in\mathcal{T}^*}p_Sz_S\right)}{f(\boldsymbol{1})\left(1-\frac{N\lambda p_{\mathcal{T}^*}}{\mu_{\mathcal{T}^*}}\right)}, 
\end{equation}
with $f$ defined in Proposition~\ref{prop:pgf}.
Using l'H\^opital's rule, the first factor indeed converges to the first factor of~\eqref{eq:non_scaled_aim} as $\lambda\uparrow\lambda^*$. The computation of the limit of the second factor occurs along the same lines as the proof of Theorem~\ref{th:ssc_noncritical }. We observe that
\begin{equation}
\lim_{\lambda\uparrow\lambda^*} \left(1-\frac{N\lambda}{\mu_{\mathcal{T}^*}}\sum\limits_{S\in\mathcal{T}^*}p_{S}z_S\right) = 0,
\end{equation}
as $\lambda^* = \mu_{\mathcal{T}^*}/(N p_{\mathcal{T}^*})$.
Moreover, for any $m$ and $\boldsymbol{S}$, we know that \eqref{eq:intermediate_limit} holds. Then we investigate 
\begin{equation}\label{eq:limit_1b}
\lim_{\lambda\uparrow\lambda^*} \left(1-\frac{N\lambda}{\mu_\mathcal{T}^*}\sum\limits_{S\in\mathcal{T}^*}p_{S}z_S\right)  \prod\limits_{j=1}^m\left(1-\frac{N\lambda\sum_{i=1}^jp_{S_i}z_{S_i}}{\mu^*(S_1,\dots,S_j)}\right)^{-1}.
\end{equation}
 On the one hand, when $\{S_1,\dots,S_j\}  \neq \mathcal{T}^*$, we have
\begin{equation}
\lim_{\lambda\uparrow\lambda^*} \left(1-\frac{N\lambda\sum_{i=1}^jp_{S_i}z_{S_i}}{\mu^*(S_1,\dots,S_j)}\right)^{-1} = \left(1-\frac{\mu_{\mathcal{T}^*}\sum_{i=1}^jp_{S_i}z_{S_i}}{p_{\mathcal{T}^*}\mu^*(S_1,\dots,S_j)}\right)^{-1} \in (0,\infty).
\end{equation}
On the other hand, when $\{S_1,\dots,S_{|\mathcal{T}^*|}\}   =\mathcal{T}^*$, we observe that 
\begin{equation}
\left(1-\frac{N\lambda}{\mu_\mathcal{T}^*}\sum\limits_{S\in\mathcal{T}^*}p_{S}z_S\right) \left(1-\frac{N\lambda\sum_{i=1}^{{|\mathcal{T}^*|}}p_{S_i}z_{S_i}}{\mu^*(S_1,\dots,S_{|\mathcal{T}^*|})}\right)^{-1} = 1.
\end{equation}
With these observations, a similar derivation as in~\eqref{eq:derivation_limit} shows that 
\begin{equation}
\lim_{\lambda\uparrow\lambda^*} \left(1-\frac{N\lambda}{\mu_\mathcal{T}^*}\sum\limits_{S\in\mathcal{T}^*}p_{S}z_S\right)f(\boldsymbol{z})
= A(\boldsymbol{z})f^*(\boldsymbol{z}).
\end{equation}
The prefactor
\begin{equation}\label{eq:A_z}
A(\boldsymbol{z}) \coloneqq  \sum\limits_{\boldsymbol{T}\in \mathcal{P}(\mathcal{T}^*)}  \prod\limits_{j=1}^{|\mathcal{T}^*|}\frac{N\lambda^*p_{T_j}z_{T_j}}{\mu(T_1,\dots,T_j)} \prod\limits_{j=1}^{|\mathcal{T}^*|-1}\left(1-\frac{N\lambda^*\sum_{i=1}^jp_{T_i}z_{T_i}}{\mu(T_1,\dots,T_j)}\right)^{-1},
\end{equation} and
$\lim_{\lambda\uparrow\lambda^*} A(\boldsymbol{z}) = A(\boldsymbol{1})$. Hence, the second factor on the right-hand side in~\eqref{eq:nonscaled_toshow} indeed converges to $\frac{f^*(\boldsymbol{z})}{f^*(\boldsymbol{1})}$ as $\lambda\uparrow\lambda^*$. 
This concludes the proof for the c.o.c.\ policy.\\

\textbf{C.o.s.\ mechanism:} In order to show the above-mentioned result, it is again sufficient to show convergence of the corresponding generating function. More precisely, for any $z_S\coloneqq \exp\left(-\left(1-\frac{\lambda}{\lambda^*}\right)t_S\right)$ for $S\in\mathcal{T}^*$, $z_S\coloneqq \exp\left(-t_S\right)$ for $S\notin\mathcal{T}^*$ and $t_S\ge 0$,
\begin{equation}\label{eq:non_scaled_aim_cos}
\mathbb{E}\left[\prod\limits_{S\in\mathcal{S}}z_S^{\tilde{Q}_S}\right] \rightarrow \left(1+\sum\limits_{S\in\mathcal{T}^*} \frac{p_S}{p_{\mathcal{T}^*}}t_S\right)^{-1}\cdot \frac{g^*(\boldsymbol{z})}{g^*(\boldsymbol{1})},
\end{equation}
as $\lambda\uparrow\lambda^*$. The first factor on the right-hand side can be recognized as the MGF of the random variable $\mathrm{Exp}(1)\left(\frac{p_S}{p_{\mathcal{T}^*}}\right)_{S\in \mathcal{T}^*}$. The second factor is the MGF of $(\tilde{Q}_S^*)_{S\notin \mathcal{T}^*}$, with $g^*$ as defined in \eqref{eq:truncated_g}. Observe that 
\begin{equation}
\mathbb{E}\left[\prod\limits_{S\in\mathcal{S}}z_S^{\tilde{Q}_S}\right] =  \frac{g(\boldsymbol{z})}{g(\boldsymbol{1})} = \frac{1-\frac{N\lambda p_{\mathcal{T}^*}}{\mu_{\mathcal{T}^*}}}{1-\frac{N\lambda }{\mu_{\mathcal{T}^*}}\sum\limits_{S\in\mathcal{T}^*}p_Sz_S} 
\frac{g(\boldsymbol{z})\left(1-\frac{N\lambda}{\mu_{\mathcal{T}^*}}\sum\limits_{S\in\mathcal{T}^*}p_Sz_S\right)}{g(\boldsymbol{1})\left(1-\frac{N\lambda p_{\mathcal{T}^*}}{\mu_{\mathcal{T}^*}}\right)}, 
\end{equation}
with $g$ defined in Proposition~\ref{prop:pgfcos}.
Using l'H\^opital's rule, the first factor indeed converges to the first factor of~\eqref{eq:non_scaled_aim_cos} as $\lambda\uparrow\lambda^*$. In order to show convergence of the second factor, one can use the same ideas as in the proof for the c.o.c.\ mechanism.
This concludes the proof.
\end{proof}

\begin{remark}\normalfont
The weak convergence result in Theorem~\ref{th:ssc_nonscaled} for the c.o.s.\ mechanism could even be extended to
\begin{equation}
\left(\left(1-\frac{\lambda}{\lambda^*}\right)(Q_S)_{S\in \mathcal{T}^*},(\tilde{Q}_S)_{S\notin \mathcal{T}^*}\right) \xrightarrow{d} \left(\mathrm{Exp}(1)\left(\frac{p_S}{p_{\mathcal{T}^*}}\right)_{S\in \mathcal{T}^*}, (\tilde{Q}_S^*)_{S\notin \mathcal{T}^*}\right),
\end{equation}
where the total number of jobs of each type in the critically loaded subsystem is considered. To prove this result we require the observation made in the proof of Theorem~\ref{statespacecollapse_general}, namely there are only a finite number of jobs in service of each of the types.
\end{remark}

\section{Simulation results for c.o.s\ mechanism}\label{app:simulations}
\subsection{Simulation set up}\label{subsec:sim_settings}
We investigate the performance of the c.o.s.\ policy for various compatibility graphs with $N=4, 10$ and 20 servers. We will mainly focus on two different classes of compatibility constraints, namely power-of-$d$ (Po$d$) and circular-$d$ constraints. The first class of constraints is outlined in Remark~\ref{remark:pod}. The latter class of assignment constraints for a fixed value of~$d$ consists of~$N$ job types, labeled as $n=1, 2,\dots, N$, each with arrival rate~$\lambda$ or $p_n = \frac{1}{N}$. A type-$n$ job is compatible with servers $n,n+1,\dots,n+d-1\mod N$.

We distinguish between the following settings:
\begin{enumerate}
\item $N=4$, Po$d$ with $d=2$. Hence, $p_S \equiv 1/\binom{N}{d} = 1/6$ for all $S\in\mathcal{S}$.
\item $N=4$, circular-$d$ with $d=2$. Hence, $p_n=1/4$ for $n=1,2,3,4$.
\item $N=4$, circular-$d$ with $d=2$ and heterogeneous arrival rates. We set  $p_n=4/10$ for $n=1,3$ and $p_n=1/10$ for $i=2,4$.
\item $N=4$, $\mathcal{S} = \left\{ \{1\},\{2\},\{3\},\{4\},\{1,2,3,4\}\right\}$. Job-type $n$ is only compatible with server $n$ for $n=1,2,3,4$ and $p_n = 7/40 = 0.175$. The fifth job type is compatible with all servers and $p_5 = 3/10$.
\item $N=10$, Po$d$ with $d=2$. Hence, $p_S \equiv 1/\binom{N}{d} = 1/45\approx 0.0222$ for all $S\in\mathcal{S}$.
\item $N=10$, circular-$d$ with $d=2$. Hence, $p_n=1/10$ for $n=1,\dots,10$.
\item $N=10$, Po$d$ with $d=8$. Hence, $p_S \equiv 1/\binom{N}{d} = 1/45\approx 0.0222$ for all $S\in\mathcal{S}$.
\item $N=10$, circular-$d$ with $d=8$. Hence, $p_n=1/10$ for $n=1,\dots,10$.
\item $N=20$, Po$d$ with $d=2$. Hence, $p_S \equiv 1/\binom{N}{d} = 1/190\approx 0.00526$ for all $S\in\mathcal{S}$.
\item $N=20$, circular-$d$ with $d=2$. Hence, $p_n=1/20= 0.05$ for $n=1,\dots,20$.
\item $N=20$, Po$d$ with $d=18$. Hence, $p_S \equiv 1/\binom{N}{d} = 1/190\approx 0.00526$ for all $S\in\mathcal{S}$.
\item $N=20$, circular-$d$ with $d=18$. Hence, $p_n=1/20 = 0.05$ for $n=1,\dots,20$.
\end{enumerate}
When the servers all operate at the same speed $\mu_n \equiv 1$ and $\lambda=0.9$, the local stability conditions in~\eqref{stabcond3} and~\eqref{stabcond4} are satisfied for all the above settings. Relying on the result in Theorem~\ref{statespacecollapse}, we will investigate how the total number of customers $Q$ in the various settings can be related to the heavy-traffic approximation $Q'\coloneqq (1-\lambda/\mu)^{-1}\text{Exp}(1)$, an exponentially distributed random variable with parameter $(1-\lambda/\mu)$. 
For each of the different settings we carry out 10 simulation runs, each run consisting of 500000 events and the initial number of jobs in the system is sampled from $(1-\lambda/\mu)^{-1}\text{Exp}(1)$.

\subsection{Simulation results}

\begin{table}[h]
\centering
\begin{tabular}{|l|c|c|c|c|}
\hline
                & setting 1        & setting 2        & setting 3        & setting 4       \\ \hline
$\mathbb{E}[Q]$ & 12.945 (0.314)   & 13.182 (0.320)   & 14.114 (0.294)   & 15.976 (0.329)   \\
$\text{Var}(Q)$ & 100.797 (12.607) & 102.386 (12.387) & 107.290 (11.596) & 110.839 (12.386) \\
$\mathbb{E}\left[\frac{Q_1}{Q}\mathds{1}\{Q>0\}\right]$ & 0.166 (0.001) & 0.249 (0.001) & 0.408 (0.002) & 0.203 (0.010) \\
$\mathbb{E}\left[\frac{Q_M}{Q}\mathds{1}\{Q>0\}\right]$  & 0.167 (0.001)    & 0.246 (0.014)    & 0.091 (0.001)    & 0.195 (0.001)  \\
\hline
                & setting 5           & setting 6            & setting 7           & setting 8  \\ \hline
$\mathbb{E}[Q]$ & 23.079 (0.370)      & 26.836 (0.436)       & 16.043 (0.331)      & 15.995 (0.333)  \\
$\text{Var}(Q)$ & 134.485 (12.474)    & 160.234 (12.616)     & 102.051 (13.263)    & 100.277 (12.901) \\
$\mathbb{E}\left[\frac{Q_1}{Q}\mathds{1}\{Q>0\}\right]$ ($\times 10^{-2}$) & 2.217 (0.0308) & 9.697 (0.0932) & 2.218 (0.0287) & 9.963 (0.0987) \\
$\mathbb{E}\left[\frac{Q_M}{Q}\mathds{1}\{Q>0\}\right]$ ($\times 10^{-2}$)   & 2.226 (0.0218) & 10.052 (0.1314) & 2.223 (0.0219) & 10.022 (0.0552)  \\
\hline
                & setting 9           & setting 10            & setting 11           & setting 12  \\ \hline
$\mathbb{E}[Q]$ & 41.060 (0.462)      & 52.330 (0.724)       & 23.759 (0.294)      & 23.759 (0.295)  \\
$\text{Var}(Q)$ & 194.779 (13.684)    & 293.639 (21.352)     & 99.234 (11.747)    & 99.245 (11.765) \\
$\mathbb{E}\left[\frac{Q_1}{Q}\mathds{1}\{Q>0\}\right]$ ($\times 10^{-2}$) & 0.531 (0.0167) & 5.053 (0.0901) & 0.529 (0.0196) & 5.002 (0.0519) \\
$\mathbb{E}\left[\frac{Q_M}{Q}\mathds{1}\{Q>0\}\right]$ ($\times 10^{-2}$)   & 0.522 (0.0187) & 5.026 (0.126) & 0.518 (0.0151) & 5.000 (0.0430)  \\
\hline
\end{tabular}
\caption{An overview of the simulation results, averaged over the various simulation runs per setting; the corresponding standard deviations are indicated between brackets.}\label{tab:sim_results}
\end{table}

\begin{figure}[t]
\centering
  \begin{tabular}{c}
  \includegraphics[scale=1]{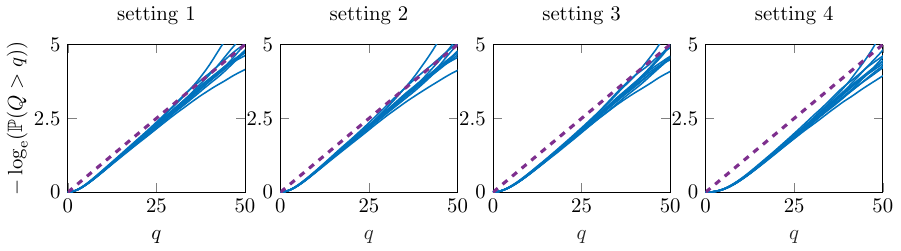}\\
  \includegraphics[scale=1]{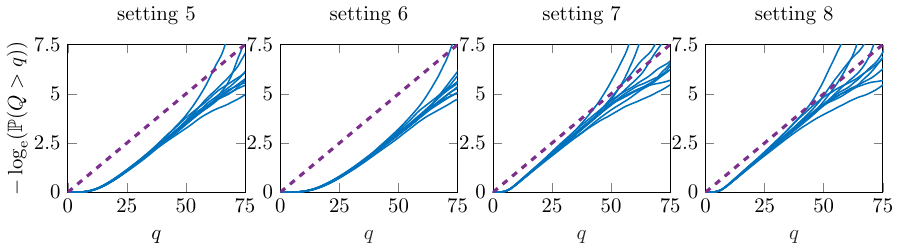} \\
  \includegraphics[scale=1]{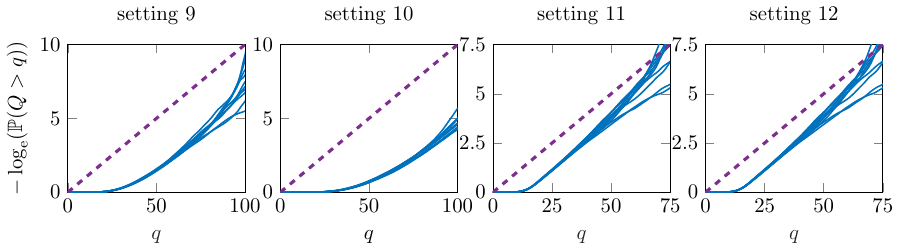}
  \end{tabular}
  \caption{A visualization of the empirical total queue length ($Q$) distribution (solid, blue lines) for the various simulation runs per setting. The dashed line represents the total queue length distribution for the heavy-traffic approximation, $Q' = (1-\lambda/\mu)^{-1}\text{Exp}(1) = 10\text{Exp}(1)$. }\label{fig:sim_results}
\end{figure}

Based on the result in Theorem~\ref{statespacecollapse}, one could expect that the distribution of $Q'=(1-\lambda/\mu)^{-1}\text{Exp}(1)$ serves as a good approximation for distribution of the total number of jobs in the system $Q$, and in particular $\mathbb{E}[Q] \approx (1-\lambda/\mu)^{-1} = 10$ and Var$(Q)\approx (1-\lambda/\mu)^{-2} = 100$. Studying the simulation results, this seems a reasonable approximation when the compatibility graph is fairly dense and symmetric. In order to demonstrate that, we summarized the simulation results in Table~\ref{tab:sim_results} and Figure~\ref{fig:sim_results}. Table~\ref{tab:sim_results} shows some performance measures, like the mean total queue length, averaged over the various simulation runs for the different settings. Figure~\ref{fig:sim_results} visualizes the  empirical cumulative distributions of the total queue length $Q$ for the various simulation runs per setting (solid blue lines), together with the cumulative distribution of the heavy-traffic approximation $Q'$ (dashed purple line).\\

Focusing on settings~1 to~4, all settings with $N=4$ servers, we clearly observe some agreement when comparing the simulations with the heavy-traffic approximation. Only around 50 jobs in the system, the theoretical approximation and the empirical distributions of the simulations are shaped differently. However, as $\mathbb{P}(Q'>50) <0.007$, it is reasonable to believe that these differences are caused by the restricted number of events per simulation, such that the tail probabilities will become less accurate. The same observation can be made for all remaining settings as well. Out of these four settings, setting~4 has a rather sparse compatibility graph and deviates more from the approximated quantities than the other three settings. In fact, this setting amounts to random assignment for $70\%$ of the jobs. So, in a way it is remarkable that this setting already shows the characteristics of a fully pooled system. Furthermore, this setting is highly asymmetrical when it comes to the set of compatible servers per arriving job, i.e, one server versus all servers.
\\

For the larger systems we see that the theoretical approximation deviates more from the obtained performance measure when the corresponding compatibility graph is less dense, i.e., comparing settings~5 and~6 with settings~7 and~8 and settings~9 and~10 with settings~11 and~12. As we set $\lambda/\mu = 9/10$, we are considering a reasonably high load, but of course it is also too low to really capture heavy-traffic behavior. 
Hence, when considering the Po$d$ or circular-$d$ setting, the accuracy of the heavy-traffic approximation increases for larger values of~$d$. Note that taking larger values of~$d$ is indeed one way to make the graph more dense, i.e, each arriving job will have access to more compatible servers. 

Another way to increase the flexibility of the system, or the density of the compatibility graph, is to increase the number of different job types. For fixed values of $d$ we notice that the heavy-traffic approximation $Q'$ is slightly more accurate for the Po$d$ setting than the circular-$d$ setting, for example by comparing settings~5 and~6. 

It has been observed before, for instance in~\cite{Weng2020}, that the compatibility graph must satisfy some density property, referred to as \textit{well connectedness}, to match the performance of a fully flexible system. It is worth emphasizing that they considered a many-server regime rather than a heavy-traffic regime. However, relying on many-server results for the JSQ($d$) policy, we can explain the shape of the empirical distributions in settings~5 and~9. For these setttings the number of compatible servers per job type, $d$, is significantly smaller than the total number of servers in the system,~$N$. As the c.o.s.\ policy amounts to the JSW~policy, and it is believed that JSW and JSQ policies tend to behave similarly, we may use the results derived by Mitzenmacher~\cite{mitzenmacher2001power}. More precisely, for the JSQ($d$) policy it is known that $\mathbb{P}(Q^{\text{JSQ($d$)}}>q) = \lambda^{\frac{d^q-1}{d-1}}$ in the many-server regime. Hence,
\[
-\log_{\mathrm e} \mathbb{P}(Q>q)= -\log_{\mathrm e} \mathbb{P}(Q^{\text{JSW($d$)}}>q) \approx -\log_{\mathrm e}(\lambda)\frac{d^q-1}{d-1},
\]
which is of course not linear in its argument $q$, but in fact has this U-shape that we also observe in the empirical distributions in settings~5 and~9. So for small, fixed values of $d$, moderate load per server, and a large number of servers~$N$, it is to be expected to observe results which do not match with our theoretical approximations. \\  



Finally, we use $\mathbb{E}\left[\frac{Q_S}{Q}\mathds{1}\{Q>0\}\right]$ to check whether the relative proportions of the various job types in the system remain unchanged compared to the arrival probabilities $p_S$, i.e., this gives an indication if the state space collapse property holds for reasonably high load. Comparing these quantities in Table~\ref{tab:sim_results} with the probabilities given in Subsection~\ref{subsec:sim_settings}, we indeed observe some agreement with the exception of setting 4. As noted above, this setting consists of a highly asymmetrical compatibility graph due to the single job type that allows for some server flexibility.

\section{Proofs of Corollaries~\ref{cor:R_i} and~\ref{cor:SojournWaiting}}\label{app:HT}
The result stated in Corollary~\ref{cor:R_i} follows from straightforward manipulations of the definitions of the MGFs of $(R_n)_{n=1,\dots,N}$, $(Q_S)_{S\in\mathcal{S}}$ and  $(\tilde{Q}_S)_{S\in\mathcal{S}}$.
\begin{proof}[Proof of Corollary~\ref{cor:R_i}.]
The stationary number of replicas at server $n$ satisfies  
\begin{equation}\label{eq:R_n}
R_n = \sum_{S:n\in S} Q_S \text{~~~and~~~}\sum_{S:n\in S} \tilde{Q}_S \le R_n \le 1+ \sum_{S:n\in S} \tilde{Q}_S
\end{equation} for the c.o.c.\ and c.o.s.\ mechanism, respectively. 
It suffices to show that 
\begin{equation}
\mathbb{E}\left[ \prod\limits_{n=1}^N \exp\left(-\left(1-\frac{\lambda}{\lambda^*}\right)t_nR_n \right)\right] \rightarrow \left(1 + \sum\limits_{n=1}^N t_nq_n\right)^{-1}
\end{equation}
for $t_n\ge 0$. For the c.o.c.\ mechanism, it can be observed that
\begin{equation}
\mathbb{E}\left[ \prod\limits_{n=1}^N \exp\left(-\left(1-\frac{\lambda}{\lambda^*}\right)t_nR_n \right)\right]  = \mathbb{E}\left[ \prod\limits_{S\in\mathcal{S}} \exp\left(-\left(1-\frac{\lambda}{\lambda^*}\right)Q_S \sum\limits_{n\in S}t_n \right)\right] 
\end{equation}
by using relation~\eqref{eq:R_n}. Due to Theorem~\ref{statespacecollapse_general}, we know that the right-hand side converges to
\[ \left(1 + \sum\limits_{S\in\mathcal{S}} p_S \sum\limits_{n\in S} t_n\right)^{-1}
\]
as $\lambda\uparrow\lambda^*$, from which the desired result follows. A similar approach, using~\eqref{eq:R_n}, \eqref{eq:convergence_waiting} and the squeeze theorem to compute the limit of interest, establishes the result for the c.o.s.\ mechanism. This concludes the proof.
\end{proof}

The distributional form of Little's law~\cite{Keilson1988} allows us to derive the result in Corollary~\ref{cor:SojournWaiting} for the sojourn and waiting time in a heavy-traffic regime.   
Some details are described below.

\begin{proof}[Proof of Corollary~\ref{cor:SojournWaiting}.]
The result will be shown for the c.o.c.\ mechanism; the proof for the c.o.s.\ mechanism is highly similar and will be omitted.
Fix a job type $S\in\mathcal{S}$. Then the PGF of $Q_S$  is equal to
\begin{equation}\label{eq:singlevarPGF}
G_{Q_S}(z) = \mathbb{E}\left[z^{Q_S}\right] = \mathbb{E}\left[z^{Q_S}\prod\limits_{\tilde{S}\in \mathcal{S}\setminus S}1^{Q_{\tilde{S}}}\right]
\end{equation}
with $|z|\le 1$.
The MGF of the sojourn time is defined as $\mathcal{L}_{V_S}(s)=\mathbb{E}\left[{\rm e}^{-sV_S}\right]$. As jobs of the same type leave the system in order of arrival, it is known from~\cite[Theorem~1]{Keilson1988} that 
\begin{equation}
G_{Q_S}(z) = \mathcal{L}_{V_S}(\lambda_S(1-z)),
\end{equation}
with $\lambda_S = N\lambda p_S$.
Let $\mathcal{T}^*, p_{\mathcal{T}^*}, \mu_{\mathcal{T}^*}, \rho_{\mathcal{T}^*}$ and $\lambda^*$ be as defined in Definition~\ref{def:critical_subset_arrrate}.
In order to deduce the limiting behavior of $(1-\frac{\lambda}{\lambda^*})V_S$, observe the following intermediate steps:
\begin{equation}
\mathcal{L}_{\left(1-\frac{\lambda}{\lambda^*}\right)V_S} (s) = \mathcal{L}_{V_S}\left(\left(1-\frac{\lambda}{\lambda^*}\right)s\right)
= G_{Q_S}\left(1-\left(1-\frac{\lambda}{\lambda^*}\right)\frac{s}{\lambda_S}\right).
\end{equation}
Using \eqref{eq:singlevarPGF} and similar arguments as in the proof of Theorem~\ref{statespacecollapse_general} while taking the limit $\lambda\uparrow\lambda^*$ results in 
\begin{equation}
\lim_{\lambda\uparrow\lambda^*} \mathcal{L}_{\left(1-\frac{\lambda}{\lambda^*}\right)V_S} (s) = 
\begin{cases}
\left(1+\rho_{\mathcal{T}^*}s\right)^{-1} & S\in\mathcal{T}^*\\
1 & S\notin \mathcal{T}^*,
\end{cases}
\end{equation}
which may be recognized as the MGFs of an exponentially distributed random variable with parameter~$\rho_{\mathcal{T}^*}^{-1} = \mu_{\mathcal{T}^*}/p_{\mathcal{T}^*}$ and the degenerate random variable `0', respectively. 

In order to show that the waiting time of a job of type $S\in\mathcal{S}$ follows the same distribution, 
it is sufficient to note that the sojourn time of a type-$S$ job is upper bounded by the sum of its waiting time and the service time offered by the first server where it starts its service. Note that this service time is independent of the waiting time and is stochastically smaller than an exponentially distributed random variable $X$ with parameter $\mu^- = \min\{\mu_n \colon n =1,\dots,N\}$. Hence, the sojourn time $V_S$ is lower bounded by the waiting time $W_S$, and upper bounded  by $W_S+X$. This results in 
\begin{equation}
\mathcal{L}_{\left(1-\frac{\lambda}{\lambda^*}\right)W_S} (s) \ge \mathcal{L}_{\left(1-\frac{\lambda}{\lambda^*}\right)V_S} (s) \ge \mathcal{L}_{\left(1-\frac{\lambda}{\lambda^*}\right)W_S} (s) \mathcal{L}_{\left(1-\frac{\lambda}{\lambda^*}\right)X} (s).
\end{equation}
Since 
\begin{equation}
\lim_{\lambda\uparrow\lambda^*} \mathcal{L}_{\left(1-\frac{\lambda}{\lambda^*}\right)X} (s) = \lim_{\lambda\uparrow\lambda^*} \frac{\mu^-}{\mu^- + s\left(1-\frac{\lambda}{\lambda^*}\right) } =1,
\end{equation}
it follows that also
\begin{equation}
\lim_{\lambda\uparrow\lambda^*} \mathcal{L}_{\left(1-\frac{\lambda}{\lambda^*}\right)W_S} (s) = \lim_{\lambda\uparrow\lambda^*} \mathcal{L}_{\left(1-\frac{\lambda}{\lambda^*}\right)V_S} (s).
\end{equation}

Since both $V_S$ and $W_S$ are non-negative random variables, convergence of the MGFs implies convergence in distribution to exponentially distributed random variables with parameter~$\mu_{\mathcal{T}^*}/p_{\mathcal{T}^*}$ or to the degenerate random variable `0', depending on whether or not job type $S$ belongs to $\mathcal{T}^*$, respectively~\cite{Feller1971}. For the redundancy c.o.s.\ policy, jobs of the same type enter the system in order of arrival so that the moment generating function of $\tilde{Q}_S$ can be related to the moment generating function of the waiting time $W_S$~\cite[Theorem~2]{Keilson1988}. A similar argument as above shows that both the (scaled) waiting time $W_S$ and the (scaled) sojourn time $V_S$ will converge to the same distribution as under the c.o.s.\ policy.
This concludes the proof of Corollary~\ref{cor:SojournWaiting}.
\end{proof}

\end{document}